\documentclass[11pt]{amsart}
\usepackage[english]{babel}
\usepackage{egothic}
\usepackage[T1]{fontenc}

\usepackage{amsmath}
\usepackage{amsthm}
\usepackage{amssymb}
\usepackage[all]{xy}
\usepackage{mathrsfs}
\usepackage{enumerate}
\usepackage{physics}
\usepackage[notcite, final, notref]{showkeys}
\usepackage{dsfont}
\usepackage{tikz-cd}
\usepackage{tikz}
\usepackage{graphicx}
\usepackage{float}
\usepackage{xcolor}

\newtheorem{theorem}{Theorem}[section]
\newtheorem{problem}[theorem]{Problem}
\newtheorem{lemma}[theorem]{Lemma}
\newtheorem{proposition}[theorem]{Proposition}
\newtheorem{corollary}[theorem]{Corollary}
\newtheorem{definition}[theorem]{Definition}

\theoremstyle{definition}
\newtheorem{remark}[theorem]{Remark}

\usepackage{xcolor}

\title[]{Polyanalytic Gaussian Radial Basis Function Kernel and Itô-Hermite polynomials}

\author[H. De Bie]{Hendrik De Bie}
\address{Clifford Research Group, Department of Electronics and Information Systems, Faculty of Engineering and
	Architecture, Ghent University, Krijgslaan 281, 9000 Gent, Belgium.}
\email{Hendrik.DeBie@UGent.be}

\author[A. De Martino]{Antonino De Martino}
\address{Politecnico di
	Milano\\Dipartimento di Matematica\\Via E. Bonardi, 9\\20133 Milano\\Italy.}
\email{antonino.demartino@polimi.it}

\author[K. Diki]{Kamal Diki}
\address{Clifford Research Group, Department of Electronics and Information Systems, Faculty of Engineering and
	Architecture, Ghent University, Krijgslaan 281, 9000 Gent, Belgium.}
\email{Kamal.Diki@UGent.be}

\begin{document}
\maketitle
\begin{abstract}
We introduce a polyanalytic extension of the Gaussian radial basis function (RBF) kernel by computing the action of the convolution operator on normalized Hermite functions. In particular, using the Zaremba-Bergman formula we derive an explicit closed form for this new reproducing kernel function. We then establish an isomorphism relating the reproducing kernel Hilbert space induced by the polyanalytic Gaussian RBF kernel with the corresponding polyanalytic Fock space. Moreover, we provide a characterization of polyanalytic Gaussian RBF spaces in terms of a Landau-type operator. In addition, we investigate the polyanalytic counterpart of the Weyl operator, which leads to applications involving the Christoffel-Darboux formula for Hermite polynomials and Mehler's kernel. Finally, we discuss the analogue of the Weyl operator in the context of the polyanalytic Gaussian RBF setting.
\end{abstract}


\tableofcontents

\section{Introduction}
\setcounter{equation}{0}
The theory of kernels and reproducing kernel Hilbert spaces (RKHSs) originated in functional analysis with the seminal work of Aronszajn \cite{Aronszajn}. A RKHS can be defined as a Hilbert space of complex-valued functions in which every evaluation functional is bounded. Therefore, by the Riesz representation theorem, for each point in the data set there exists a unique element of the Hilbert space such that its inner product with any function $f$ in the RKHS reproduces the value of $f$ at that point. This unique element is called the \textit{reproducing kernel} of the RKHS. For a general introduction to the theory of reproducing kernel Hilbert spaces, we refer to \cite{PR2016,Saitoh}, and for probability and statistics perspectives, see \cite{Berlinet}. It is important to note that these mathematical methods have turned out to be fundamental in practical applications, particularly in machine learning. For instance, kernel functions are employed to measure similarity between data points and provide the theoretical framework for solving classification and regression problems \cite{STJ}. A large class of learning algorithms rely on RKHSs and are known as kernel methods, among which support vector machines (SVMs) are some of the most popular \cite{SC2008}. Moreover, kernel functions are also used in mathematical physics, for example in the study of coherent states in quantum mechanics \cite{Gaz} (see also \cite{S}). \\

A fundamental example of reproducing kernel Hilbert spaces is the so-called Fock space, or Segal-Bargmann space introduced by Bargmann in \cite{Bargmann1961}. This space has significant applications in quantum mechanics \cite{Hall2013} (see also \cite{Folland, zhu}). It consists of entire functions that are square integrable with respect to the Gaussian measure. Over the years, many extensions of the classical Fock space have been studied, including those in the classical analytic case \cite{ACK, ADDS2025, CLSWY2020, zhu, Z1}, in $q$-calculus context \cite{ACKS, Arik}, in the case of Mittag-Leffler type kernels \cite{ AD2023, AT2025, RRD2018}, in the polyanalytic setting \cite{Asampling, AF, ACDSS1, FH, GS2024,  Vas}, and in various hypercomplex frameworks \cite{ACSS2014,ADD2025, ADSMN, CD2017, DMD2, DGS2019, DG2017, DKS2019}.  \\  

Another crucial example of reproducing kernels is the \textit{Gaussian radial basis function (RBF) kernel}, which is one of the widely used kernels in SVMs. It also plays an important role in neural networks \cite{HKK1990}. A detailed description of reproducing kernel Hilbert spaces of holomorphic functions associated with the Gaussian RBF kernel in the setting of several complex variables has been given in \cite{SC2008, SD2006}. The idea of connecting the Gaussian RBF kernel and Fock space theory was explored for the first time in \cite{ACDSS}. This work showed how constructions from the analytic Fock space setting could be translated into the analytic Gaussian RBF setting, leading to the development of counterparts of well-known operators in quantum mechanics. For example, analogues of the creation, annihilation, and Weyl operators, as well as Fourier and Segal-Bargmann transforms, were obtained in \cite{ACDSS}. Later, extensions of these results in quaternionic analysis and several complex variables were presented in \cite{DMD2024}. \\

The main objective of this paper is to construct a counterpart of the Gaussian RBF kernel and its associated RKHS within the framework of polyanalytic function theory. For an introduction to this theory in complex analysis, we refer to \cite{AF, Balk1991, VasBook}. In this context, our construction relies on two main steps. First, we compute the convolution of modulated Hermite functions in Lemmas \ref{ACHA} and \ref{cor1}. Then, we apply the Zaremba-Bergman formula (see the review \cite{S}) to obtain the polyanalytic kernel, see Propositions \ref{An4} and \ref{truerep}. These computations yield the new polyanlytic Gaussian RBF kernel of order $N\geq 1$, given by 

$$
\displaystyle K_{RBF,N}(z,w):=e^{-\frac{1}{4}(z-\overline{w})^2} L^{1}_{N-1}\left(\frac{|z-w|^2}{2}\right), \quad \text{ for all } z,w\in \mathbb{C},
$$
where $L^{1}_{N-1}$ denotes the generalized Laguerre polynomial of order $N-1$. This polyanalytic kernel extends the standard (analytic) Gaussian RBF kernel denoted by $K_{RBF}$ (\cite{SC2008,SD2006}), which corresponds to $N=1$. In particular, it can be expressed as 
$$
\displaystyle K_{RBF,N}(z,w)=K_{RBF}(z,w) L^{1}_{N-1}\left(\frac{|z-w|^2}{2}\right), \quad  \text{ for all } z,w\in \mathbb{C}.$$ 
\newpage Although we do not consider the multivariable case in this paper, we observe that the polyanalytic kernel on $\mathbb R^d$ can be written in terms of the Euclidean norm $||\cdot||_2$ as

$$\displaystyle K_{RBF,N}(x,y)=e^{-\frac{1}{4}||x-y||^2_2} L^{1}_{N-1}\left(\frac{||x-y||^2_2}{2}\right) =K_{RBF}(x,y) L^{1}_{N-1}\left(\frac{||x-y||^2_2}{2}\right), \quad  \text{ for all } x,y\in \mathbb{R}^d.
$$ 

The following flowchart summarizes the main ideas involved in the construction of the polyanalytic Gaussian RBF kernel:
\begin{figure}[H]
   \centering
    \includegraphics[width=1.1\linewidth]{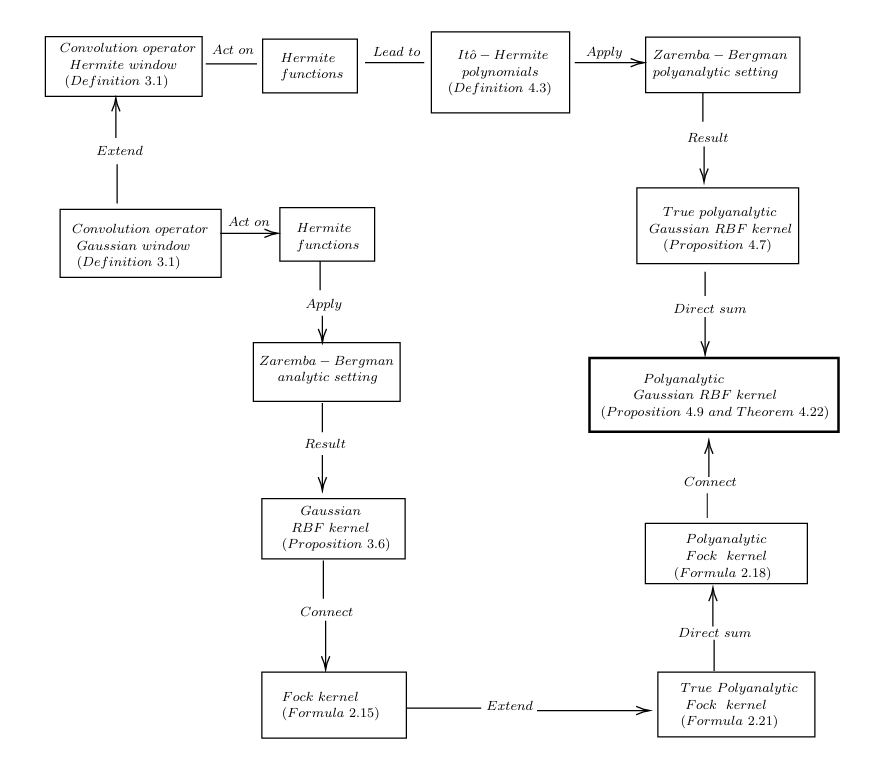}
    \label{figure}
	\caption{Flowchart of proposed approaches for constructing the polyanalytic Gaussian RBF kernels.}
\end{figure}

The idea of computing the Hermite convolution product arose while investigating superoscillations via the Short-time Fourier transform (STFT) in \cite{ADDS2024}. It turns out that the resulting convolution formula is closely related to the well-known Itô-Hermite polynomials, originally introduced in the context of stochastic processes in \cite{I} (see also \cite{IS}). We refer the reader to the book \cite{Hida1980} for further explanations on these polynomials and their applications in complex white noise analysis. Moreover, it is worth noting that this convolution product of modulated Hermite functions is naturally related to the so callled \textit{special Hermite functions}, which is a special family of eigenfunctions considered in \cite{Strichartz1989}, see also the book \cite{Thangavelu1993}. \\ \\

Polyanalytic Fock spaces can also be characterized as eigenspaces of a certain Laudau operator, see \cite[Section 5.1]{AF}. This characterization can be observed by comparing polyanalytic Fock kernels with the reproducing kernels of the Landau eigenspaces computed in \cite{Mouayn}. Moreover, the Itô-Hermite polynomials naturally appear in this context as eigenfunctions of a specific Laundau operator (see \cite{IC}). Inspired by this perspective, we introduce a Landau-type operator within the framework of polyanalytic Gaussian RBF spaces, which can also be related to a Schrödinger type operator with a magnetic field.  \\ \\ We also investigate the counterpart of Weyl operators within polyanalytic function theory. These operators were recently studied on polyanalytic Fock spaces in connection with quantum harmonic analysis \cite{FH}. Our approach introduces these operators as compositions of time-frequency shifts with the true polyanalytic Bargmann transform and its inverse. Based on these tools, we introduce the analogues of the Weyl operators in the polyanalytic RBF setting.\\ \\

The structure and main results of the paper are organized as follows: Section 2 reviews classical notions on the theory of Fock spaces in both the analytic and polyanalytic settings. It also recalls preliminary results on analytic Gaussian RBF kernels and a convolution formula of Hermite functions that will be used in the sequel. Section 3 introduces the convolution operator with Gaussian and Hermite windows. In particular, we calculate the action of the convolution operator with Gaussian window on normalized Hermite functions. Then, we observe that the Zaremba-Bergman formula associated with this specific operator yields the classical Gaussian (analytic) RBF kernel. Section 4 is devoted to the study of the polyanalytic Gaussian RBF kernels using the convolution operator with a Hermite window and its associated Zaremba-Bergman formula. The polyanalytic Gaussian RBF reproducing kernel Hilbert spaces are presented and connected to their Fock polyanalytic analogues. Section 5 discusses an operator approach allowing us to introduce these spaces as eigenspaces of a suitable Landau-type operator. In section 6, we study the counterpart of the Weyl operator in the polyanalytic context, by considering the composition of polyanalytic Bargmann transforms with the time-frequency shifts. As a result of these calculations, we give two applications involving the Christoffel-Darboux formula and Mehler's kernel. An analogue of the Weyl operator on the polyanalytic Gaussian RBF space is also presented. Finally, the paper's appendix includes some straightforward extensions in the case of parametric spaces, a generating function for Itô-Hermite polynomials, and an alternative proof for the convolution product of Hermite functions.

%\lipsum[1-4]

%\includegraphics[width=0.5\textwidth, center]{F.jpg}
%\includegraphics[height=6cm, width=5cm]{Picture.png}

%\begin{figure}[h]
%	\centering
%	\includegraphics[width=0.8\textwidth]{Poly.jpg}
%	\caption{Polyanalytic Gaussian Radial Basis Function Kernel}
%\end{figure}

\section{Preliminaries}
 \setcounter{equation}{0}
We begin by reviewing basic definitions and properties of \textit{analytic and polyanalytic Fock spaces, analytic and polyanalytic Segal-Bargmann transforms, Gaussian radial basis function (RBF) kernel, the Weyl operator, and complex Itô-Hermite polynomials.} 
\subsection{Analytic case}
First, let us introduce the classical (parametric) analytic Fock space using its geometric description \cite{Folland, zhu}:
\begin{definition}\label{Fockdef}
Let $\alpha>0$. An entire function $f: \mathbb{C} \to \mathbb{C}$ belongs to the \textit{(parametric) Fock space}, denoted by $ \mathcal{F}_\alpha(\mathbb{C})$, if
$$ \| f \|_{\mathcal{F}_\alpha (\mathbb{C})}^2= \frac{\alpha}{\pi} \int_{\mathbb{C}} | f(z)|^2 e^{-\alpha |z|^2} d \lambda(z) < \infty,$$
where $ d \lambda(z)=dx dy$ is the classical Lebesgue measure for $z=x+iy$.
\end{definition}
The space $\mathcal{F}_\alpha(\mathbb{C})$ is equipped with the inner product

\begin{equation}
\langle f,g \rangle_{\mathcal{F}_\alpha(\mathbb{C})}=\frac{\alpha}{\pi}\int_\mathbb{C}\overline{g(z)}f(z)e^{-\alpha|z|^2} d \lambda(z),
\end{equation}
for all $f,g \in  \mathcal{F}_\alpha(\mathbb{C})$. The Fock space is a reproducing kernel Hilbert space (RKHS), with reproducing kernel 
\begin{equation}\label{Kf}
K_\alpha(z,w)=K_{\alpha,w}(z):=e^{\alpha z \overline{w}}, \quad \text{ for all }  z,w \in \mathbb{C}.
\end{equation}

Furthermore, the reproducing kernel property ensures that for every $f\in\mathcal{F}_\alpha(\mathbb{C})$ and $w\in \mathbb{C} $, the value of $f(w)$ can be expressed as
\begin{equation}
f(w)=\displaystyle  \frac{\alpha}{\pi}\int_{\mathbb{C}} \overline{K_\alpha(z,w)} f(z) e^{-\alpha |z|^2}d \lambda(z)=\langle f, K_{\alpha, w}\rangle_{ \mathcal{F}_\alpha(\mathbb{C})}.
\end{equation}

For a fixed $w\in\mathbb{C}$, the \textit{normalized Fock kernel} is defined as
\begin{equation}
k_\alpha(z,w)=k_{\alpha,w}(z):=\dfrac{K_\alpha(z,w)}{\sqrt{K_\alpha(w,w)}}=e^{\alpha(z\overline{w}-\frac{|w|^2}{2})}, \quad \text{ for all } z\in\mathbb{C}.
\end{equation}

The family of normalized Fock kernels $(k_{\alpha, w})_{w\in\mathbb{C}}$ is referred to as \textit{coherent states} in the terminology of quantum mechanics, as they form eigenstates of the annihilation operator. See \cite{Gaz, GHS2019} for further details and explanations. 
\begin{remark}
In most results in this paper, we choose $\alpha=1$ to simplify our calculations. In particular, for $\alpha=1$, we use the simplified notation for the Fock space $ \mathcal{F}(\mathbb{C}):=\mathcal{F}_1(\mathbb{C})$ and for the Fock kernel $K(z,w):=K_1(z,w)=e^{z\overline{w}}$. When the parameter $\alpha$ is relevant, it will be explicitly stated in the corresponding result. 
\end{remark}
We now recall the Hermite polynomials, which will play an essential role in the sequel.
 \begin{definition}\label{HermitePol}
The Hermite polynomials are defined via their generating function
\begin{equation}\label{HermiteGF}
\displaystyle  e^{2xz-z^2}=\sum_{n=0}^{\infty} \frac{H_n(x)}{n!}z^n,
\end{equation}
for every $z\in\mathbb{C}$ and $x\in \mathbb{R}$.  The Hermite functions are defined by multiplying the Hermite polynomials $H_n$ with the Gaussian function, i.e.,
\begin{equation}
h_n(x)=e^{-\frac{x^2}{2}}H_n(x), \quad n=0,1,\cdots.
\end{equation}
\end{definition}

\begin{remark}
The Hermite polynomials appear as the Taylor coefficients of the holomorphic function $e^{2xz - z^2}$. \end{remark}
An orthonormal basis of the Hilbert space $L^2(\mathbb{R})$ is given by the normalized Hermite functions, defined as follows (see \cite{Lebedev1972, Szego})
\begin{equation}
\psi_n(x):=\frac{h_n(x)}{||h_n||_{L^2}}=\frac{h_n(x)}{\sqrt{n!2^n \sqrt{\pi}}},\quad  n=0,1,2, \cdots
\end{equation}
The Segal-Bargmann transform, introduced in \cite{Bargmann1961}, is defined as follows:
\begin{definition}\label{SBkernel}
The Segal-Bargmann kernel is defined by the following expression
$$\displaystyle   A(z,x)=A_z(x) := \sum_{n=0}^{\infty}\psi_n(x)\frac{z^n}{\sqrt{n!}}=\pi^{-1/4} e^{-\frac{1}{2}(x^2+z^2)+\sqrt{2} zx},$$
for every $z\in \mathbb{C}$ and $x\in\mathbb{R}$. The Segal-Bargmann transform of a function $\psi\in L^2(\mathbb{R})$ is given by

\begin{align*}
\displaystyle (\mathcal{B}\psi)(z)&:=\langle \psi, A_{\overline{z}} \rangle_{L^2(\mathbb{R})} \\
&=\pi^{-1/4} \int_{\mathbb{R}} e^{-\frac{1}{2}(x^2+z^2)+\sqrt{2}zx}\psi(x)dx,
\end{align*}
where $z\in \mathbb{C}$, and $A_z$ is the Segal-Bargmann kernel.
\end{definition}

It is worth noting that the Segal-Bargmann transform establishes a connection between the Schrödinger representation and the Bargmann-Fock representation through the following fundamental result (see \cite{Bargmann1961}).

\begin{theorem}
The Segal-Bargmann transform $\mathcal{B}$ defines a unitary operator that maps the Hilbert space $L^2(\mathbb{R})$ onto the Fock space $\mathcal{F}(\mathbb{C})$. Furthermore, $\mathcal{B}$ transforms the orthonormal basis $(\psi_k)_{k\in\mathbb{N}_0}$ onto an orthonormal basis of $\mathcal{F}(\mathbb{C})$, given by
\begin{equation}
(\mathcal{B}\psi_k)(z)=\frac{z^k}{\sqrt{k!}}:=e_k(z), \quad k\in \mathbb N_0.
\end{equation}
\end{theorem}

 An important operator on the Fock space is the \textit{Weyl operator}, which is defined by multipying the normalized Fock kernel with the complex translation operator. This operator connects the Schrödinger and Bargmann-Fock representations of the Heiseinberg group \cite{GK2001}, and is defined as follows:
\begin{definition}[Weyl operator]
For $a\in\mathbb{C}$, the Weyl operator $\mathcal{W}_{a}:\mathcal{F}(\mathbb{C})\longrightarrow \mathcal{F}(\mathbb{C}),$ is defined by
\begin{equation}
\label{cv}
\mathcal{W}_a f(z):=k_a(z)f(z-a), \quad f\in \mathcal{F}{(\mathbb{C})},  z, a\in \mathbb{C}.
\end{equation}
\end{definition}

The Weyl operator $\mathcal{W}_a$ is unitary on the Fock space $\mathcal{F}(\mathbb{C})$, with inverse and adjoint given by

\begin{equation}
(\mathcal{W}_a)^{-1}=(\mathcal{W}_a)^{*}=\mathcal{W}_{-a}.
\end{equation}
Moreover, the Weyl operators satisfy the following semi-group-like relation
\begin{equation}\label{Weylsg}
\mathcal{W}_a\mathcal{W}_b=\exp(- i  \Im(a\bar{b}))\mathcal{W}_{a+b}, \quad a,b\in\mathbb{C}.
 \end{equation}
\begin{remark}
If $a,b\in\mathbb{R}$, the Weyl operators satisfy the classical semi-group property
\begin{equation}\label{Weylsg1}
\mathcal{W}_a\mathcal{W}_b=\mathcal{W}_{a+b}.
 \end{equation}
\end{remark} 

We conclude this section by recalling the (analytic) Gaussian radial basis function (RBF) kernel in a complex variable, along with the associated reproducing kernel Hilbert spaces, as introduced in \cite{SD2006} (see also \cite{SC2008}). 

\begin{definition}[Gaussian RBF kernel] Let $\gamma >0$. The function 
\begin{equation}
\label{rbfkernel}
K^\gamma_{RBF}(z,w)=\exp\left(-\dfrac{(z-\overline{w})^2}{\gamma^2}\right), \quad \text{ for all } z,w\in\mathbb{C},
\end{equation}

is called the Gaussian RBF kernel, or simply the RBF kernel, with width parameter $\displaystyle \frac{1}{\gamma}$. The special case $ \gamma = 2 $ will be denoted by $K_{\mathrm{RBF}}(z, w)$.

\end{definition}
\begin{remark}
For real variables $x$ and $y$ in $\mathbb{R}$, the kernel reduces to the standard real-valued Gaussian RBF kernel
$$K^\gamma_{RBF}(x,y)=\exp\left(-\dfrac{(x-y)^2}{\gamma^2}\right),$$
which is commonly used in kernel methods such as support vector machines (SVMs).
\end{remark}
The reproducing kernel Hilbert spaces (RKHSs) associated with the complex Gaussian RBF kernels $(K_\gamma)_{\gamma>0}$ can be introduced via the the following geometric description: 

\begin{definition}[Gaussian RBF RKHS] 
\label{RBF}
	Let $\gamma >0$.  An entire function $f:\mathbb{C}\longrightarrow \mathbb{C}$ belongs to the \textit{Gaussian RBF space} (or \textit{RBF space}, for short), denoted by $\mathcal{H}_{\gamma}^{RBF}(\mathbb{C})$ (or simply $\mathcal{H}_\gamma)$, if
\begin{equation}
\displaystyle ||f||_{\mathcal{H}_\gamma}^2:=\dfrac{2}{\pi\gamma^2}\int_{\mathbb{C}}|f(z)|^2\exp\left(\frac{(z-\overline{z})^2}{\gamma^2}\right) d\lambda(z)<\infty,
\end{equation}
where $d\lambda(z)=dxdy$ is the Lebesgue measure with respect to the complex variable $z=x+iy$.
\end{definition}

Further results connecting the Gaussian RBF space with the classical Bargmann-Fock space can be found in \cite{ACDSS}. For instance, an entire function $f:\mathbb{C}\longrightarrow \mathbb{C}$ belongs to the Gaussian RBF space $\mathcal{H}_\gamma$ if and only if there exists a unique function $g\in\mathcal{F}_{\frac{2}{\gamma^2}}$such that $$f(z)=\exp(-\frac{z^2}{\gamma^2})g(z), \quad \text{ for all } z\in\mathbb{C}.$$Therefore, there is an isometric isomorphism connecting the Gaussian RBF and Fock spaces, given by the following multiplication operator

\begin{equation}
\label{RBFfock}
\mathcal{M}_\gamma:\mathcal{H}_\gamma\longrightarrow \mathcal{F}_{\frac{2}{\gamma^2}}, \quad \mathcal{M}_\gamma[f](z):=\exp(\frac{z^2}{\gamma^2})f(z), \quad \text{ for all } f\in \mathcal{H}_\gamma, z\in \mathbb{C}.
\end{equation}

We refer the reader to \cite{DMD2024} for extensions of such connections to the case of several complex variables and to the quaternionic setting.

\subsection{Polyanalytic case}
We briefly review basic concepts and properties of polyanalytic functions as needed for this paper. For a general introduction to this function theory, we refer the reader to \cite{AF, Balk1991,VasBook}.

\begin{definition}

Let $\Omega$ be a domain in $\mathbb C$ and let $N\in \mathbb N$. A complex valued function $f:\Omega\subset \mathbb{C}\longrightarrow \mathbb{C}$ of class $\mathcal{C}^N$ is said to be \textit{polyanalytic of order} $N$  if it belongs to the kernel of the $N$-th power of the classical Cauchy-Riemann operator, that is, $$\displaystyle \frac{\partial^N}{\partial \overline{z}^N}f(z)=0, \quad \text{ for all }  z\in\Omega.$$
The space of all such functions defined on $\Omega$ is denoted by $H_N(\Omega)$. 
\end{definition}
\begin{remark}
 A function that is polyanalytic of order $N$ on the entire complex plane $\mathbb{C}$ is called \textit{poly-entire} of order $N$. If we consider $N=1$ in the above definition we get the classical definition of holomorphic functions.
\end{remark}
A key property of this function theory is that any polyanalytic function of order $N$ can be decomposed in terms of $N$ analytic functions. Specifically, any function $f\in H_N(\Omega)$ can be written in the form (see \cite[p. 9, Formula (1.1)]{Balk1991})
\begin{equation}\label{Polydeccomplex}
f(z)=\displaystyle \sum_{k=0}^{N-1}\overline{z}^kf_k(z),
\end{equation}
where each $f_k$ is an analytic function in $\Omega$. This decomposition has been used recently in \cite{CDKSS} to extend classical duality theorems to the polyanalytic setting.
\subsubsection{Polyanalytic and true polyanalytic Fock spaces} 

The polyanalytic Fock space of order $N$, discussed in \cite[pp.169-170]{Balk1991}, extends the classical analytic Fock space (see Definition \ref{Fockdef}) by considering poly-entire functions of order $n$ that are square-integrable with respect to the normalized Gaussian measure. Specifically, we have
\begin{definition}
\label{paraFock}
Let $N\in \mathbb N$ and $\alpha>0$. The (parametric) polyanalytic Fock space of order $N$, denoted $\mathcal{F}_{\alpha,N}(\mathbb{C})$, is defined by
$$\mathcal{F}_{\alpha,N}(\mathbb{C}):=\left\lbrace f\in H_N(\mathbb{C})\mid  \frac{\alpha}{\pi}\int_{\mathbb{C}}|f(z)|^2e^{-\alpha|z|^2}d\lambda(z)<\infty \right\rbrace,$$
where $d\lambda(z)$ is the classical Lebesgue measure on $\mathbb C$. In particular, when $\alpha=1$, we use the simplified notation for the polyanalytic Fock space $ \mathcal{F}_N(\mathbb{C}):=\mathcal{F}_{1,N}(\mathbb{C})$.
\end{definition}

The reproducing kernel associated with the reproducing kernel Hilbert space $\mathcal{F}_{\alpha,N}(\mathbb{C})$ is given explicitly by (see \cite{Balk1991})
\begin{equation}\label{Kn}
K_N^\alpha(z,w)=e^{\alpha z\overline{w}}\displaystyle \sum_{k=0}^{N-1}\frac{(-1)^k}{k!}{N \choose k+1}\alpha^{k}|z-w|^{2k}, \quad \text{ for all } z,w\in\mathbb{C}.
\end{equation}

Furthermore, the polyanalytic kernel \eqref{Kn} can be expressed in terms of generalized Laguerre polynomials (see \cite{AF})
 
 \begin{equation}
 	\label{kernelP}
K_N^\alpha(z,w)=e^{\alpha z\overline{w}} L^{1}_{N-1}(\alpha|z-w|^2), \quad \text{for all } z, w \in \mathbb{C},
\end{equation}
where $L^{\beta}_{k}(x)$ denotes the generalized Laguerre polynomial of order $k$ and parameter $\beta$, defined by

\begin{equation}\label{Laguerre}
L^\beta_k(x):=\sum_{j=0}^{k}(-1)^j {k+\beta\choose k-j  } \frac{x^j}{j!}.
\end{equation}
In particular, for $\alpha=1$, we use the following notation for the polyanalytic Fock kernel
$$K_N(z,w)=e^{z\overline{w}} L^{1}_{N-1}(|z-w|^2), \quad \text{for all } z, w \in \mathbb{C}.$$
We now recall the so-called \textit{true polyanalytic Fock spaces}, introduced by Vasilevski in \cite{Vas} (see also \cite{Asampling, AF} for further discussion). 
\begin{definition}
\label{FT}
Let $\ell \in \mathbb{N}$ and $\alpha>0$.
A function $ f: \mathbb{C} \to \mathbb{C}$ belongs to the true (parametric) polyanalytic Fock space $ \mathcal{F}_{T,\alpha}^{(\ell)} (\mathbb{C})$ if and only if
$$\frac{\alpha}{\pi} \int_{\mathbb{C}}|f(z)|^2 e^{-\alpha |z|^2} \, d \lambda(z) < \infty,$$
and there exists an entire function $g$ such that
\begin{equation}\label{Bpoly2}
f(z)= \frac{1}{\alpha^{\ell-1}\sqrt{(\ell-1)!}} e^{\alpha |z|^2} \frac{\partial ^{\ell-1}}{\partial z^{\ell-1} }(e^{-\alpha |z|^2} g(z)), \quad \text{ for all } z\in \mathbb{C}.
\end{equation}
When $\alpha=1$, we use the simplified notation for the true polyanalytic Fock space $ \mathcal{F}_T^{(\ell)}(\mathbb{C}):=\mathcal{F}_{T,1}^{(\ell)}(\mathbb{C})$ and denote by $\mathsf{K}_{\ell,T}(z,w)$ the corresponding kernel.
\end{definition}
Each space $\mathcal{F}_{T,\alpha}^{(\ell)}(\mathbb C) $ is a reproducing kernel Hilbert space with reproducing kernel 
\begin{equation}
	\label{tpoly}
\mathsf{K}_{\ell,T}^\alpha(z,w)=e^{\alpha z\overline{w}} L^{0}_{\ell-1}(\alpha |z-w|^2), \quad \text{for all }z, w \in \mathbb{C}, \quad \ell=1,\cdots, N.
\end{equation}
The particular case $\alpha=1$ is denoted by $\mathsf{K}_{\ell,T}(z,w)$. The relation between the polyanalytic Fock kernel $K_N^\alpha$ and the true polyanalytic Fock kernels $(\mathsf{K}_{\ell, T}^\alpha)_{\ell=1,\cdots, N}$ is given by the following formula 

\begin{equation}
\displaystyle K_N^\alpha(z,w)=\sum_{\ell=1}^{N}\mathsf{K}_{\ell,T}^\alpha (z,w), \quad \text{ for all } z,w \in \mathbb{C}.
\end{equation}

As a result, the space $\mathcal{F}_{\alpha,N}(\mathbb{C})$ can be obtained as an orthogonal direct sum of the true polyanalytic Fock spaces $\mathcal{F}_{T, \alpha}^{(\ell)}(\mathbb{C})$, $\ell=1,...,N$, i.e.,
\label{R3}

\begin{equation}\label{VasDecom}
 \mathcal{F}_{\alpha, N}(\mathbb{C})= \bigoplus_{\ell=1}^{N} \mathcal{F}_{T,\alpha}^{(\ell)}(\mathbb{C}).
\end{equation}
The decomposition \eqref{VasDecom} of the polyanalytic Fock space $ \mathcal{F}_{\alpha, N}(\mathbb{C})$ in terms of the true polynalytic Fock spaces $\mathcal{F}_{T,\alpha}^{(\ell)}(\mathbb{C})$, $\ell=1,...,n,$ was obtained by Vasilevski in \cite{Vas}(see also \cite{AF}). Furthermore, we also have 
\begin{equation}
\label{l2fock}
L^2\left(\mathbb{C},e^{-\alpha |z|^2}\right)= \bigoplus_{\ell=1}^{\infty } \mathcal{F}_{T,\alpha}^{(\ell)}(\mathbb{C}).
\end{equation}

The true polyanalytic Segal-Bargmann transform ($\alpha=1$) was introduced by Vasilevski in \cite{Vas} as a unitary operator $B_\ell$ mapping $L^2(\mathbb R)$ onto the true polyanalytic Fock space  $\mathcal{F}_{T}^{(\ell)}(\mathbb{C})$, defined by

\begin{equation}
\label{truepoly}
(B_\ell \psi)(z):= \int_\mathbb{R}H_{\ell-1}\left(\frac{z+\overline{z}}{\sqrt{2}}-x\right) A_z(x) \psi(x) dx,\quad \text{ for all } z\in \mathbb{C}, \psi\in L^2(\mathbb R),
\end{equation}
where $H_{\ell-1}$ denotes the Hermite polynomial of order $\ell-1$ and $A_z$ is the Segal-Bargmann kernel introduced in Definition \ref{SBkernel}. 
\begin{remark}

The inverse of the true polyanalytic Bargmann transform is given by
\begin{equation}
\label{inverse}
(B_\ell^{-1} \psi)(x):= \int_\mathbb{C}H_{\ell-1}\left(\frac{z+\overline{z}}{\sqrt{2}}-x\right) \overline{A_z(x)} g(z) \frac{d \lambda(z)}{\pi},\quad \text{ for all } x \in \mathbb{R}, g \in \mathcal{F}^{(\ell)}_T(\mathbb C),
\end{equation}
\end{remark}
It is important to note that the true polyanalytic Bargmann transform can be introduced via the representation \eqref{Bpoly2} by choosing the entire function $g$ to be the classical Segal-Bargmann transform of the input function $\psi$, i.e $g=\mathcal{B}\psi)$ (see \cite{AF}).
Finally, one may define the \textit{full polyanalytic Bargmann transform}, which maps each vector $\mathbf \varphi=(\varphi_1, \cdots, \varphi_N)\in L^2(\mathbb R, \mathbb C^N)$ to the following polyanalytic function of order $N$ in $\mathcal{F}_N(\mathbb C)$:
\begin{equation}
\label{fullpoly}
\displaystyle (\mathbf{B}_N\varphi)(z):= \sum_{\ell=1}^{N} (B_\ell \varphi_\ell)(z), \quad \text{ for all } z\in \mathbb C.
\end{equation}
\begin{remark}
In the special case $N= \ell = 1 $, the true polyanalytic and polyanalytic Fock space reduce to the classical analytic Fock space, recovering the Segal--Bargmann transform introduced earlier.
$$\mathcal{F}_{1, \alpha}(\mathbb C)=\mathcal{F}_{T,\alpha}^{(1)}(\mathbb C)=\mathcal{F}_\alpha(\mathbb C),$$
and $$\mathbf{B}_1=B_1=\mathcal B.$$
\end{remark}

In \cite{AF}, it has been proved that for $\ell \in \mathbb{N}$, the true polyanalytic Fock space coincides with the so-called \emph{Landau level}, defined as

\begin{equation}
	\label{mouy}
	\mathcal{A}_{\ell,\alpha}(\mathbb{C}):= \left\{ f \in L^2 \left(\mathbb{C}, e^{- \alpha |z|^2}\right) \,: \quad \tilde{\Delta} f=\alpha (\ell-1)f\right\}, \quad \alpha>0
\end{equation}

where 
\begin{equation}
\label{magnetic}
\widetilde{\Delta}:=-\frac{d^2}{dz d \bar{z}}+ \alpha \bar{z} \frac{d}{d \bar{z}},
\end{equation}

 is the so-called magnetic Laplacian, i.e.,
\begin{equation}
\label{mouy1}
\mathcal{F}_{T,\alpha}^{(\ell)}(\mathbb{C})=\mathcal{A}_{\ell, \alpha}(\mathbb{C}).
\end{equation}

\subsubsection{Complex Itô-Hermite polynomials} 
Important examples of polyanalytic functions in the complex setting are the \textit{complex Hermite polynomials}, introduced by Itô in 1952 for applications in stochastic processes (see \cite{I}). These polynomials have also been studied in the case of two complex variables, particularly in relation to combinatorics and integral operators (see \cite{GHS2019}).
\begin{definition}
Let $n,m \in \mathbb{N}_0$ and $\alpha>0$. The complex Hermite polynomials are defined by
\begin{equation}
	\label{hermcom}
	H_{n , m}^\alpha(z,\overline{z})=\alpha^\ell \sum_{j=0}^{min(n, m)} (-1)^j j! \binom{n}{j} \binom{m}{j} \alpha^{n-j} z^{m-j} \overline{z}^{n-j},
\end{equation}
for all $z \in \mathbb{C}$. 
When $ \alpha = 1 $, we simplify the notation for the complex Hermite polynomials by writing $ H_{n,m}(z, \bar{z}) = H_{n,m}^1(z, \bar{z})$.

\end{definition}

\begin{remark}
The complex Hermite polynomials are polyanalytic of order $n+1$.
\end{remark}

\begin{remark}
The complex Hermite polynomials $H_{n , m}^\alpha(z,\overline{z})$ can be extended to holomorphic functions of two complex variables by considering
\begin{equation}
	\label{compexhermite}
	H_{n , m}^\alpha(z,w)= \alpha^\ell \sum_{j=0}^{min(n, m)} (-1)^j j! \binom{n}{j} \binom{m}{j} \alpha^{n-j}z^{m-j} w^{n-j},
\end{equation}
for all $ z,w\in \mathbb{C}$.
\end{remark}
The complex Hermite polynomials can be obtained by applying powers of a specific differential operator to monomials, as shown in \cite[Proposition 2.1]{Intissars} (see also \cite[Appendix A]{DMD2} for the parametric version). Specifically, we have:
\begin{equation}
	\label{oper}
	\left(-\frac{d}{dz}+\alpha \bar{z}\right)^n (\alpha z)^m=H_{n,m}^\alpha(z, \bar{z}), \qquad \alpha>0.
\end{equation}
The complex Hermite polynomials satisfy the following orthogonality relation:

\begin{equation}
	\label{ort}
	\int_{\mathbb{C}} H_{n,m}^\alpha(z, \bar{z})\overline{H_{n',m'}^\alpha(z, \bar{z})} e^{-\alpha | z|^2} d \lambda(z)= \pi  n! m! \alpha^{n+m-1}\delta_{n,n'} \delta_{m,m'}.
\end{equation}
Other important properties of the complex Hermite polynomials are listed below:
\begin{itemize}
\item The Rodrigues formula:
	\begin{equation}
	\label{rodri}
	H_{n,m}^\alpha(z, \bar{z})=(-1)^{n+m} e^{\alpha |z|^2} \frac{\partial^{n+m}}{\partial z^n \partial \bar{z}^m} \left(e^{-\alpha|z|^2}\right), \quad \alpha>0.
	\end{equation}

	\item The exchange of the indexes of the complex Hermite polynomials leads to
	\begin{equation}
	\label{exx}
	H_{n,m}(z, \bar{z})=\overline{H_{m,n}(z, \bar{z})}.
	\end{equation}
\item For $\alpha >0$ we have
\begin{equation}
	\label{dabr}
	\frac{d}{d\bar{z}} H_{n,m}^\alpha(z, \bar{z})=n  \alpha H_{n-1,m}^\alpha(z, \bar{z}), \quad \alpha>0.
\end{equation}
\item For $w \in \mathbb{C}$, $m \in \mathbb{N}$ and $\alpha >0$ we have
\begin{equation}
	\label{gen}
	\sum_{n=0}^{\infty} \frac{H_{m,n}^\alpha(z, \bar{z})H_{n,m}^\alpha(w, \bar{w})}{n! \alpha^n}= \alpha^m m! L_m^{0}(\alpha|z-w|)e^{\alpha z \bar{w}}, 
\end{equation}
where $L_m^{0}$ are Laguerre  polynomials. Although formula \eqref{gen} is considered common knowledge, as it follows from the Zaremba-Bergman formula, a proof is provided in Appendix \ref{App2} for completeness.

\item Let $u \in \mathbb{R}$ and $z \in \mathbb{C}$, then we have (see \cite[Proposition 3.4]{Ghanmi})
\begin{equation}
\label{genHer}
\sum_{n=0}^{\infty}\frac{z^n}{n!} H_{m,n}\left(\frac{u}{\sqrt{2}}, \frac{u}{\sqrt{2}}\right)= \left(\frac{u}{\sqrt{2}}-z\right)^m e^{\frac{uz}{\sqrt{2}}}.
\end{equation}
\item The Itô-Hermite polynomials can be written in terms of Laguerre polynomials (see \cite{I}):
\begin{equation}
\label{relHL}
 H_{n,m}(z, \bar{z})=(-1)^{\min(n,m)} \min(n,m)! e^{i \theta(n-m)}|z|^{|n-m|}L_{\min(n,m)}^{|n-m|}(|z|^2), \quad z=|z|e^{i \theta}.
\end{equation}

\end{itemize}

Given a time-domain signal $\varphi:\mathbb{R} \to \mathbb{C}$, and fixed $x,\omega\in \mathbb{R}$, we define the \textit{translation} and \textit{modulation} operators by:
\begin{equation}
(\tau_xf)(t)=f(t-x),  \quad M_\omega f(t)=e^{i\omega t}f(t).
\end{equation}
These two operators satisfy the \textit{canonical commutation rule}, given by
\begin{equation}
\label{commu}
\tau_xM_\omega=e^{- i x\omega}M_\omega \tau_x.
\end{equation}
As a consequence of the previous identity, we observe that the translation and modulation operators commute if and only if $ x \omega \in 2\pi \mathbb{Z}$. The operators $\tau_xM_\omega$ and $M_\omega \tau_x$ are referred to as \textit{time-frequency shifts}. 
\\ \\ Finally, we recall a useful formula (see \cite[Proposition 5.7 and Corollary 5.15]{ADDS2024}), which will be crucial for our computations in this paper.
\begin{lemma}\label{ACHA}
Let $k,m\in \mathbb{N}_0$ and $u,x,\lambda\in \mathbb R$. Then the following identity holds
\begin{align*}
\displaystyle (M_xh_k*M_uh_m)(\lambda)&=i^{m+k}e^{-\frac{x^2}{2}-\frac{u^2}{2}}\int_{\mathbb{R}} e^{-t^2+t(u+x+i\lambda)}H_k(x-t)H_m(u-t)dt\\
&= \sqrt{\pi}i^{m-k} 2^{\frac{k+m}{2}} e^{- \frac{\lambda^2}{4}+\frac{i \lambda (x+u)}{2}}   e^{-\frac{(x-u)^2}{4}} H_{k,m}\left(\frac{u-x+i\lambda}{\sqrt{2}},\frac{u-x-i\lambda}{\sqrt{2}}\right).
\end{align*}
\end{lemma}
\begin{remark}
By comparing the integral formula in Lemma \ref{ACHA} with the definition of the special Hermite functions given in \cite{Strichartz1989, Thangavelu1993}, one can see directly that these special Hermite functions can be expressed in terms of the complex Itô-Hermite polynomials.
\end{remark}
The special case of $u=x=0$ in Lemma \ref{ACHA} yields a representation of the complex Itô-Hermite polynomials as a convolution of the classical Hermite functions. Specifically, see \cite{ADDS2024}, we obtain
\begin{lemma}[Convolution of Hermite functions]
	\label{cor1}
	Let $m,n \in \mathbb{N}_0$. The convolution of two Hermite functions is given by
	\begin{equation}
		\label{conv}
		(h_m*h_n)(\lambda)= \sqrt{\pi} 2^{\frac{m+n}{2}} e^{- \frac{\lambda^2}{4}} H_{m,n} \left( \frac{\lambda}{\sqrt{2}}, \frac{\lambda}{\sqrt{2}}\right), \quad \lambda \in \mathbb{R}.
	\end{equation}
	Equivalently, in terms of normalized Hermite functions $\psi_k$,we have 
	\begin{equation}
		\label{norm}
		(\psi_m*\psi_n)(\lambda)=\frac{1}{\sqrt{m!n!}} e^{- \frac{\lambda^2}{4}} H_{m,n} \left( \frac{\lambda}{\sqrt{2}}, \frac{\lambda}{\sqrt{2}}\right), \quad \lambda \in \mathbb{R}.
	\end{equation}
\end{lemma}

This identity shows that the convolution of Hermite functions corresponds to the evaluation of the complex Itô-Hermite polynomials on the diagonal. 
\begin{remark}
A proof of the convolution formula \eqref{conv} can be found in \cite[Corollary 5.8]{ADDS2024}. In the appendix of this paper, we provide an alternative proof based on an algebraic approach.
\end{remark}
 \section{Convolution operator and Gaussian RBF kernel: analytic setting}
 \setcounter{equation}{0}
 In this section, we show that the Gaussian RBF kernel can be computed using the \textit{Zaremba-Bergman formula}, see \cite{S}, which is associated with the action of a convolution operator (with a Gaussian window) on the normalized Hermite functions. 
 \begin{definition}[Convolution Operator]\label{ConvolutionOperator}
Let $g\in L^1(\mathbb{R})$ be fixed. The associated convolution operator is defined by
$$\displaystyle (T_gf)(u)=\int_\mathbb{R}g(u-x)f(x)dx=(g*f)(u).$$
When the window function is $g=\psi_m$ for $m\in\mathbb{N}_0$, we denote the operator by $T_{\psi_m}=T_m$, that is, 
$$\displaystyle  (T_{m}f)(u):=(T_{\psi_m}f)(u)=\int_\mathbb{R}\psi_m(u-x)f(x)dx=(\psi_m*f)(u).$$
\end{definition}
In what follows, we focus on the case $g(u)=\psi_0(u)$ and consider the corresponding convolution operator $T_0$, which is denoted simply by $T$.
\begin{remark}
The convolution operator $T_g$ is a bounded linear operator from $L^2(\mathbb{R})$ into itself and it satisfies the norm inequality $$||T_g||\leq ||g||_{L^1},$$ where $||T_g||$ denotes the operator norm.

\end{remark}

We begin by computing the action of the convolution operator $T$ on the Segal-Bargmann kernel $A(z,x)$, see Definition \ref{SBkernel}.
\begin{lemma}\label{TBarg}
For all $z\in \mathbb{C}$ and $u\in\mathbb{R}$, we have
\begin{equation}
(TA_{z})(u)=e^{\frac{zu}{\sqrt{2}}-\frac{u^2}{4}}=k_{\frac{u}{\sqrt{2}}}(z).
\end{equation}
\end{lemma}
\begin{proof}
By definition of the convolution operator $T$, we have
\begin{align*}
\displaystyle (TA_{z})(u)&=\int_{\mathbb{R}}\psi_0(u-x)A_{z}(x)dx\\
&=\frac{e^{-\frac{z^2}{2}}}{\sqrt{\pi}}\int_{\mathbb{R}}e^{-\frac{(u-x)^2}{2}}e^{-\frac{x^2}{2}+\sqrt{2}zx}dx\\
&=\frac{1}{\sqrt{\pi}}e^{-\frac{1}{2}(z^2+u^2)}\int_{\mathbb{R}}e^{-x^2+x(u+\sqrt{2}z)}dx.
\end{align*}
Using the Gaussian integral formula, we obtain

\begin{align*}
\displaystyle (TA_{z})(u)&=e^{-\frac{1}{2}(z^2+u^2)}e^{\frac{(u+\sqrt{2}z)^2}{4}}\\
&=e^{\frac{zu}{\sqrt{2}}-\frac{u^2}{4}}\\
&=k_{1,\frac{u}{\sqrt{2}}}(z).
\end{align*}

\end{proof}

\begin{proposition}
\label{kamal1}
Let $g(u)=\psi_0(u)$, where $\psi_0$ is the ground state function. Then, for every $n\in \mathbb{N}_0$ and $u\in \mathbb R,$ the convolution of $\psi_0$ with $\psi_n$ satisfies
$$(T\psi_n)(u):=(\psi_0*\psi_n)(u)=\frac{u^n}{\sqrt{n!2^n}}e^{-\frac{u^2}{4}}. $$
\end{proposition}
\begin{proof}
First, we apply the convolution operator $T$ to the Segal-Bargmann kernel $A_z$. On one hand, by Lemma \ref{TBarg} we have 
$$\displaystyle (TA_z)(u)=e^{\frac{zu}{\sqrt{2}}-\frac{u^2}{4}}.$$
Expanding the right-hand side as a power series in $z$, we obtain
 $$\displaystyle  (TA_z)(u)=\sum_{n=0}^{\infty}z^n\frac{u^n}{n!\sqrt{2^n}}e^{-\frac{u^2}{4}}.$$

 On the other hand, since  $T$ is a bounded operator on $L^2(\mathbb{R})$, we also have $$\displaystyle(TA_z)(u)=\sum_{n=0}^{\infty} z^n\frac{(T\psi_n)(u)}{\sqrt{n!}}.$$

  Therefore, by analyticity in the complex variable $z$ we identify the coefficients and obtain

 $$(T\psi_n)(u)=\frac{u^n}{\sqrt{n!2^n}}e^{-\frac{u^2}{4}}, \quad \text{ for all } u\in \mathbb{R},  n\in\mathbb N_0. $$

\end{proof}

\begin{remark}
	The holomorphic extension of the family of functions $$(T\psi_n)(z)=(\psi_0*\psi_n)(z)=\frac{z^n}{\sqrt{n!2^n}}e^{-\frac{z^2}{4}}$$ form an orthonormal basis of the reproducing kernel Hilbert space associated with the complex Gaussian RBF kernel $K_{RBF}$.
\end{remark}

Consider now the holomorphic extension of $(T\psi_n)(u)$, and let us compute the associated Zaremba formula. The calculations lead to the complex Gaussian RBF kernel introduced in \cite{SC2008, SD2006}. Namely, we have the following result

 \begin{proposition} \label{ZarBerRBF}
 The Zaremba-Bergman formula associated with the holomorphic extension of $(T\psi_n)_n$ is given by
 \begin{equation}
 \displaystyle \sum_{n=0}^{\infty}(T\psi_n)(z)\overline{(T\psi_n)(w)}= e^{-\frac{(z-\overline{w})^2}{4}}=K_{RBF}(z,w),
 \end{equation}
 for all $z,w\in \mathbb{C}$. In other words, we have
  \begin{equation}\label{RA}
 \displaystyle \sum_{n=0}^{\infty}(\psi_0*\psi_n)(z)\overline{(\psi_0*\psi_n)(w)}= e^{-\frac{(z-\overline{w})^2}{4}}=K_{RBF}(z,w),
 \end{equation}
 \end{proposition}
 \begin{proof}
By Proposition \ref{kamal1} we have
 \begin{align*}
 \displaystyle  \sum_{n=0}^{\infty}(T\psi_n)(z)\overline{(T\psi_n)(w)}&=\sum_{n=0}^{\infty} \frac{z^n}{\sqrt{n!2^n}}e^{-\frac{z^2}{4}} \frac{\overline{w}^n}{\sqrt{n!2^n}}e^{-\frac{\overline{w}^2}{4}} \\
 &= \sum_{n=0}^{\infty}  \frac{(z\overline{w})^n}{n!2^n}e^{-\frac{(z^2+\overline{w}^2)}{4}} \\
 &=e^{\frac{z\overline{w}}{2}} e^{-\frac{(z^2+\overline{w}^2)}{4}} \\
 &= e^{-\frac{(z-\overline{w})^2}{4}}\\
 &=K_{RBF}(z,w.)
 \end{align*}
 \end{proof}
\begin{remark}
When we restrict the formula \eqref{RA} to the real line, we can recover the standard Gaussian RBF kernel as follows

$$ \sum_{n=0}^{\infty}(\psi_0*\psi_n) \otimes (\psi_0*\psi_n)  (x,y)= e^{-\frac{(x-y)^2}{4}}=K_{RBF}(x,y),\quad x,y\in \mathbb R.$$
\end{remark}
The normalized reproducing kernel of the Fock space $\mathcal{F}(\mathbb C)$ can be obtained as a generating function associated with the family $(\psi_0*\psi_n)_{n\geq 0}$.
\begin{proposition}
For any $z\in\mathbb C$ and $u\in \mathbb R,$ we have 
\begin{equation}
\displaystyle \sum_{n=0}^{\infty}\frac{z^n}{\sqrt{n!}}(\psi_0*\psi_n)(u)=e^{\frac{zu}{\sqrt{2}}-\frac{u^2}{4}}:=k_{1,\frac{u}{\sqrt{2}}}(z).
\end{equation}
\end{proposition}
\begin{proof}
This follows directly from the boundedness of the convolution operator $T$ on $L^2(\mathbb{R})$ and from the fact that $$T\psi_n=\psi_0*\psi_n, \quad n\in\mathbb N_0. $$
\end{proof}

Let $m\in \mathbb{N}_0$ be fixed. We introduce the generating function $\Phi_m(z,u)$ associated with the convolution of normalized Hermite functions $T_m\psi_n=\psi_m*\psi_n$, namely,

\begin{equation}
\displaystyle \Phi_m(z,u):=  \sum_{n=0}^{\infty}\frac{z^n}{\sqrt{n!}} (T_m\psi_n)(u)=\sum_{n=0}^{\infty}\frac{z^n}{\sqrt{n!}}(\psi_m*\psi_n)(u), \quad \text{ for all } z\in \mathbb C, u\in \mathbb R.
\end{equation}

\begin{theorem}\label{Wem}
The generating function $\Phi_m(z,u)$ can be obtained by computing the action of the Weyl operator on the Fock orthonomal basis elements $e_m(z)=\frac{z^m}{\sqrt{m!}}$. Specifically, we have 
\begin{equation}
\Phi_{m}(z,u)=(-1)^m\mathcal{W}_{\frac{u}{\sqrt{2}}}e_m(z), \quad \text{ for all } z\in\mathbb C, u\in \mathbb R, m\in \mathbb N_0.
\end{equation}
Moreover, for $w \in \mathbb{C}$ we have
\begin{equation}
\label{sumher}
\sum_{m=0}^{\infty}\frac{w^m}{\sqrt{m!}}\Phi_m(z,u)= e^{w(\frac{u}{\sqrt{2}}-z)}k_{1,\frac{u}{\sqrt{2}}}(z).
\end{equation}
\end{theorem}

\begin{proof}
Let $m\in \mathbb N_0$ be fixed. We begin by expressing $\Phi_m(z,u)$ as a generating function of the complex Hermite polynomials $(H_{m,n})_{n\in \mathbb{N}_0}$. Indeed, we have
\begin{align*}
\displaystyle 
\Phi_{m}(z,u)&=\sum_{n=0}^{\infty}\frac{z^n}{\sqrt{n!}}(T_m\psi_n)(u)\\
&=\sum_{n=0}^{\infty}\frac{z^n}{\sqrt{n!}}(\psi_m*\psi_n)(u)\\
&=\frac{1}{\sqrt{\pi m! 2^m}} \sum_{n=0}^{\infty} \frac{1}{n!}\left(\frac{z}{\sqrt{2}}\right)^n(h_m*h_n)(u)
\end{align*}

We now apply the Hermite convolution formula \eqref{conv}, and by \eqref{genHer} we obtain
\begin{eqnarray}
	\nonumber
\Phi_{m}(z,u)&=&\frac{e^{-\frac{u^2}{4}}}{\sqrt{\pi m!}} \sum_{n=0}^{\infty}\frac{z^n}{n!} H_{m,n}\left(\frac{u}{\sqrt{2}}, \frac{u}{\sqrt{2}}\right)\\
\label{step}
&=&\frac{(-1)^m}{\sqrt{m!}} \left(z-\frac{u}{\sqrt{2}}\right)^m k_{1,\frac{u}{\sqrt{2}}}(z)\\
\nonumber
&=& (-1)^m\mathcal{W}_{\frac{u}{\sqrt{2}}}e_m(z).
\end{eqnarray}
Finally formula \eqref{sumher} follows by \eqref{step}.

\end{proof}
\begin{corollary}
The following identity holds 
\begin{equation}
\displaystyle \sum_{m,n=0}^{\infty}\frac{w^mz^n}{\sqrt{m!n!}}(\psi_m*\psi_n)(u)= e^{w(\frac{u}{\sqrt{2}}-z)}k_{1,\frac{u}{\sqrt{2}}}(z).
\end{equation}
\end{corollary}
\begin{proof}
This formula is a direct consequence of Theorem \ref{Wem}.
\end{proof}

We conclude this section by briefly discussing the \textit{analytic extension} of the following convolution formula 
\begin{equation}\label{C}
(T_m\psi_n)(u)=(\psi_m*\psi_n)(u)=\frac{e^{-\frac{u^2}{4}}}{\sqrt{m!n!}}H_{m,n}\left(\frac{u}{\sqrt{2}},\frac{u}{\sqrt{2}}\right), \quad \text{ for all } m,n\in \mathbb N_0.
\end{equation}

There is a unique analytic extension of the function $(T_m\psi_n)(u)$ to the complex variable $z$, which leads to the following:
\begin{definition}[Analytic Convolution Kernel]\label{ACK}
Let $m,n\in \mathbb{N}_0$. For every $z\in \mathbb{C}$, define
\begin{align*}
(T_{m,a}\psi_n)(z):=(\psi_m*\psi_n)(z)&=\frac{1}{\sqrt{\sqrt{\pi}m!2^m}}\int_{\mathbb{R}}e^{-\frac{1}{2}(z-y)^2}H_m(z-y)\psi_n(y)dy\\
&=\frac{1}{\sqrt{m!n!}} e^{- \frac{z^2}{4}} H_{m,n} \left( \frac{z}{\sqrt{2}}, \frac{z}{\sqrt{2}}\right).\\
\end{align*}
For all $z,w\in \mathbb C,$ we introduce the analytic convolution kernel function 
\begin{align*}\label{analytic}
\displaystyle \mathbf{K}_m(z,w)&=\displaystyle \sum_{n=0}^{\infty }(T_{m,a}\psi_n)(z) \overline{(T_{m,a}\psi_n)(w)}\\
&=e^{-\frac{z^2}{4}-\frac{\overline{w}^2}{4}}\sum_{n=0}^\infty \frac{1}{m!n!} H_{m,n} \left( \frac{z}{\sqrt{2}}, \frac{z}{\sqrt{2}}\right) \overline{H_{m,n} \left( \frac{w}{\sqrt{2}}, \frac{w}{\sqrt{2}}\right)}.\\
\end{align*}
\end{definition}

\begin{remark}
In the case $m=0$, the analytic convolution kernel $\mathsf K_m$ reduces to the standard Gaussian radial basis function kernel $K_{RBF}$, as obtained in Proposition \ref{ZarBerRBF}.
The study of the analytic convolution extension kernel given by $\mathbf{K}_m$ will be the subject of a subsequent paper. In the following sections of the present work, we focus on the \textit{polyanalytic extension} of the convolution formula \eqref{C}, by considering the polyanalytic complex Hermite polynomials $H_{m,n} \left( \frac{z}{\sqrt{2}}, \frac{\overline{z}}{\sqrt{2}}\right)$ instead of the analytic ones used in Definition \ref{ACK}.
\end{remark}

\section{Polyanalytic Gaussian RBF kernels}
 \setcounter{equation}{0}
 
Given, for instance, the real-valued function
$$
f(x) = 1 - x^2, \quad x \in \mathbb{R},
$$
there are several ways to extend this function to the complex domain. In this paper we focus on these two extensions:
\begin{enumerate}
	\item The analytic (or holomorphic) extension:
$$
	F(z) = 1 - z^2,\quad z\in \mathbb{C}.
$$
	\item The polyanalytic extension:
$$
	G(z, \bar{z}) = 1 - z\bar{z}, \quad z\in \mathbb{C}.
$$
\end{enumerate}
In the previous section, we explored the connection between the convolution operator and the analytic extension of the RBF kernel. In this section we want to tackle the following problems

\begin{problem}
	\label{p1}
Let $m\in \mathbb{N}_0$. Investigate the polyanalytic extension of  the Gaussian RBF kernel by considering the Zaremba formula associated with the family of functions $$(T_m\psi_n)(x)=(\psi_m*\psi_n)(x),$$ where $T_m$ is the convolution operator introduced in Definition \ref{ConvolutionOperator} by choosing $g(x)=\psi_m(x)$, the normalized Hermite function.
\end{problem}
\begin{problem}
	\label{pnew}
What is the function space for which the polyanalytic extension of the Gaussian RBF kernel serves as the reproducing kernel?
\end{problem}

We start focusing on Problem \eqref{p1}. The polyanalytic extension of the convolution $(\psi_m*\psi_n)(x)$ is defined as follows.

\begin{definition}[Polyanalytic Extension]\label{An3}
Let $m,n\in \mathbb{N}_0$. Define
$$(T_{m,p}\psi_n)(z)=(\psi_m*\psi_n)(z,\overline{z}):=\frac{1}{\sqrt{m!n!}} e^{- \frac{z^2}{4}} H_{n,m} \left( \frac{z}{\sqrt{2}}, \frac{\overline{z}}{\sqrt{2}}\right), \quad z\in \mathbb{C}.$$
The function $(T_{m,p}\psi_n)(z)$ is polyanalytic of order $m+1$.
\end{definition}

We now introduce the following true polyanalytic kernel function

\begin{definition}
Let $\ell \in \mathbb{N}$. The \textit{true polyanalytic Gaussian RBF kernel} of order $\ell$ is defined as

\begin{equation}
 \mathsf{K}_{\ell,T}(z,w):=\sum_{n=0}^{\infty}(T_{\ell-1,p}\psi_n)(z)\overline{(T_{\ell-1,p}\psi_n)(w)}=\sum_{n=0}^{\infty} (\psi_{\ell-1}*\psi_n)(z,\overline{z}) \overline{(\psi_{\ell-1}*\psi_n)(w,\overline{w})},
\end{equation}
for all $z,w\in \mathbb{C}.$
\end{definition}

\begin{remark}
For $\ell=1$, we recover the classical Gaussian RBF kernel. The representation of $\psi_{\ell-1}*\psi_n$ in terms of the Itô-Hermite polynomials (see \cite{ADDS2024}) allows us to construct a polyanalytic extension of the Gaussian RBF kernel. This extension is closely related to the polyanalytic Fock kernel.
\end{remark}

As a first step, we consider the restriction of the kernel function $\mathsf{K}_{\ell,T}$ to the real case, given by

\begin{equation}
\displaystyle \mathsf{K}_{\ell,T}(x,y):=\sum_{n=0}^{\infty}(T_{\ell-1}\psi_n)(x)(T_{\ell-1}\psi_n)(y)=\sum_{n=0}^{\infty} (\psi_{\ell-1}*\psi_n)(x) (\psi_{\ell-1}*\psi_n)(y),
\end{equation}
for all $x,y\in \mathbb{R}.$
\begin{proposition}
Let $\ell\in \mathbb{N}$. Then, for all $x,y \in \mathbb{R}$, we have  
\begin{align*}
\displaystyle \mathsf{K}_{\ell,T}(x,y)&=K_{RBF}(x,y) L^{0}_{\ell-1}\left(\frac{|x-y|^2}{2}\right), 
\end{align*}
where $K_{RBF}$ denotes the classical Gaussian RBF kernel for $\gamma=2$, see \eqref{rbfkernel}.
\end{proposition}
\begin{proof}
By the convolution of the normalized Hermite functions, see \eqref{norm}, we have
\begingroup\allowdisplaybreaks
\begin{align*}
\displaystyle \mathsf{K}_{\ell,T}(x,y)&=\sum_{n=0}^{\infty}(T_{\ell-1}\psi_n)(x)(T_{\ell-1}\psi_n)(y)\\
&=\sum_{n=0}^{\infty} (\psi_{\ell-1}*\psi_n)(x) (\psi_{\ell-1}*\psi_n)(y)\\
&=  e^{- \frac{x^2}{4}- \frac{y^2}{4}}\sum_{n=0}^{\infty} \frac{H_{\ell-1,n} \left( \frac{x}{\sqrt{2}}, \frac{x}{\sqrt{2}}\right) H_{n,\ell-1} \left( \frac{y}{\sqrt{2}}, \frac{y}{\sqrt{2}}\right)}{(\ell-1)!n!} 
\end{align*}
\endgroup
Therefore, a direct application of formula \eqref{gen} leads to 

\begin{align*}
\displaystyle \mathsf{K}_{\ell,T}(x,y)&= e^{- \frac{x^2}{4}- \frac{y^2}{4}} e^{\frac{xy}{2}} L_{\ell-1}^{0}\left(\frac{|x-y|^2}{2}\right)\\
&=e^{-\frac{1}{4}(x-y)^2}L_{\ell-1}^{0}\left(\frac{|x-y|^2}{2}\right)\\
&=K_{RBF}(x,y) L^{0}_{\ell-1}\left(\frac{|x-y|^2}{2}\right).
\end{align*}
\end{proof}

Inspired by the above computations, we can provide an extension of the Zaremba–Bergman formula to the polyanalytic setting.

\begin{proposition}\label{An4}
Let $\ell\in \mathbb{N}$. Then, it holds that 
\begin{equation}
\displaystyle \mathsf{K}_{\ell,T}(z,w)=e^{-\frac{1}{4}(z-\overline{w})^2}L_{\ell-1}^{0}\left(\frac{|z-w|^2}{2}\right),
\end{equation}
for all $z,w\in \mathbb{C}$.
\end{proposition}
\begin{proof}
By \eqref{norm} and \eqref{exx} we have
\begin{align*}
\displaystyle \mathsf{K}_{\ell,T}(z,w)&= \sum_{n=0}^{\infty} (\psi_{\ell-1}*\psi_n)(z,\overline{z}) \overline{(\psi_{\ell-1}*\psi_n)(w,\overline{w})}\\
&=  e^{- \frac{z^2}{4}- \frac{\overline{w}^2}{4}}\sum_{n=0}^{\infty} \frac{1}{(m-1)!n!}H_{n,\ell-1} \left( \frac{z}{\sqrt{2}}, \frac{\overline{z}}{\sqrt{2}}\right) \overline{H_{n,\ell-1} \left( \frac{w}{\sqrt{2}}, \frac{\overline{w}}{\sqrt{2}}\right)}\\
&= e^{- \frac{z^2}{4}- \frac{\overline{w}^2}{4}}\left(\overline{\sum_{n=0}^{\infty} \frac{1}{(m-1)!n!} H_{\ell-1,n} \left( \frac{z}{\sqrt{2}}, \frac{\overline{z}}{\sqrt{2}}\right) H_{n,\ell-1} \left( \frac{w}{\sqrt{2}}, \frac{\overline{w}}{\sqrt{2}}\right)}\right).
\end{align*}
Hence, applying formula \eqref{gen} we obtain
\begin{align*}
\displaystyle \mathsf{K}_{\ell,T}(z,w)&=e^{- \frac{z^2}{4}- \frac{\overline{w}^2}{4}} e^{\frac{z\overline{w}}{2}} L_{\ell-1}^{0}\left(\frac{|z-w|^2}{2}\right)\\
&=e^{-\frac{1}{4}(z-\overline{w})^2}L_{\ell-1}^{0}\left(\frac{|z-w|^2}{2}\right)\\
&=K_{RBF}(z,w) L^{0}_{\ell-1}\left(\frac{|z-w|^2}{2}\right).
\end{align*}

\end{proof}
We now introduce another type of polyanalytic Gaussian RBF kernel, using the \emph{polyanalytic extension à la Balk}; see \cite{Balk1991}.

\begin{definition}
Let $N\in \mathbb{N}$. The \textit{polyanalytic Gaussian RBF kernel} of order $N$  is defined as

\begin{equation}
\displaystyle K_{RBF,N}(z,w):=\sum_{m=1}^{N}\mathsf{K}_{m,T}(z,w),
\end{equation}
for all $z,w\in \mathbb{C}.$
\end{definition}
\begin{proposition}
\label{truerep}
Let $M\in \mathbb{N}$. It holds that 
\begin{equation}
\displaystyle K_{RBF,N}(z,w)=K_{RBF}(z,w) L^{1}_{N-1}\left(\frac{|z-w|^2}{2}\right),
\end{equation}
for all $z,w\in \mathbb{C}$.
\end{proposition}
\begin{proof}
Using the well-known formula for the generalized Laguerre polynomials (see \cite{Lebedev1972}, p.29)
\begin{equation}
\label{l1}
\sum_{m=1}^{N}L^0_{m-1}(x)=L^1_{N-1}(x),
\end{equation}
we have 
\begin{align*}
\displaystyle K_{RBF,N}(z,w)&=\sum_{m=1}^{N}\mathsf{K}_{m,T}(z,w)\\
&=K_{RBF}(z,w)\sum_{m=1}^{N}L^{0}_{m-1}\left(\frac{|z-w|^2}{2}\right)\\
&=K_{RBF}(z,w) L^{1}_{N-1}\left(\frac{|z-w|^2}{2}\right).
\end{align*}
\end{proof}
\begin{remark}
For $N=1$ we recover the classical Gaussian RBF kernel function.
\end{remark}
The polyanalytic Gaussian RBF kernel can be expressed using the polyanalytic Fock space kernel as follows:
\begin{proposition}\label{FRBF}
Let $N\in \mathbb{N}$. Then, it holds that 
\begin{equation}\label{FRBFF}
\displaystyle K_{RBF,N}(z,w)=e^{-\frac{z^2}{4}}K^{\frac{1}{2}}_N(z,w)e^{-\frac{\overline{w}^2}{4}},
\end{equation}
for all $z,w\in \mathbb{C}$. Here $K^{\frac{1}{2}}_N(z,w)$ is the polyanalytic Fock kernel of order $N$, see \eqref{kernelP}.
\end{proposition}
\begin{proof}
This follows from the expression of the polyanalytic Fock kernel with  $\alpha=\frac{1}{2}$ given by  $$\displaystyle K^{\frac{1}{2}}_N(z,w)=e^{\frac{z\overline{w}}{2}} L^{1}_{N-1}\left(\frac{|z-w|^2}{2}\right),$$ 
for all $z,w\in \mathbb{C}$
\end{proof}
\begin{remark}
It is easy to see that the kernel function $K_{RBF,N}(z,w)$ is polyanalytic of order $N$ with respect to the variable $z$ for every fixed parameter $w$.
\end{remark}
Now, we turn our attention to the Problem \ref{pnew}. Our goal is to show that the functions $ \mathsf{K}_{m,T}(z, w) $ and $ K_{\mathrm{RBF}, N}(z, w) $ serve as reproducing kernels for the following function spaces.

\begin{definition}[True polyanalytic RBF space]
	\label{trurbf}
	Let $\gamma>0$ and $\ell \in \mathbb{N}$.
	A function $f: \mathbb{C} \to \mathbb{C}$ belongs to the true polyanalytic RBF space if and only if
	$$\frac{2}{\gamma^2 \pi}\int_{\mathbb{C}} |f(z)|^2 e^{\frac{(z-\bar{z})^2}{\gamma^2}}d\lambda(z)<\infty,$$
	and there exists a holomorphic function $h(z)$ such that
\begin{equation}
	\label{cond2}
	f(z)=  \frac{\gamma^{2(\ell-1)}}{2^{\ell-1}}\sqrt{\frac{1}{(\ell-1)!}} e^{\frac{2|z|^2}{\gamma^2}-\frac{z^2}{\gamma^2}} \frac{\partial^{\ell-1}}{\partial z^{\ell-1}} \left( h(z) e^{-\frac{2|z|^2}{\gamma^2}}\right), \quad \forall z \in \mathbb{C}.
\end{equation}
	We denote this space as $\mathcal{H}_{T, \gamma}^{(\ell)}(\mathbb{C})$.
\end{definition}

In the polyanalytic setting it is also possible to define another type of extension of the classical RBF space.

\begin{definition}(Polyanalytic Gaussian RBF space)
\label{polyRBF}
Let $\gamma>0$ and  $N\in \mathbb{N}$. A polyanalytic function $f:\mathbb{C}\longrightarrow \mathbb{C}$ belongs to the polyanalytic Gaussian RBF space of order $N$, denoted by $\mathcal{H}_{N, \gamma}(\mathbb{C})$ (or simply $\mathcal{H}_{N, \gamma})$ if we have
\begin{equation}
\displaystyle ||f||_{\mathcal{H}_{N, \gamma}}^2:=\dfrac{2}{\gamma^2\pi}\int_{\mathbb{C}}|f(z)|^2\exp\left(\frac{(z-\overline{z})^2}{\gamma^2}\right) d\lambda(z)<\infty,
\end{equation}
where $d\lambda(z)=dxdy$ is the Lebesgue measure with respect to the variable $z=x+iy$.
\end{definition}

\begin{remark}
In the sequel, when $ \gamma = 2$, the true polyanalytic RBF space and the polyanalytic Gaussian RBF space will be denoted by $ \mathcal{H}_{T}^{(\ell)}(\mathbb{C})$ and $ \mathcal{H}_{N}(\mathbb{C})$, respectively.

\end{remark}

\begin{remark}
The true polyanalytic RBF space and the polyanalytic RBF space coincide only when $ N = \ell = 1 $. In this case, they both reduce to the classical RBF space; see Definition~\ref{RBF}.

	$$ \mathcal{H}_{0, \gamma}(\mathbb{C})=\mathcal{H}_{T, \gamma}^{(1)}(\mathbb{C})= \mathcal{H}_\gamma^{RBF}(\mathbb{C}).$$
\end{remark}

To establish a connection between the true polyanalytic RBF space and the polyanalytic Gaussian RBF space, it is first necessary to relate the polyanalytic Fock space to the polyanalytic RBF space. To this end, we begin with the following important preliminary result:

\begin{proposition}
	\label{l2iso}
	let $\gamma >0$.
	The multiplication operator $ \mathbf{M}_\gamma: L^2(\mathbb{C}, e^{\frac{(z-\bar{z})^2}{\gamma^2}}) \to  L^2\left(\mathbb{C}, e^{- \frac{2|z|^2}{\gamma^2}}\right)$  given by
	$$ \mathbf{M}_\gamma f(z)= e^{\frac{z^2}{\gamma^2}} f(z), \qquad f \in L^2\left(\mathbb{C}, e^{\frac{(z-\bar{z})^2}{\gamma^2}}\right), \, z \in \mathbb{C}$$
is a unitary operator and its inverse is given by
$$\mathbf{M}^{-1}_\gamma g(z):=e^{-\frac{z^2}{\gamma^2}}g(z), \qquad g \in L^2(\mathbb{C}, e^{- \frac{2|z|^2}{\gamma^2}}).$$
\end{proposition}
\begin{proof}
First we prove that the operator
$$ \mathbf{M}^*_\gamma g(z)= e^{- \frac{z^2}{\gamma^2}} g(z), \qquad g \in L^2 \left(\mathbb{C}, e^{\frac{(z-\bar{z})^2}{\gamma^2}}\right)$$
is the adjoint of $ \mathbf{M}_\gamma$. Let $f_1 \in L^2(\mathbb{C}, e^{\frac{(z-\bar{z}^2)}{\gamma^2}})$ and $f_2 \in L^2(\mathbb{C}, e^{-\frac{2|z|^2}{\gamma^2}})$, then we have
\begin{eqnarray*}
\langle \mathbf{M}_\gamma f_1, f_2 \rangle_{L^2\left(\mathbb{C}, e^{- \frac{2|z|^2}{\gamma^2}}\right)}&=& \int_{\mathbb{R}^2} \overline{\mathbf{M}_\gamma f_1(z)} f_2(z) e^{- \frac{2|z|^2}{\gamma^2}} d \lambda(z)\\
&=& \int_{\mathbb{R}^2} e^{\frac{\overline{z^2}}{\gamma^2}} \overline{f_1(z)} f_2(z) e^{- \frac{2|z|^2}{\gamma^2}} d \lambda(z)\\
&=& \int_{\mathbb{R}^2} e^{\frac{\overline{z^2}}{\gamma^2}+\frac{z^2}{\gamma^2}-\frac{z^2}{\gamma^2}} \overline{f_1(z)} f_2(z)e^{- \frac{2|z|^2}{\gamma^2}} d \lambda(z)\\
&=& \int_{\mathbb{R}^2} \overline{f_1(z)}e^{- \frac{z^2}{\gamma^2}}f_2(z) e^{\frac{(z-\bar{z})^2}{\gamma^2}}d \lambda(z)\\
&=& \langle f_1, \mathbf{M}^*_\gamma f_2 \rangle_{L^2\left(\mathbb{C}, e^{\frac{(z-\bar{z})^2}{\gamma^2}}\right)}.
\end{eqnarray*}
This proves that the operator $ \mathbf{M}^*_\gamma$ is the adjoint of $ \mathbf{M}_\gamma$. Thus it is clear that
$$ \mathbf{M}_\gamma \mathbf{M}^*_\gamma g=g, \qquad g \in L^2 \left(\mathbb{C}, e^{\frac{(z-\bar{z})^2}{\gamma^2}}\right), \qquad \mathbf{M}^*_\gamma \mathbf{M}_\gamma f=f, \qquad f \in L^2 \left(\mathbb{C}, e^{- \frac{2|z|^2}{\gamma^2}}\right).$$
Moreover, it is straightforward that $ \mathbf{M}_\gamma \mathbf{M}^*_\gamma= \mathbf{M}^*_\gamma\mathbf{M}_\gamma$. This proves that the operator $ \mathbf{M}_\gamma$ is unitary. Finally, since the operator $ \mathbf{M}_\gamma$ is unitary we have that $ \mathbf{M}^{*}_\gamma=\mathbf{M}^{-1}_\gamma$.

\end{proof}

\begin{theorem}(polyanalytic RBF-Fock isomorphism)\label{FRBFchara}
		\label{isopar}
Let $\gamma>0$ and $N\in\mathbb{N}$. Then the multiplication operator given by
\begin{equation}
	\mathcal{M}_Nf(z):=\exp(\frac{z^2}{\gamma^2})f(z), \quad  f\in \mathcal{H}_{N, \gamma}(\mathbb{C}), \quad z\in \mathbb{C}.
\end{equation}
is a unitary operator between the polyanalytic gaussian RBF and Fock space with parameter $\alpha=\frac{2}{\gamma^2}$.
\end{theorem}
\begin{proof}
Set $\displaystyle g(z)=\exp(\frac{z^2}{\gamma^2})f(z)$ for every $z\in\mathbb{C}$. Note that $g$ is a polyanalytic function of order $N$ as a product of an entire function with a polyanalytic function of order $N$. Finally since the space $ \mathcal{H}_N(\mathbb{C})$ is the space of polyanalytic functions that belongs to the space $L^2\left(\mathbb{C}, e^{\frac{(z-\bar{z})^2}{\gamma^2}}\right)$ and $ \mathcal{F}_{\frac{2}{\gamma^2}, N}(\mathbb{C})$ is the space of polyanalytic function that belongs to $L^2\left(\mathbb{C}, e^{-\frac{2|z|^2}{\gamma^2}}\right) $, the result follows by Proposition \ref{l2iso}.

\end{proof}

We now have all the necessary tools to establish a connection between the true polyanalytic RBF space and the polyanalytic RBF space.

\begin{theorem}
	\label{sums}
Let $\gamma>0$ and	$N \in \mathbb{N}$.
	The polyanalytic RBF space  $ \mathcal{H}_N(\mathbb{C})$ is a direct sum of true polyanalytic RBF space $ \mathcal{H}_{T, \gamma}^{(\ell)}(\mathbb{C})$, with $\ell=1$,..., $N$, i.e.
\begin{equation}
\label{summ}
\mathcal{H}_{N, \gamma}(\mathbb{C})= \bigoplus_{\ell=1}^N \mathcal{H}_{T, \gamma}^{(\ell)}(\mathbb{C}).
\end{equation}
\end{theorem}
\begin{proof}
	We prove the result by double inclusion. Let $f \in \mathcal{H}_{N, \gamma}(\mathbb{C})$. By Theorem \ref{isopar} we have that $\mathcal{M}_N f \in  \mathcal{F}_{\frac{2}{\gamma^2},N}(\mathbb{C})$. By \eqref{VasDecom} we get that $ \mathcal{M}_N f \in \bigoplus_{\ell=1}^N \mathcal{F}_{T, \frac{2}{\gamma^2}}^{(\ell)}(\mathbb{C})$. This implies that there exist unique $f_\ell \in \mathcal{F}_{T, \frac{2}{\gamma^2}}^{(\ell)}(\mathbb{C})$, with $\ell=1$,..., $N$, such that
	$$ \mathcal{M}_N f(z)= \sum_{\ell=1}^{N} f_k(z).$$
	By the definition of the multiplication operator $ \mathcal{M}_N$ we have
	\begin{equation}
		\label{summs}
		f(z)= \sum_{\ell=1}^{N}f_\ell(z) e^{-\frac{z^2}{\gamma^2}}.
	\end{equation}
	Now, we set $g_\ell(z):=f_\ell(z) e^{-\frac{z^2}{4}}$, for any $\ell=1,...,N$. By using the fact that $f_\ell \in \mathcal{F}_{T, \frac{2}{\gamma^2}}^{(\ell)}(\mathbb{C})$ we have
	$$ \frac{2}{\gamma^2\pi}\int_{\mathbb{C}} | g_\ell(z)|^2 e^{\frac{(z-\bar{z})^2}{\gamma^2}}= \frac{2}{\gamma^2\pi} \int_{\mathbb{C}} |f_\ell(z)|^2 e^{-\frac{2|z|^2}{\gamma^2}}<\infty,$$
	and by \eqref{Bpoly2} we have that there exists a holomorphic function $h$ such that
	$$ g_\ell(z)= f_\ell(z) e^{-\frac{z^2}{\gamma^2}}= \frac{\gamma^{2(\ell-1)}}{2^{\ell-1}}\sqrt{\frac{1}{(\ell-1)!}} e^{\frac{2|z|^2}{\gamma^2}-\frac{z^2}{\gamma^2}} \frac{\partial^{\ell-1}}{\partial z^{\ell-1}} \left( h(z) e^{-\frac{2|z|^2}{\gamma^2}}\right).$$
	This proves that $g_\ell \in \mathcal{H}^{(\ell)}_{T,\gamma}(\mathbb{C})$, for any $\ell=1,...,N$, by Definition \ref{trurbf}. Hence by \eqref{summs} we have $f \in \bigoplus_{k=0}^N \mathcal{H}_{T,\gamma}^{(\ell)}(\mathbb{C})$. 
	\\Now, let $F \in \bigoplus_{\ell=1}^N \mathcal{H}_{T, \gamma}^{(\ell)}(\mathbb{C})$, then there exists $p_\ell \in \mathcal{H}_{T, \gamma}^{(\ell)}(\mathbb{C})$, with $\ell=1,...,N$ such that
	$$ F(z)= \sum_{\ell=1}^{\infty} p_\ell(z).$$
	We multiply both sides of the above equality by $ e^{\frac{z^2}{\gamma^2}}$ and we have
	\begin{equation}
		\label{summ2}
		F(z)e^{\frac{z^2}{\gamma^2}}= \sum_{\ell=1}^{\infty} p_\ell(z)e^{\frac{z^2}{\gamma^2}}.
	\end{equation}
	We set $h_\ell(z):=p_\ell(z)e^{\frac{z^2}{\gamma^2}}$ for $\ell=1,...,N$. By using the fact that $p_k$ belongs to the true RBF we have
	$$ \frac{2}{\gamma^2 \pi}\int_{\mathbb{C}} | h_\ell(z)|^2 e^{-\frac{2|z|^2}{\gamma^2}}= \frac{1}{2 \pi} \int_{\mathbb{C}} |p_\ell(z)| e^{\frac{(z-\bar{z})^2}{\gamma^2}}<\infty,$$
	and by \eqref{cond2} there exists a holomorphic function $h(z)$ such that
	$$h_\ell(z)=p_\ell(z)e^{\frac{z^2}{\gamma^2}}= 2^{\ell-1} \sqrt{\frac{1}{(\ell-1)!}} e^{\frac{ 2|z|^2}{\gamma^2}} \frac{\partial^{\ell-1}}{\partial z^{\ell-1}} \left( h(z) e^{-\frac{2|z|^2}{\gamma^2}}\right).$$
	Thus by Definition \ref{FT} we have that $h_\ell(z) \in \mathcal{F}_{T, \frac{2}{\gamma^2}}^{(\ell)}(\mathbb{C})$. Therefore by \eqref{summ2} we have $ F(z)e^{\frac{z^2}{\gamma^2}} \in \bigoplus_{\ell=1}^N \mathcal{F}_{T, \frac{2}{\gamma^2}}^{(\ell)}(\mathbb{C})$. By \eqref{VasDecom} we have that $  F(z)e^{\frac{z^2}{\gamma^2}} \in \mathcal{F}_{\frac{2}{\gamma^2}, N}(\mathbb{C})$. Finally, by Theorem \ref{isopar} we have $F \in \mathcal{H}_N(\mathbb{C})$.
\end{proof}

The connection between the true RBF space and the polyanalytic RBF space leads to the following relation between the true RBF space and the true Fock space.

\begin{proposition}
	\label{iso}
	Let $\gamma>0$ and $\ell \in \mathbb{N}$. Then there exists an isometric isomorphism between the true polyanalytic RBF and true Fock space given by the multiplication operator $\mathbf{M}: \mathcal{H}_{T, \gamma}^{(\ell)}(\mathbb{C}) \to \mathcal{F}_{T, \frac{2}{\gamma^2}}^{(\ell)}(\mathbb{C})$ defined by 
	
\begin{equation}
\label{iso3}
\mathbf{M}_{\gamma} f(z)= e^{\frac{z^2}{\gamma^2}} f(z), \qquad  f \in \mathcal{H}_{T, \gamma}^{(\ell)}(\mathbb{C}), \, z \in \mathbb{C}.
\end{equation}
\end{proposition}
\begin{proof}
The result follows by using arguments similar to those used in the proofs of Theorems~\ref{sums} and~\ref{isopar}.

\end{proof}

The infinite direct sum of true polyanalytic RBF spaces is deeply connected to the space $ L^2\left(\mathbb{C}, e^{\frac{(z - \bar{z})^2}{\gamma^2}} \right) $, as demonstrated in the following result.

\begin{proposition}
	Let $\gamma>0$. For and $z \in \mathbb{C}$ we have
	\begin{equation}
		\label{l2ineq}
		L^2\left(\mathbb{C}, e^{\frac{(z-\bar{z})^2}{\gamma^2}}\right)=\bigoplus_{\ell=1}^\infty \mathcal{H}_{T, \gamma}^{(\ell)}(\mathbb{C}).
	\end{equation}
\end{proposition}
\begin{proof}
	We prove the equality \eqref{l2ineq} by double inclusion. Let $f \in L^2\left(\mathbb{C}, e^{\frac{(z-\bar{z})^2}{\gamma^2}}\right)$. Then by Proposition \ref{l2iso} we have that there exists $g \in L^2\left(\mathbb{C}, e^{-\frac{2|z|^2}{\gamma^2}}\right)$ such that
	$$ f(z)= e^{-\frac{z^2}{\gamma^2}}g(z).$$
	By \eqref{l2fock} there exists unique $h_\ell \in \mathcal{F}_{T, \frac{2}{\gamma^2}}^{(\ell)}(\mathbb{C})$, for any $\ell \in \mathbb{N}$, such that
	$$ g(z)= \sum_{\ell=1}^{\infty} h_\ell(z).$$
	So we can write
	\begin{equation}
		\label{sums1}
		f(z)= \sum_{\ell=1}^{\infty} e^{-\frac{z^2}{\gamma^2}} h_\ell(z).
	\end{equation}
	So by Proposition \ref{iso}, for any $k \in \mathbb{N}$, we have $e^{-\frac{z^2}{\gamma^2}} h_k(z) \in \mathcal{H}_{T, \gamma}^{(\ell)}(\mathbb{C})$. Thus by \eqref{sums1} we get that $ f \in \oplus_{\ell=1}^\infty \mathcal{H}_{T}^{(\ell)}(\mathbb{C})$. Now, we prove the other equality. Let $f \in \oplus_{\ell=1}^\infty \mathcal{H}_{T,\gamma}^{(\ell)}(\mathbb{C})$. Then there exist unique $g_\ell \in \mathcal{H}_{T, \gamma}^{(\ell)}(\mathbb{C})$ for any $k \in \mathbb{N}$ such that
	$$ f(z)= \sum_{\ell=1}^{\infty} g_\ell(z).$$
	By Proposition \ref{iso} we deduce that there exist $h_k \in F_{T, \frac{2}{\gamma^2}}^{(\ell)}(\mathbb{C})$, for any $k \in \mathbb{N}$, such that
	\begin{equation}
		\label{step0}
		f(z)= \sum_{\ell=1}^{\infty}e^{\frac{z^2}{\gamma^2}}h_\ell(z).
	\end{equation}
	By \eqref{l2fock} we have 
	$$ \sum_{\ell=1}^{\infty}h_\ell(z) \in L^2(\mathbb{C}, e^{-\frac{2|z|^2}{\gamma^2}}).$$
	Finally by \eqref{step0} and Proposition \ref{l2iso} we get $f \in L^2\left(\mathbb{C}, e^{\frac{(z-\bar{z})^2}{\gamma^2}}\right)$.
\end{proof}

The above results pave the way to provide a final solution to the Problem~\ref{pnew}.

\begin{theorem}(Polyanalytic-RBF reproducing kernel property)\label{RKP-RBF}
	The polyanalytic RBF space $\mathcal{H}_{N}$  is a reproducing kernel Hilbert space with reproducing kernel $K_{RBF, N}(z,w)$. The reproducing property is given by the following integral representation
	\begin{equation}\label{RPRBF}
		\displaystyle f(w)=\dfrac{1}{2\pi}\int_\mathbb{C} f(z)\overline{K_{RBF, N}(z,w)}\exp\left(\frac{(z-\overline{z})^2}{4}\right)d\lambda(z), \quad f\in\mathcal{H}_N, w\in\mathbb{C}.
	\end{equation}
	Moreover, the true polyanalytic RBF space is a reproducing kernel Hilbert space whose reproducing kernel is given by
	\begin{equation}
		\label{rep}
		\mathsf{K}_{ \ell,T}(z,w)=L_{\ell-1}^0\left(\frac{|z-w|^2}{2}\right)e^{- \frac{(z-\bar{w})^2}{4}}, \qquad \ell=1,...,N
	\end{equation}
	where $L_{n}^0(x)$ are the Laguerre polynomials, see \eqref{Laguerre}.
\end{theorem}
\begin{proof}
	By using the expression of $K_{RBF, N}(z,w)$, see Proposition \ref{FRBF}, we get \[ \begin{split}
		\frac{1}{2\pi }\int_\mathbb{C} f(z)\overline{K_{RBF, N}(z,w)}&e^{\frac{(z-\overline{z})^2}{4}}d\lambda(z) \\&=  \frac{1}{2\pi }\int_\mathbb{C} f(z)e^{-\frac{\overline{z}^2+w^2}{4}}\overline{K^{\frac{1}{2}}_N(z,w)}
		e^{\frac{(z-\overline{z})^2}{4}}d\lambda(z)\\
		& =\frac{1}{2\pi}e^{-\frac{w^2}{4}}\int_\mathbb{C} e^{\frac{z^2}{4}} f(z)\overline{K^{\frac{1}{2}}_N(z,w)}e^{-\frac{|z|^2}{2}}d\lambda(z)\\
		& = e^{-\frac{w^2}{4}} \frac{1}{2\pi}\int_\mathbb{C} (\mathcal{M}_Nf)(z)\overline{K^{\frac{1}{2}}_N(z,w)}e^{-\frac{|z|^2}{2}}d\lambda(z).
	\end{split}
	\]
	
	However, we already know from Theorem \ref{FRBFchara} that $\mathcal{M}_N$  is an isometric isomorphism between the polyanalytic Fock and RBF spaces. Since $f\in\mathcal{H}_N$ then we have $ \mathcal{M}_Nf\in\mathcal{F}_{1/2,N}$ . We now apply the polyanalytic Fock reproducing kernel property and the explicit expression of $\mathcal{M}_N$ to conclude that
	$$ \frac{1}{2\pi }\int_\mathbb{C} f(z)\overline{K_{RBF, N}(z,w)}e^{\frac{(z-\overline{z})^2}{4}}d\lambda(z)= e^{-\frac{w^2}{4}}(\mathcal{M}_N f)(w)=f(w).$$
	
	Thus, we have shown that $f$ satisfies the reproducing kernel property. Finally, the fact that \eqref{rep} is the reproducing kernel of the true RBF space follows from Theorem~\ref{sums} and the property~\eqref{l1} of the Laguerre polynomials.
	
\end{proof}

\begin{remark} \label{An}
By employing arguments similar to those used in Theorem \ref{RKP-RBF}, one can show that  
\begin{equation} \label{KtrueRBF}
\mathsf{K}_{ \ell,T}^\gamma(z,w):=L_{\ell-1}^0\!\left(\frac{2|z-w|^2}{\gamma^2}\right)
e^{- \frac{(z-\bar{w})^2}{\gamma^2}}, \qquad z,w \in \mathbb{C},
\end{equation}  
and  
\begin{equation}
\label{KpolyRBF}
K_{RBF, N}^\gamma(z,w):=L_{N-1}^1\!\left(\frac{2|z-w|^2}{\gamma^2}\right)
e^{- \frac{(z-\bar{w})^2}{\gamma^2}}, \qquad z,w \in \mathbb{C},
\end{equation}
are the reproducing kernels of the true (parametric) polyanalytic RBF space and the Gaussian RBF polyanalytic space, respectively, for a generic parameter $\gamma > 0$.

\end{remark}

\section{An operator approach to polyanalytic Gaussian RBF kernels}
 \setcounter{equation}{0}

In this section, our goal is to investigate the connection between the classical RBF space (see Definition \ref{RBF}) and the true polyanalytic RBF space (see Definition \ref{trurbf}). This connection can be established by suitably constructing an operator that relates to a Schrödinger-type operator. We begin by recalling the counterparts of the creation and annihilation operators in the context of the RBF space; see~\cite{ACDSS}.

\begin{definition}[Creation and annihilation operator for the RBF space]
The creation operator for the RBF space, denoted by $M_z^{RBF}$, is defined in a way that the following diagram is commutative:
	\begin{equation}
	\label{diag1}
	\xymatrix{
		\mathcal{H}_{2}(\mathbb{C}) \ar[r]_{M_z^{RBF}} \ar[d]_{\mathcal{M}_2} & \mathcal{H}_{2}(\mathbb{C})  \\ \mathcal{F}_{1/2}(\mathbb{C}) \ar[r]_{M_z} & \mathcal{F}_{1/2}(\mathbb{C}) \ar[u]_{\mathcal{M}^{-1}_2},
	}
\end{equation}
where the operator $\mathcal{M}_2$ is the operator that connects the RBF space and the Fock space, see \eqref{RBFfock}.
	The annihilation operator for the RBF space, denoted by $D_z^{RBF}$, is defined in a way that the following diagram is commutative:
	\begin{equation}
	\label{diag2}
	\xymatrix{
		\mathcal{H}_{2}(\mathbb{C}) \ar[r]_{D_z^{RBF}} \ar[d]_{\mathcal{M}_2} & \mathcal{H}_{2}(\mathbb{C})  \\ \mathcal{F}_{1/2}(\mathbb{C}) \ar[r]_{D_z} & \mathcal{F}_{1/2}(\mathbb{C}) \ar[u]_{\mathcal{M}^{-1}_2}
	}
\end{equation}
\end{definition}

Now, by using the definition of  the operator $\mathcal{M}_2$ we can provide an explicit expression of the operators $M_z^{RBF}$ and $D_z^{RBF}$.

\begin{proposition}
For $f \in \mathcal{H}_2(\mathbb{C})$ the creation operator for the RBF space can be written as
	\begin{equation}
		\label{Mrbf}
		M_z^{RBF} f(z)=z f(z).
	\end{equation}
	The annihilation operator for the RBF space, for $f \in \mathcal{H}_2(\mathbb{C})$, can be written as
	\begin{equation}
		\label{Drbf}
		D_z^{RBF}f(z)= \frac{d}{dz} f(z)+ \frac{z}{2}f(z).
	\end{equation}
\end{proposition}

\begin{proof}
	We start proving \eqref{Mrbf}. By the diagram \eqref{diag1} and the definition of the operator $\mathcal{M}_2$ we have
	\begin{eqnarray*}
		M_{z}^{RBF}f(z)&=& \left(\mathcal{M}^{-1}_2M_z \mathcal{M}_2\right)f(z)\\
		&=&\mathcal{M}^{-1}_2 M_z \left( e^{\frac{z^2}{4}} f(z)\right)\\
		&=& e^{-\frac{z^2}{4}} z e^{\frac{z^2}{4}}f(z)\\
		&=& z f(z).
	\end{eqnarray*}
	Now, we prove \eqref{Drbf}. By the diagram \eqref{diag2} and the definition of the operator $\mathcal{M}_2$ we have
	\begin{eqnarray*}
		D_{z}^{RBF}f(z)&=& \left(\mathcal{M}^{-1}_2D_z \mathcal{M}_2\right)f(z)\\
		&=&\mathcal{M}^{-1}_2 \frac{d}{dz} \left( e^{\frac{z^2}{4}} f(z)\right)\\
		&=& e^{-\frac{z^2}{4}} \frac{d}{dz}\left( e^{\frac{z^2}{4}}  f(z)\right)\\
		&=& \frac{d}{dz}f(z)+ \frac{z}{2} f(z).
	\end{eqnarray*}
\end{proof}

\begin{remark}
	Since the creation operator in the RBF case coincides with the one in the Fock space we keep denoting  $M_z^{RBF} f(z)$ by $M_zf(z)$.
\end{remark}

Now, we provide an orthonormal basis for the true polyanalytic RBF space.

\begin{theorem}
	\label{bases}
	Let $\ell \in \mathbb{N}$ fixed.
	 An orthonormal basis of the true polyanalytic RBF space $\mathcal{H}_{T, 2}^{(\ell)}(\mathbb{C})$ is given by
	$$e_{\ell,m}(z):= \frac{2^{\frac{\ell+m}{2}}}{\sqrt{ (\ell-1)! m!} }H_{\ell-1,m}^{1/2}(z, \bar{z}) e^{-\frac{z^2}{4}}, \quad m\in \mathbb{N}_0.$$ 
	
\end{theorem}
\begin{proof}
	First of all we have to compute the following integral
	\begin{equation}
		\label{ort1}
		\langle e_{\ell,m}, e_{\ell,m'} \rangle_{\mathcal{H}_{T, 2}^{(\ell)}}= \frac{1}{2\pi} \int_{\mathbb C} e_{\ell,m}(z, \bar{z}) \overline{ e_{\ell,m'}(z, \bar{z})}e^{\frac{(z-\bar{z})^2}{4}}d \lambda(z),
	\end{equation}
	where $m' \in \mathbb{N}_0$.
			By \eqref{ort} we have
	\begin{eqnarray}
		\nonumber
		\langle e_{\ell,m}, e_{\ell,m'} \rangle_{\mathcal{H}_{T}^{(\ell)}}&=&   \frac{2^{\ell+m-2}}{\pi (\ell-1)! m!}\int_{\mathbb{C}} H_{\ell-1,m}^{1/2}(z, \bar{z}) \overline{H_{\ell-1,m'}^{1/2}(z, \bar{z})} e^{-\frac{z^2}{4}}e^{-\frac{\bar{z}^2}{4}}e^{\frac{(z-\bar{z})^2}{4}}d \lambda(z)\\
		\nonumber
				&=& \frac{2^{\ell+m-2}}{\pi (\ell-1)! m!}\int_{\mathbb{C}} H_{\ell-1,m}^{1/2}(z, \bar{z}) \overline{H_{\ell-1,m'}^{1/2}(z, \bar{z})} e^{-\frac{|z|^2}{2}} d \lambda(z)\\
				\label{ort11}
		&=& \delta_{m,m'}.
	\end{eqnarray}
	This implies that $e_{\ell,m}(z, \bar{z})$ belongs to the true  polyanalytic RBF space and two of them  are orthogonal. Moreover by \eqref{rodri} we have
$$
e_{\ell,m}(z, \bar{z})=\frac{2^{\frac{\ell+m}{2}} (-1)^{\ell-1+m}}{\sqrt{(\ell-1)! m!}} \frac{\partial^{\ell-1+m}}{\partial z^{\ell-1} \partial \bar{z}^m} \left(e^{-\frac{|z|^2}{2}}\right).
$$
Now, by using the fact that $\frac{\partial^m}{\partial z^m} \left(e^{-\frac{|z|^2}{2}}\right)= z^m e^{-\frac{|z|^2}{2}}(-1)^m 2^m$ we get
$$e_{\ell,m}(z, \bar{z})=\frac{(-1)^{\ell-1}}{\sqrt{m! 2^m}} \frac{2^{\frac{\ell}{2}}}{\sqrt{(\ell-1)!}} e^{\frac{|z|^2}{2}-\frac{z^2}{4}} \frac{\partial^{\ell-1}}{\partial z^{\ell-1}} \left(z^m e^{-\frac{|z|^2}{2}}\right).$$
Together with \eqref{ort11}, this implies that $e_{\ell,m}(z,\bar{z})$ belongs to the true polyanalytic RBF space.

Now our goal is to prove that $e_{\ell,m}(z,\bar{z})$ forms a basis. Let $f \in \mathcal{H}_{T}^{(\ell)}(\mathbb{C})$, then the function $g(z):= e^{\frac{z^2}{4}}f(z)$ belongs to the true polyanalytic Fock space, see Proposition \ref{iso}. So we can write
	\begin{equation}
		\label{exp11}
		f(z)= \sum_{m=0}^{\infty} H_{\ell-1,m}^{1/2}(z, \bar{z})e^{- \frac{z^2}{4}}a_m= \sum_{m=0}^{\infty}\sqrt{(\ell-1)!m!} 2^{- \frac{n+m}{2}}e_{\ell,m}(z, \bar{z}) a_m, \quad \{a_m\}_{m=0}^\infty \subseteq \mathbb{C}.
	\end{equation}
	This means that $e_{\ell,m}(z, \bar{z})$ are generators. Now, we prove that $e_{\ell,m}(z, \bar{z})$ are independent. By using formula \eqref{ort1} and the expansion in series of the function $f$ given by \eqref{exp11}, for $k=1,...\ell$, we get
	\begin{eqnarray*}
		\langle f, e_{\ell,m}(z) \rangle_{\mathcal{H}_{T}^{(\ell)}}&=& \frac{1}{2\pi} \int_{\mathbb{C}} \overline{f(z)} e_{\ell,m}(z, \bar{z}) e^{\frac{(z-\bar{z})^2}{4}}d \lambda(z)\\
		&=& \frac{1}{2\pi} \int_{\mathbb{C}} \left(\sum_{m'=0}^{\infty} \sqrt{m'\ell}2^{- \frac{\ell'+m'}{2}}\overline{ e_{\ell, m'}(z, \bar{z})} a_{m'}\right)e_{\ell,m}(z, \bar{z}) e^{\frac{(z-\bar{z})^2}{4}}d \lambda(z)\\
		&=&\sum_{m'=0}^{\infty} \sqrt{m'\ell}2^{- \frac{\ell+m'}{2}}\left[ \frac{1}{2\pi} \int_{\mathbb{C}} \overline{ e_{\ell, m'}(z, \bar{z})} e_{\ell,m}(z, \bar{z}) e^{\frac{(z-\bar{z})^2}{4}}d \lambda(z)\right]a_{m'}\\
		&=& m! 2^{-m}a_m.
	\end{eqnarray*}
	We can exchange the series and the integral since the series expansion \eqref{exp11} converges uniformly on the ball $B(0,r)$, with $r>0$. Therefore $\langle f, e_{\ell,m} \rangle_{\mathcal{H}_{T}^{(\ell)}}=0$ if and only if for all $\ell\in \mathbb N$ and $m \in \mathbb{N}_0$ we have that $a_m=0$, this means that $f=0$.
\end{proof} 

\begin{remark}
For $\ell=1$, we have
$$ e_{1,m}=\frac{z^m}{\sqrt{m!} 2^{m/2}} e^{- \frac{z^2}{4}},$$
that is the orthonormal basis of the complex Gaussian RBF space given in \cite{SC2008, SD2006}.
\end{remark}

Now, we study the connection between the holomorphic RBF space and the true polyanalytic RBF space. This can be achieved by the following operator:

\begin{definition}
	Define the operator $\mathcal{A}_{z, \bar{z}}$ by
	\begin{equation}
		\label{opem}
		\mathcal{A}_{z, \bar{z}}:=\left(-D_z^{RBF}+\frac{M_{\bar{z}}^{RBF}}{2}\right)=\left(- \frac{d}{dz}- \frac{M_z}{2}+\frac{M_{\bar{z}}}{2}\right).
	\end{equation}
\end{definition} 

\begin{proposition}
Let $\ell \in \mathbb{N}$ fixed. Then we have

	\begin{equation}
	\label{form}
	\mathcal{A}_{z, \bar{z}}^{\ell-1} \left( z^m e^{- \frac{z^2}{4}}\right)= 2^m H_{\ell-1,m}^{1/2}(z, \bar{z})e^{- \frac{z^2}{4}}, \quad m \in \mathbb{N}_0.
\end{equation}
Moreover, the operator $ \mathcal{A}_{z, \bar{z}}^{\ell-1}$ maps the analytic Gaussian RBF space $\mathcal{H}_{2}(\mathbb{C})$ to the true polyanalytic Gaussian RBF space $\mathcal{H}_{T, 2}^{(\ell)}(\mathbb{C})$.
\end{proposition}
\begin{proof}
First, we prove \eqref{form} by induction on $\ell$. For $\ell=1$  the formula is trivial, since $H_{0,m}^{1/2}(z, \bar{z})= \left(\frac{z}{2}\right)^m.$
We suppose formula \eqref{form} is valid for $n$ and we prove it for $n+1$. By \eqref{oper} we have		
\begingroup\allowdisplaybreaks
\begin{eqnarray*}
	\mathcal{A}_{z, \bar{z}}^{\ell} \left( z^m e^{-\frac{z^2}{2}}\right)&=& \mathcal{A}_{z, \bar{z}} \left( 2^m H_{\ell-1,m}^{1/2}(z, \bar{z})e^{- \frac{z^2}{4}}\right) \\\
	&=&2^m \left[- \frac{d}{dz}\left( H_{\ell-1,m}^{1/2}(z, \bar{z}) e^{-\frac{z^2}{4}}\right)- \frac{z}{2} H_{\ell-1,m}^{1/2}(z, \bar{z}) e^{-\frac{z^2}{4}}+ \frac{\bar{z}}{2} H_{\ell-1,m}^{1/2}(z, \bar{z}) e^{-\frac{z^2}{4}}\right] \\
	&=&2^m \left[- \frac{d}{dz} \left(H_{\ell-1,m}^{1/2}(z, \bar{z})\right)+\frac{z}{2} H_{\ell-1,m}^{1/2}(z, \bar{z}) \right.\\
	&& \left. -\frac{z}{2} H_{\ell-1,m}^{1/2}(z, \bar{z}) +\frac{\bar{z}}{2} H_{\ell-1,m}(z, \bar{z}) \right]e^{-\frac{z^2}{4}}\\
	&=&2^m \left[\left(\frac{\bar{z}}{2}- \frac{d}{dz}\right)(H_{\ell-1,m}^{1/2}(z, \bar{z}))\right]e^{-\frac{z^2}{4}}\\
	&=& 2^m \left[\left(\frac{\bar{z}}{2}- \frac{d}{dz}\right)^{\ell}\left(\frac{z}{2}\right)^m\right]e^{-\frac{z^2}{4}}\\
	&=& 2^mH_{\ell,m}^{1/2}(z, \bar{z}) e^{-\frac{z^2}{4}}.
\end{eqnarray*}
\endgroup
This proves \eqref{form}.
	Let $f \in \mathcal{H}_{2}(\mathbb{C})$, then we have the following expansion in series:
	\begin{equation}
		\label{expp}
		f(z)= \sum_{m=0}^{\infty} z^m e^{- \frac{z^2}{4}}a_m, \qquad \{a_m\}_{m \in \mathbb{N}_0} \subseteq \mathbb{C}.
	\end{equation}
	We apply the operator $ \mathcal{A}_{z, \bar{z}}^{\ell-1}$  to the above expansion of the function $f$, and by \eqref{form} we get
	
	$$\mathcal{A}_{z, \bar{z}}^{\ell-1}f(z)= \sum_{m=0}^{\infty} 2^mH_{\ell-1,m}^{1/2}(z, \bar{z}) e^{-\frac{z^2}{4}} a_m:=\sum_{m=0}^{\infty}e_{\ell,m}(z, \bar{z}) b_m(\ell):=F(z)$$
	where $b_m(\ell):=\sqrt{(\ell-1)! m!} 2^{\frac{\ell-m}{2}}a_m$. By Proposition \ref{bases} we get that the function $F(z)$ belongs to the true polyanalytic RBF space $\mathcal{H}_{T, 2}^{(\ell)}(\mathbb{C})$, so this proves the result.

\end{proof}

Applying the $(\ell-1)$-th power of the operator $\mathcal{A}_{z, \bar{z}}$ to the reproducing kernel of the RBF space yields the reproducing kernel of the true polyanalytic RBF space, as demonstrated in the following result. 
\begin{proposition}
	Let $\ell \in \mathbb{N}$. Then for $z$, $w \in \mathbb{C}$ we have
	$$ \mathcal{A}_{z, \bar{z}}^{\ell-1} \overline{\mathcal{A}_{w, \bar{w}}^{\ell-1}}\left( e^{- \frac{(z-\bar{w})^2}{4}} \right)=\frac{(\ell-1)!}{2^{\ell-1}}L_{\ell-1}^0\left(\frac{|z-w|^2}{2}\right)e^{- \frac{(z-\bar{w})^2}{4}}.$$
\end{proposition}
\begin{proof}
	By \cite[Prop. 3.5]{ACDSS} we know the following expansion in series for the reproducing kernel of the RBF space:	
	\begin{equation}
		\label{exp}
		e^{- \frac{(z-\bar{w})^2}{4}}= \sum_{m=0}^{\infty} \frac{1}{2^m m!} \left(z^m e^{-\frac{z^2}{4}}\right)\overline{\left(w^m e^{-\frac{w^2}{4}}\right)}.
	\end{equation}
	Now, we apply  the operators $\mathcal{A}_{z, \bar{z}}^{\ell-1}$ and $\overline{\mathcal{A}_{w, \bar{w}}^{\ell-1}}$ to \eqref{exp}. First, we observe that any complex valued function $f$ satisfies 
	\begin{equation}\label{AN2}
	\overline{ \mathcal{A}_{z, \bar{z}}} (\bar{f})= \overline{\mathcal{A}_{z, \bar{z}}f}. 
	\end{equation}
	
	Then, by \eqref{gen} (with $\alpha=\frac{1}{2}$), \eqref{AN2}, and \eqref{form} we get
	\begin{eqnarray*}
		\mathcal{A}_{z, \bar{z}}^{\ell-1} \overline{\mathcal{A}_{w, \bar{w}}^{\ell-1}} \left( e^{- \frac{(z-\bar{w})^2}{4}} \right)&=& \sum_{m=0}^{\infty} \frac{1}{  2^m m!} \mathcal{A}_{z, \bar{z}}^{\ell-1}\left(z^m e^{-\frac{z^2}{4}}\right)\overline{\mathcal{A}_{w, \bar{w}}^{\ell-1}\left(w^m e^{-\frac{w^2}{4}}\right)}\\
		&=&\sum_{m=0}^{\infty} \frac{2^m}{m!} H_{\ell-1,m}^{1/2}(z, \bar{z}) e^{-\frac{z^2}{4}} \overline{H_{\ell-1,m}^{1/2}(w, \bar{w}) e^{-\frac{w^2}{4}}}\\
		&=& \left(\sum_{m=0}^{\infty} \frac{2^m}{m!}  H_{\ell-1,m}^{1/2}(z, \bar{z}) \overline{H_{m,\ell-1}^{1/2}(w, \bar{w})}  \right) e^{-\frac{z^2}{4}} e^{-\frac{\bar{w}^2}{4}}\\
		&=&\frac{(\ell-1)!}{2^{\ell-1}}L^{0}_{\ell-1}\left(\frac{|z-w|^2}{2}\right)e^{- \frac{(z-\bar{w})^2}{4}}.
	\end{eqnarray*}

\end{proof}

In order to give a characterization of the true polyanalytic RBF space in terms of the operator $\mathcal{A}_{z, \bar{z}}$ we need to study its adjoint.

\begin{theorem}
	\label{adj1}
The adjoint of the operator $\mathcal{A}_{z, \bar{z}}$ on $L^2\left(e^{\frac{(z-\bar{z})^2}{4}}, \mathbb{C}\right)$ is given by
	\begin{equation}
		\label{adj}
		\mathcal{A}_{z, \bar{z}}^{*}= \frac{d}{d \bar{z}}.
	\end{equation}
	Moreover we have
	\begin{equation}
		\label{app}
		\mathcal{A}_{z, \bar{z}} \mathcal{A}_{z, \bar{z}}^{*}=- \frac{d^2}{dz d \bar{z}}-\frac{z}{2}\frac{d}{d \bar{z}}+ \frac{\bar{z}}{2} \frac{d}{d \bar{z}}.
	\end{equation}
\end{theorem}

\begin{proof}
	Let $u$, $v \in L^2\left(e^{\frac{(z-\bar{z})^2}{4}}, \mathbb{C}\right)$. By integrating by parts we have
	\begingroup\allowdisplaybreaks
	\begin{eqnarray*}
		&& \int_{\mathbb{C}}[\mathcal{A}_{z, \bar{z}}u(z)] \overline{v}(z)e^{\frac{(z-\bar{z})^2}{4}}d\lambda(z)\\
		&&=\int_{\mathbb{C}} \left[\left(-\frac{d}{dz}-\frac{M_z}{2}+ \frac{M_{\bar{z}}}{2}\right)u(z)\right]\overline{v(z)} e^{\frac{(z-\bar{z})^2}{4}}d\lambda(z)\\
		&&=\int_{\mathbb{R}^2} \left[\left(-\frac{1}{2}\frac{\partial}{\partial x}+\frac{i}{2} \frac{\partial}{\partial y}-iy\right)u(x,y)\right]\overline{v(x,y)} e^{-y^2}dxdy\\
		&&= -\frac{1}{2}\int_{\mathbb{R}^2} \left(\frac{\partial}{\partial x} u(x,y)\right)\overline{v(x,y)} e^{-y^2}dxdy+\frac{i}{2} \int_{\mathbb{R}^2} \left(\frac{\partial}{\partial y} u(x,y)\right)\overline{v(x,y)} e^{-y^2}dxdy\\
		&&-i \int_{\mathbb{R}^2} u(x,y) \overline{v(x,y)} y e^{- y^2}dxdy\\
		&&= \frac{1}{2}\int_{\mathbb{R}^2} u(x,y) \frac{\partial}{\partial x}[\overline{v(x,y)}] e^{- y^2}dxdy-\frac{i}{2} \int_{\mathbb{R}^2} u(x,y) \frac{\partial}{\partial y}\left[\overline{v(x,y)} e^{-y^2}\right]dxdy\\
		&&-i \int_{\mathbb{R}^2} u(x,y) \overline{v(x,y)} y e^{-y^2}dxdy\\
		&&= \frac{1}{2}\int_{\mathbb{R}^2} u(x,y) \frac{\partial}{\partial x}[\overline{v(x,y)}] e^{-y^2}dxdy-\frac{i}{2} \int_{\mathbb{R}^2} u(x,y) \frac{\partial}{\partial y}\left[\overline{v(x,y)}\right] e^{-y^2}dxdy\\
		&&+i\int_{\mathbb{R}^2} u(x,y) \overline{v(x,y)} y e^{-y^2}dxdy-i \int_{\mathbb{R}^2} u(x,y) \overline{v(x,y)} y e^{-y^2}dxdy\\
		&&= \int_{\mathbb{R}^2} u(x,y) \frac{1}{2} \left(\overline{\frac{\partial}{\partial x}v(x,y)+i\frac{\partial}{\partial y}v(x,y)}\right)e^{-y^2}dxdy\\
		&&= \int_{\mathbb{C}} u(z) \overline{\frac{d}{d \bar{z}}v(z)} e^{\frac{(z-\bar{z})^2}{4}}d\lambda(z)\\
		&&=\int_{\mathbb{C}} u(z)\overline{[\mathcal{A}_{z, \bar{z}}^* v(z)]} d \lambda(z).
	\end{eqnarray*}
	\endgroup
	
	This proves \eqref{adj}. Finally formula \eqref{app} follows by combining \eqref{adj} and \eqref{opem}.
	
\end{proof}

Due to the importance that the operator $\mathcal{A}_{z, \bar{z}} \mathcal{A}_{z, \bar{z}}^{*}$ will play in the sequel of the paper we establish the following notation:
$$ \widetilde{\Box}:=\mathcal{A}_{z, \bar{z}} \mathcal{A}_{z, \bar{z}}^{*}, \quad z \in \mathbb{C}.$$

The goal of the remaining part of this section is to study the main properties of the operator $\widetilde{\Box} $, which will be of fundamental importance for characterizing the true polyanalytic RBF space $\mathcal{H}_{T, 2}^{(\ell)}(\mathbb{C})$. In the sequel the following notation will be useful:
	$$  \widetilde{H}_{\ell,m}(z, \bar{z}):=H_{\ell-1,m}^{1/2}(z, \bar{z}) e^{-\frac{z^2}{4}}.$$

\begin{proposition}
	\label{eigenvalue}
	Let $m \in \mathbb{N}_0$, $\ell \in \mathbb{N}$ and $z \in \mathbb{C}$. Then $\widetilde{H}_{n,m}(z, \bar{z})$ are the eigenfunctions of the operator $\widetilde{\Box}$, i.e.:
	$$ \widetilde{\Box} [\widetilde{H}_{\ell,m}(z, \bar{z})]= \frac{\ell-1}{2} \widetilde{H}_{\ell,m}(z, \bar{z}).$$	
\end{proposition}
\begin{proof}
We start considering the case $\ell=1$. In this case by the definition of the complex Hermite polynomials, see \eqref{hermcom}, we have
\begin{eqnarray*}
\widetilde{\Box}\left(\widetilde{H}_{1,m}(z, \bar{z})\right)&=&\widetilde{\Box}\left(H_{0,m}^{1/2}(z, \bar{z}) e^{-\frac{z^2}{4}} \right)\\
&=&  \left(- \frac{d^2}{dz d \bar{z}}-\frac{z}{2}\frac{d}{d \bar{z}}+ \frac{\bar{z}}{2} \frac{d}{d \bar{z}} \right)\left[\left(\frac{z}{2}\right)^me^{- \frac{z^2}{4}}\right]\\
&=&0,
\end{eqnarray*}
the last equality follows from the fact that the function $\left(\frac{z}{2}\right)^me^{- \frac{z^2}{4}}$ is holomorphic. This proves the result for $\ell=1$. 
\\Now, we suppose that $\ell \geq 2$. We apply first the operator $\mathcal{A}_{z, \bar{z}}^{*}$ to the function $\widetilde{H}_{\ell,m}(z, \bar{z})$. By Theorem \ref{adj1} and \eqref{dabr} we get
	$$	\mathcal{A}_{z, \bar{z}}^{*} [\widetilde{H}_{\ell,m} (z, \bar{z})]=\frac{d}{d\bar{z}} \left(H_{\ell-1,m}^{1/2}(z, \bar{z})\right)e^{-\frac{z^2}{4}}=\frac{\ell-1}{2} H_{\ell-2,m}^{1/2}(z, \bar{z})e^{-\frac{z^2}{4}}.$$
	Now, we apply the operator $\mathcal{A}_{z, \bar{z}}$ to $\mathcal{A}_{z, \bar{z}}^{*} \widetilde{H}_{\ell,m}(z, \bar{z})$. By applying two times formula \eqref{form} we obtain
\begingroup\allowdisplaybreaks
	\begin{eqnarray*}
		\mathcal{A}_{z, \bar{z}} \mathcal{A}_{z, \bar{z}}^{*} [\widetilde{H}_{\ell,m}(z, \bar{z})]&=&\frac{\ell-1}{2} \mathcal{A}_{z, \bar{z}}\left(H_{\ell-2,m}^{1/2}(z, \bar{z}) e^{- \frac{z^2}{2}}\right)\\
		&=&\frac{\ell-1}{2^{m+1}} \mathcal{A}_{z, \bar{z}} \mathcal{A}_{z, \bar{z}}^{\ell-2} \left(z^m e^{-\frac{z^2}{4}}\right)\\
		&=& \frac{\ell-1}{2^{m+1}} \mathcal{A}_{z, \bar{z}}^{\ell-1} \left(z^m e^{-\frac{z^2}{4}}\right)\\
		&=& \frac{\ell-1}{2} H_{\ell-1,m}^{1/2}(z, \bar{z}) e^{-\frac{z^2}{4}}\\
		&=&\frac{\ell-1}{2}  \widetilde{H}_{\ell,m}(z, \bar{z}).
	\end{eqnarray*}
	\endgroup
\end{proof}

It is interesting to observe that the operator $\widetilde{\Box}$ can be related to the following type of Schrödinger operator:
\begin{definition}
 We define the operator $\mathcal{L}$ by
	\begin{equation}
		\label{soperator}
		\mathcal{L}:=-\frac{1}{4}\Delta-  \frac{iy}{2} \frac{\partial}{\partial x}+\frac{y^2}{4}, \quad  \Delta:= \frac{\partial^2}{\partial x^2}+\frac{\partial^2}{\partial y^2}.
	\end{equation}
\end{definition}

\begin{remark}
From the paper of Sondheimer and Wilson, see \cite{SW}, we recall that the hamiltonian of a free electron in constant magnetic field is given by
\begin{equation}\label{opee}
\mathcal{H}= - \frac{h^2}{8 \pi^2 m} \Delta+ \frac{iqh}{2 \pi  mc } \textbf{A} \cdot \nabla+ \frac{q^2}{2mc^2}\textbf{A}^2,
\end{equation}
where $\textbf{A}= \frac{1}{2} \textbf{H} \times \textbf{r}$ (with $\textbf{r}=(x,y,z)$) is the vector potential, $-q$ is the charge and $m$ is the mass of the electron and the remaining symbols have the usual meaning. If we take the magnetic field along the $z$ direction, then the vector potential becomes $\textbf{A}=\left(\frac{1}{2}Hy,0,0\right)$, and so the operator in \eqref{opee} can be written as
\begin{equation}
	\label{ope}
\mathcal{H}=- \frac{h^2}{8 \pi^2 m} \Delta- \frac{iqhH}{4 \pi  mc } y \frac{\partial}{\partial x}+ \frac{q^2H^2}{8m c^2} y^2.
\end{equation}
By taking $c=m=H=1$, $h=\sqrt{2} \pi$ and $q= \sqrt{2}$ in \eqref{ope} we get exactly the operator defined in \eqref{soperator}. Similar type of operators have been considered in \cite{CGS} to study evolution of superoscillations.
\end{remark}

In order to study the eigenfunctions and eigenvalues of $\mathcal{L}$, we first need to establish a connection with the operator \( \widetilde{\Box} \). To achieve this, we introduce the notion of the \emph{RBF ground state transformation}, which is motivated by the measure associated with the RBF space (see Definition~\ref{RBF}).

\begin{definition}
	Consider the Gaussian function
	\begin{equation}
		\label{functiong}
		g(y)= \frac{e^{-\frac{y^2}{2}}}{\sqrt{2\pi}}, \quad \text{ for all }y\in \mathbb R.
	\end{equation}
	We define the RBF-ground state transformation as the operator 
\begin{equation}
\label{GT}
G\colon L^2(\mathbb{R}^2, dxdy)\to L^2(\mathbb{R}^2, g^2 dxdy), \qquad f\mapsto Gf:= \frac{f}{g}.
\end{equation}
	
\end{definition}

\begin{lemma}
	\label{prop}
The operator $G$ is an isometric-isomorphism, which inverse is given by
$$G^{-1}\colon L^2(\mathbb{R}^2, g^2 dxdy) \to L^2(\mathbb{R}^2, dxdy), \qquad f\mapsto f \cdot g.$$
\end{lemma}

Now, we have all the tools to connect the operators $\mathcal{L}$, $\widetilde{\Box}$ and the RBF-ground state transformation. For convenience in the following proof we denote by $\partial_x^2$ and $\partial_y^2$ the derivatives $\frac{\partial^2}{\partial x^2}$ and $\frac{\partial^2}{\partial y^2}$.

\begin{theorem}
The following connection between the RBF-ground state transform, the operator $\mathcal{L}$ and $\widetilde{\Box}$ holds
	\begin{equation}
		\label{relope}
		\widetilde{\Box}=G\circ \left(\mathcal{L}-\frac{1}{4}\right) \circ G^{-1}.
	\end{equation}
\end{theorem}
\begin{proof}
	Let $u \in L^2(\mathbb{R}^2,dxdy)$. By \eqref{app} and the RBF-ground state transformation we have
	\begingroup\allowdisplaybreaks
	\begin{eqnarray}
		\nonumber
		\widetilde{\Box} (Tu)&=& \widetilde{\Box} \left(\frac{u}{g}\right)\\
		\nonumber
		&=& \sqrt{2\pi} \widetilde{\Box}  \left[ u(x,y)e^{\frac{y^2}{2}}\right]\\
		\nonumber
		&=&\sqrt{2\pi} \left[-\frac{1}{4} \partial_x^2 -\frac{1}{4} \partial_y^2 -\frac{iy}{2} \left(\partial_x+i \partial_y \right)\right]\left( u(x,y)e^{\frac{y^2}{2}}\right)\\
		\nonumber
		&=& \sqrt{2\pi} \left[- \frac{1}{4} (\partial_x^2 u(x,y)) e^{\frac{y^2}{2}}-\frac{1}{4} \partial_y \left((\partial_yu(x,y))e^{\frac{y^2}{2}}+ yu(x,y) e^{\frac{y^2}{2}}\right) \right.\\
		\nonumber
		&&\left. -\frac{iy}{2} \left( (\partial_x u(x,y))e^{\frac{y^2}{2}}+i (\partial_y u(x,y)) e^{\frac{y^2}{2}}+i yu(x,y) e^{\frac{y^2}{2}}\right) \right]\\
		\nonumber
		&=&\sqrt{2\pi} \left[- \frac{1}{4} \partial_x^2 u(x,y)-\frac{1}{4} \partial_y^2u(x,y)-\frac{y}{2} \partial_y u(x,y)-\frac{u(x,y)}{4}\right.\\
		\nonumber
		&&\left. -\frac{y^2}{4}u(x,y)-\frac{i y}{2} \partial_x u(x,y)+\frac{y}{2} \partial_y u(x,y)+ \frac{y^2}{2}u(x,y) \right]e^{\frac{y^2}{2}}\\
		\label{point1}
		&=& \sqrt{2\pi} \left[- \frac{1}{4} \Delta + \frac{y^2}{4} u(x,y)-\frac{i y}{2} \partial_x u(x,y)-\frac{1}{4}u(x,y) \right] e^{\frac{y^2}{2}}.
	\end{eqnarray}
	\endgroup
	Now by the definition of the operator $\mathcal{L}$, see \eqref{soperator} and the RBF-ground state transformation (see \eqref{GT}) we have
	\begin{eqnarray}
		\nonumber
		G \left(\mathcal{L}-\frac{1}{2}\right)u&=& G \left[- \frac{1}{4} \Delta + \frac{y^2}{4} u(x,y)-\frac{i y}{2} \partial_x u(x,y)-\frac{1}{4}u(x,y) \right]\\
		\label{point2}
		&=&\sqrt{2\pi} \left[- \frac{1}{4} \Delta + \frac{y^2}{4} u(x,y)-\frac{i y}{2} \partial_x u(x,y)-\frac{1}{4}u(x,y) \right]e^{\frac{y^2}{2}}.
	\end{eqnarray}
	Since \eqref{point1} and \eqref{point2} coincide, we have

	$$\widetilde{\Box} \circ G =G \circ \left(\mathcal{L}-\frac{1}{4}\right).$$
	
Hence, by applying the operator $ G^{-1}$ on the left-hand side of the above equation, we obtain equation~\eqref{relope}.

\end{proof}
The above result is crucial to figure out the eigenvalues and eigenfunctions of the operator $\mathcal{L}$.

\begin{theorem}
 The operator $ \mathcal{L}$ has eigenvalues $\frac{1}{2}\left(\ell-1\right)$, with $\ell \in \mathbb{N}$, and corresponding eigenfunctions are $g(y) \widetilde{H}_{\ell,m}(z,\bar{z})$, with $m \in \mathbb{N}_0$.
\end{theorem}
\begin{proof}
	By Proposition \ref{eigenvalue} we know that $\widetilde{H}_{\ell,m}(z, \bar{z})$ are eigenfunctions of the operator $\widetilde{\Box}$. So by \eqref{relope} we have
	$$ G \left(\mathcal{L}-\frac{1}{4}\right) G^{-1} \left(\widetilde{H}_{\ell,m}(z, \bar{z})\right)=\frac{\ell-1}{2}\widetilde{H}_{\ell,m}(z, \bar{z}).$$
	We apply on the left the operator $G^{-1}$, and we get
	$$  \left(\mathcal{L}-\frac{1}{4}\right) G^{-1} \left(\widetilde{H}_{\ell,m}(z, \bar{z})\right)=\frac{\ell-1}{2} G^{-1}\left(\widetilde{H}_{\ell,m}(z, \bar{z})\right).$$
	By Lemma \ref{prop} we have
	$$ \left( \mathcal{L}-\frac{1}{4}\right) (g(y)\widetilde{H}_{\ell,m}(z, \bar{z}))=\frac{\ell-1}{2} g(y)\widetilde{H}_{\ell,m}(z, \bar{z}).$$
	Thus we have
	$$ \mathcal{L}[g(y)\widetilde{H}_{\ell,m}(z, \bar{z})]=\frac{1}{2}\left(\ell-1 \right)\left(g(y)\widetilde{H}_{\ell,m}(z, \bar{z})\right).$$
	
	This proves the result.
\end{proof}

Now, we provide a specific definition for the eigenspaces of the operator $\widetilde{\Box}$.

\begin{definition}
	Let $\ell \in \mathbb{N}$. We define the RBF Landau level as
	$$ \mathbf{H}^\ell(\mathbb{C}):= \biggl \{f\in L^{2}\left(\mathbb{C}, e^{\frac{(z-\bar{z})^2}{4}}\right), \qquad 	\widetilde{\Box}(f(z))=\frac{\ell-1}{2}  f(z) \biggl\}.$$
\end{definition}

The above space  have a deep connection with the true polyanalytic RBF space.

\begin{theorem}
	Let $\ell \in \mathbb{N}$. Then we have that the true RBF space and the space $\mathbf{H}^\ell(\mathbb{C})$ coincide i.e.: 
	\begin{equation}
		\label{impeq}
		\mathbf{H}^\ell(\mathbb{C})=\mathcal{H}_{T}^{(\ell)}(\mathbb{C}).
	\end{equation}
\end{theorem}
\begin{proof}
	We prove the result by double inclusion. Let $f \in \mathbf{H}^\ell(\mathbb{C})$. By definition it is clear that $f \in L^2\left(\mathbb{C}, e^{\frac{(z-\bar{z})^2}{4}}\right)$.
	So by Proposition \ref{l2iso} there exists $g \in L^2\left(\mathbb{C},  e^{-\frac{|z|^2}{2}}\right)$ such that
	\begin{equation}
		\label{fun}
		f(z)=g(z) e^{-\frac{z^2}{4}}.
	\end{equation}
By hypothesis we have $ \widetilde{\Box} f(z) = \frac{\ell-1}{2} f(z)$, and recalling the definition of the magnetic Laplacian see \eqref{magnetic}, it follows that

	\begin{eqnarray*}
		\widetilde{\Delta} g(z)&=& \left(-\frac{d}{dz} \frac{d}{d\bar{z}}+\frac{\bar{z}}{2} \frac{d}{d\bar{z}}\right)(f(z)e^{ \frac{z^2}{4}})\\
		&=&- \frac{d}{d\bar{z}} \left(\frac{d}{dz} f(z)+ \frac{z}{2} f(z)\right)e^{\frac{z^2}{4}}+ \frac{\bar{z}}{2} \left(\frac{d}{d\bar{z}}f(z)\right)e^{\frac{z^2}{4}}\\
		&=& \left(- \frac{d}{d\bar{z}} \frac{d}{dz} f(z)- \frac{z}{2} \frac{d}{d\bar{z}}f(z)+\frac{\bar{z}}{2} \frac{d}{d\bar{z}}f(z)\right) e^{\frac{z^2}{4}}\\
		&=& [\widetilde{\Box}f(z)]e^{\frac{z^2}{4}}\\
		&=& \frac{\ell-1}{2} f(z)e^{\frac{z^2}{4}}\\
		&=&\frac{\ell-1}{2} g(z).
	\end{eqnarray*}
	Since $g \in L^2 \left(\mathbb{C}, e^{-\frac{|z|^2}{2}}\right)$ and  $\widetilde{\Delta} g(z)=\frac{\ell-1}{2} g(z)$ we get that $ g \in \mathcal{A}_{\ell, \frac{1}{2}}(\mathbb{C})$, see \eqref{mouy}. By \eqref{mouy1} we have that $g \in \mathcal{F}_{T, \frac{1}{2}}^{(\ell)}(\mathbb{C})$. So by Proposition \ref{iso} and \eqref{fun} we get that $f \in \mathcal{H}_{T}^{(\ell)}(\mathbb{C})$. Now, we prove the other inclusion. Let $h \in \mathcal{H}_{T}^{(\ell)}(\mathbb{C})$, so by definition $h \in L^2\left(\mathbb{C}, e^{\frac{(z-\bar{z})^2}{4}}\right)$. By Proposition \ref{iso} we know that there exist a function $p \in \mathcal{F}^{(\ell)}_{T, \frac{1}{2}}(\mathbb{C})$ such that
	$$ h(z)= p(z) e^{-\frac{z^2}{4}}.$$
	Since $p  \in  \mathcal{A}_{\ell, \frac{1}{2}}(\mathbb{C})$, see \eqref{mouy1}, we get that $ \widetilde{\Delta}p(z)=\frac{\ell-1}{2}p(z)$. So, we have 
	\begin{eqnarray*}
		\frac{\ell-1}{2}h(z)&=& \frac{\ell-1}{2}p(z) e^{- \frac{z^2}{4}}\\
		&=& \left(\widetilde{\Delta}p(z)\right)e^{- \frac{z^2}{4}}\\
		&=&- \frac{d^2}{dz d \bar{z}} (p(z)) e^{- \frac{z^2}{4}}+ \frac{z}{2} \left(\frac{d}{d \bar{z}} p(z)\right) e^{-\frac{z^2}{4}}+ \frac{\bar{z}}{2} \left(\frac{d}{d \bar{z}} g(z)\right) e^{-\frac{z^2}{4}}- \frac{z}{2} \left(\frac{d}{d\bar{z}} p(z)\right) e^{-\frac{z^2}{4}}\\
		&=&- \frac{d}{d\bar{z}}\left[\left(\frac{d}{dz} p(z)\right) - \frac{z}{2} p(z)\right] e^{- \frac{z^2}{2}}+\frac{\bar{z}}{2} \left(\frac{d}{d \bar{z}} p(z)\right) e^{-\frac{z^2}{4}}- \frac{z}{2} \left(\frac{d}{d\bar{z}} p(z)\right) e^{-\frac{z^2}{4}}\\
		&=&- \frac{d^2}{dz d \bar{z}} \left( p(z) e^{- \frac{z^2}{4}}\right)+ \frac{\bar{z}}{2}\frac{d}{d \bar{z}} \left( p(z)e^{-\frac{z^2}{4}}\right) - \frac{z}{2} \frac{d}{d\bar{z}}\left( p(z)e^{-\frac{z^2}{4}}\right) \\
		&=& \left(- \frac{d^2}{dz d \bar{z}}+\frac{\bar{z}}{2}\frac{d}{d \bar{z}}- \frac{z}{2} \frac{d}{d\bar{z}}\right)\left( p(z)e^{-\frac{z^2}{4}}\right)\\
		&=&\widetilde{\Box} (h(z))
	\end{eqnarray*}
	This implies that by definition we have $h \in \mathbf{H}^\ell(\mathbb{C})$.
\end{proof}

The above result implies a new characterization of the analytic Gaussian RBF space.
\begin{corollary}
The analytic Gaussian RBF space can be characterized as follows
	$$ \mathcal{H}_{2}(\mathbb{C}):= \biggl \{f\in L^{2}\left(\mathbb{C}, e^{\frac{(z-\bar{z})^2}{4}}\right), \qquad 	\widetilde{\Box}(f(z))=0\biggl\}.$$
	
\end{corollary}

\begin{proof}
	The result follows by taking $\ell=1$ in \eqref{impeq}.
\end{proof}

 \section{Weyl operators for polyanalytic functions}
\setcounter{equation}{0}
In this section, we study  the Weyl operator in the polyanalytic setting. While such operator has already been considered in~\cite{FH}, the key novelty of our approach lies in providing an explicit construction. To this end, we analyze compositions of the true polyanalytic Bargmann transform $B_n $, its inverse, and time-frequency shift operators. These results are then applied to derive new formulas involving the Christoffel-Darboux formula and Mehler's kernel. Finally, in the last part of the section, we present a construction of a Weyl operator based on the true polyanalytic RBF space.

\subsection{Polyanalytic Weyl operators and time-frequency shifts}
First, we introduce the \textit{true polyanalytic Weyl operator} with a real parameter $a$ as the composition of the true polyanalytic Bargmann transform, its inverse, and the translation operator $\tau_a$:
\begin{definition}[Polyanalytic Weyl operator, Step I]
\label{pwI}
Let $a\in \mathbb{R}$ and $\ell\in\mathbb N$. 

The true polyanalytic Weyl operator of order $\ell$, denoted by $W^{a}_{\ell}$, is the operator making the following diagram commutative:
$$
W^{a}_{\ell}: \xymatrix{
	\mathcal{F}^{(\ell)}_T(\mathbb{C})\ar[r] \ar[d]_{B_{\ell}^{-1}} & \mathcal{F}^{(\ell)}_T(\mathbb{C})  \\ L^2(\mathbb{R}) \ar[r]_{\tau_a} & L^2(\mathbb{R}) \ar[u]_{B_\ell}
}$$
In other words,
$$W_{\ell}^a:=B_\ell \tau_a B^{-1}_{\ell},$$
where $B_\ell$ is the true polyanalytic Bargmann transform of order $\ell$, see \eqref{truepoly}
\end{definition}
\begin{lemma}\label{L1}
Let $\ell \in \mathbb N$, $a\in \mathbb R$,  $f\in \mathcal{F}^{(\ell)}_T(\mathbb C),$ and $z\in \mathbb C$. Then

$$ (W^a_\ell f)(z)=\int_{\mathbb{C}}f(w) I_{a,\ell-1}(z,w)e^{-|w|^2}\frac{d\lambda(w)}{\pi},$$

where  $$I_{a,\ell-1}(z,w):=\int_\mathbb{R}H_{\ell-1}\left(\frac{z+\overline{z}}{\sqrt{2}}-x\right) H_{\ell-1}\left(\frac{w+\overline{w}}{\sqrt{2}}-x+a\right) A_z(x) \overline{A_w(x-a)} dx,\quad z,w\in \mathbb{C},$$
where $A_z$ is the Segal Bargman kernel introduced in Definition \ref{SBkernel}.
\end{lemma}
\begin{proof}
For every $f\in \mathcal{F}^{(\ell)}_T(\mathbb C)$ by using \eqref{inverse}, define the function $\psi_f$ in $L^2(\mathbb{R})$
\begin{align*}
 \psi_f(x)&=\tau_aB^{-1}_{\ell}f(x)\\
&=B^{-1}_{\ell}f(x-a)\\
&=\int_\mathbb{C}H_{\ell-1}\left(\frac{w+\overline{w}}{\sqrt{2}}-x+a\right)\overline{A_w(x-a)}f(w)e^{-|w|^2}\frac{d\lambda(w)}{\pi}
\end{align*}
By the definition of the true Bargmann transform, see \eqref{truepoly}, we have
\begin{align*}
\displaystyle (W^a_\ell f)(z)&=(B_\ell\psi_f)(z)\\
&=\int_\mathbb{R}H_{\ell-1}\left(\frac{z+\overline{z}}{\sqrt{2}}-x\right)A_z(x)\psi_f(x) dx\\
&=\int_\mathbb{R}H_{\ell-1}\left(\frac{z+\overline{z}}{\sqrt{2}}-x\right)A_z(x)\left(\int_\mathbb{C}H_{\ell-1}\left(\frac{w+\overline{w}}{\sqrt{2}}-x+a\right)\overline{A_w(x-a)}f(w)e^{-|w|^2}\frac{d\lambda(w)}{\pi}\right) dx
\end{align*}
Applying Fubini's theorem to exchange the order of integration, we get

$$
\displaystyle  (W^a_\ell f)(z)=\int_{\mathbb C} f(w)\left(\int_\mathbb{R}H_{\ell-1}\left(\frac{z+\overline{z}}{\sqrt{2}}-x\right) H_{\ell-1}\left(\frac{w+\overline{w}}{\sqrt{2}}-x+a\right) A_z(x) \overline{A_w(x-a)} dx\right) e^{-|w|^2}\frac{d\lambda(w)}{\pi}.
$$

This simplifies to 

$$
\displaystyle  (W^a_\ell f)(z)=\int_{\mathbb C} f(w) I_{a,\ell-1}(z,w) e^{-|w|^2}\frac{d\lambda(w)}{\pi},
$$

which completes the proof.
\end{proof}

Now, the main  problem is to compute the following integral function.
\begin{definition}
Let $m\in \mathbb N_0$ and $a\in \mathbb R$. Define the integral function $I_{a,m}(z,w)$ as
\begin{equation}
I_{a,m}(z,w):= \int_\mathbb{R}H_m\left(\frac{z+\overline{z}}{\sqrt{2}}-x\right) H_m\left(\frac{w+\overline{w}}{\sqrt{2}}-x+a\right) A_z(x) \overline{A_w(x-a)} dx,
\end{equation}
for every $z,w\in \mathbb C.$
\end{definition}

\begin{lemma}\label{L2} Let $a\in\mathbb{R}$ and  $m\in \mathbb N_0$.Then

$$I_{a,m}(z,w)= m! 2^m e^{-\frac{a^2}{4}+\frac{a}{\sqrt{2}}z}e^{(z-\frac{a}{\sqrt{2}})\overline{w}} L^0_m\left(\left|z-\frac{a}{\sqrt{2}}-w \right|^2\right),$$
for every $z,w\in \mathbb{C}$.
\end{lemma}
\begin{proof}
By the definition of the Segal Bargmann kernel, see Definition \ref{SBkernel}, we have
$$A_z(t) \overline{A_w(t-a)}=\frac{1}{\sqrt{\pi}}e^{-\frac{1}{2}(z^2+t^2)+\sqrt{2}zt}   e^{-\frac{1}{2}(\overline{w}^2+(t-a)^2)+\sqrt{2}\overline{w}(t-a)}.$$
Set $$\displaystyle b=\frac{z+\overline{z}}{\sqrt{2}} \text{ and } c=\frac{w+\overline{w}}{\sqrt{2}}+a.$$ Then
\begin{eqnarray}
\nonumber I_{a,m}(z,w)&=&\int_\mathbb{R}H_m\left(\frac{z+\overline{z}}{\sqrt{2}}-t\right) H_m\left(\frac{w+\overline{w}}{\sqrt{2}}+a-t\right) A_z(t) \overline{A_w(t-a)} dt\\
\nonumber
&=&\frac{1}{\sqrt{\pi}}e^{-\frac{1}{2}(z^2+\overline{w}^2+a^2)-\sqrt{2}\overline{w} a}\int_\mathbb{R} e^{-t^2+t\left[a+\sqrt{2}(z+\overline{w})\right]}H_m\left(b-t\right) H_m\left(c-t\right)dt.\\
\label{F1}
&=&\frac{1}{\sqrt{\pi}}e^{-\frac{1}{2}(z^2+\overline{w}^2+a^2)-\sqrt{2}\overline{w} a} g_{a,m}(z,w)
%&=e^{-\frac{1}{2}(z^2+\overline{w}^2+a^2)-\sqrt{2}\overline{w} a} g_{a,m}(z,w),\\
\end{eqnarray}

where we define 
\begin{equation}
 g_{a,m}(z,w):=\int_\mathbb{R} e^{-t^2+t\left[a+\sqrt{2}(z+\overline{w})\right]}H_m\left(b-t\right) H_m\left(c-t\right)dt.
 \end{equation}
 By the first equality of Lemma \ref{ACHA}, by substituting the parameters $b,c,x,u, \beta$ as follows

$$\displaystyle b=x=\frac{z+\overline{z}}{\sqrt{2}}, \quad  c=u=\frac{w+\overline{w}}{\sqrt{2}}+a, \quad \lambda=-\frac{i}{\sqrt{2}}(z-w+\overline{(w-z)}),$$

we get

\begin{align*}
 \displaystyle g_{a,m}(z,w)&=\int_\mathbb{R} e^{-t^2+t\left[a+\sqrt{2}(z+\overline{w})\right]}H_m\left(b-t\right) H_m\left(c-t\right)dt\\
 &=(-1)^m e^{\frac{b^2}{2}+\frac{c^2}{2}} (M_bh_m*M_ch_m)(\lambda).
\end{align*}
The second equality of Lemma  \ref{ACHA} allows us to express $g_{a,m}(z,w)$ in terms of the complex Ito-Hermite polynomials leading to 
\begin{equation}
 \displaystyle g_{a,m}(z,w)=(-1)^m \sqrt{\pi} 2^{m} e^{\frac{b^2}{2}+\frac{c^2}{2}} e^{- \frac{\lambda^2}{4}+\frac{i \lambda (b+c)}{2}}   e^{-\frac{(b-c)^2}{4}} H_{m,m}\left(\frac{c-b+i\lambda}{\sqrt{2}},\frac{c-b-i\lambda}{\sqrt{2}}\right).
\end{equation}
By \eqref{relHL} we have

$$ H_{m,m}\left(\frac{c-b+i\lambda}{\sqrt{2}},\frac{c-b-i\lambda}{\sqrt{2}}\right)=(-1)^m m!L^0_m\left(\left| \frac{c-b+i\lambda}{\sqrt{2}} \right|^2\right),\quad m\in \mathbb N.$$

Observe that  
$$\left| \frac{c-b+i\lambda}{\sqrt{2}} \right|^2=\frac{1}{2}(c-b+i\lambda)(c-b-i\lambda)=\frac{1}{2}[\lambda^2+(b-c)^2].$$
We can easily check that 
$$\lambda^2=-\frac{1}{2} \left[ (z-w)^2+\overline{(z-w)}^2-2|z-w|^2\right],$$ 
and  $$(b-c)^2=\frac{1}{2}\left[ 2a^2+(w-z)^2+\overline{(w-z)}^2+2|w-z|^2+2\sqrt{2}a(w-z) +2\sqrt{2}a\overline{(w-z)} \right].$$
Thus, it follows that
$$\lambda^2+(b-c)^2=a^2+2|z-w|^2+\sqrt{2}a\left(w-z+\overline{(w-z)}\right).$$
The computations above yield
$$ \left| \frac{c-b+i\lambda}{\sqrt{2}} \right|^2=\frac{a^2}{2}+|z-w|^2+\frac{a}{\sqrt{2}}\left[w-z+\overline{(w-z)}\right]=\left|z-\frac{a}{\sqrt{2}}-w \right|^2.$$

As a direct consequence of \eqref{relHL}, we have  
\begin{align*} 
\displaystyle H_{m,m}\left(\frac{c-b+i\lambda}{\sqrt{2}},\frac{c-b-i\lambda}{\sqrt{2}}\right)&=(-1)^m m!L^0_m\left(\left|z-\frac{a}{\sqrt{2}}-w \right|^2\right).
\end{align*}

Hence

\begin{align*}
g_{a,m}(z,w)&=(-1)^m \sqrt{\pi} 2^{m} e^{\frac{b^2}{2}+\frac{c^2}{2}} e^{- \frac{\lambda^2}{4}+\frac{i \lambda (b+c)}{2}}   e^{-\frac{(b-c)^2}{4}} H_{m,m}\left(\frac{c-b+i\lambda}{\sqrt{2}},\frac{c-b-i\lambda}{\sqrt{2}}\right)\\
&=\sqrt{\pi} m! 2^m e^{\frac{b^2}{2}+\frac{c^2}{2}} e^{- \frac{\lambda^2}{4}+\frac{i \lambda (b+c)}{2}}   e^{-\frac{(b-c)^2}{4}} L^0_m\left(\left|z-\frac{a}{\sqrt{2}}-w \right|^2\right) \\
&=\sqrt{\pi} m! 2^m e^{\frac{b^2}{4}+\frac{c^2}{4}+\frac{bc}{2}+\frac{i\lambda (b+c)}{2}-\frac{\lambda^2}{4}} L^0_m\left(\left|z-\frac{a}{\sqrt{2}}-w \right|^2\right)\\
&=\sqrt{\pi} m! 2^m e^{\frac{(b+c)^2}{4}+\frac{i\lambda(b+c)}{2} -\frac{\lambda^2}{4}}L^0_m\left(\left|z-\frac{a}{\sqrt{2}}-w \right|^2\right)\\
&=\sqrt{\pi} m! 2^m e^{\frac{(b+c+i\lambda)^2}{4}}L^0_m\left(\left|z-\frac{a}{\sqrt{2}}-w \right|^2\right).
\end{align*}
Note that $b+c+i\lambda=a+\sqrt{2}(z+\overline{w})$. This
leads to
\begin{align*}
\displaystyle  \displaystyle g_{a,m}(z,w)&=\sqrt{\pi} m! 2^m e^{\frac{\left[a+\sqrt{2}(z+\overline{w})\right]^2}{4}}L^0_m\left(\left|z-\frac{a}{\sqrt{2}}-w \right|^2\right)\\
&=\sqrt{\pi} m! 2^m e^{\frac{a^2}{4}+\frac{(z+\overline{w})^2}{2}+\frac{a}{\sqrt{2}}(z+\overline{w})}L^0_m\left(\left|z-\frac{a}{\sqrt{2}}-w \right|^2\right)\\
&=\sqrt{\pi} m! 2^m e^{\frac{a^2}{4}+\frac{z^2}{2}+\frac{\overline{w}^2}{2}+z\overline{w}+\frac{a}{\sqrt{2}}\left(z+\overline{w}\right) } L^0_m\left(\left|z-\frac{a}{\sqrt{2}}-w \right|^2\right).
\end{align*}
Finally, we substitute $g_{a,m}(z,w)$ with its expression from \eqref{F1} and get 

\begin{align*}
\displaystyle  I_{a,m}(z,w)
&= m! 2^m  e^{-\frac{1}{2}(z^2+\overline{w}^2+a^2)-\sqrt{2}\overline{w} a} e^{\frac{a^2}{4}+\frac{z^2}{2}+\frac{\overline{w}^2}{2}+z\overline{w}+\frac{a}{\sqrt{2}}\left(z+\overline{w}\right) } L^0_m\left(\left|z-\frac{a}{\sqrt{2}}-w \right|^2\right) \\
&=m! 2^m e^{-\frac{a^2}{4}+\frac{a}{\sqrt{2}}z}e^{(z-\frac{a}{\sqrt{2}})\overline{w}} L^0_m\left(\left|z-\frac{a}{\sqrt{2}}-w \right|^2\right).
\end{align*}

\end{proof}
We can easily rewrite Lemma \ref{L2} using the true-polyanalytic reproducing kernel $\mathsf{K}_{T,m}$ as follows:
\begin{corollary}\label{C1}
Let $a\in\mathbb{R}$ and $m\in\mathbb N_0$. Then
\begin{equation}
\label{binte}
I_{a,m}(z,w)=m! 2^m e^{-\frac{a^2}{4}+\frac{za}{\sqrt{2}}}\mathsf{K}_{m+1,T}\left(z-\frac{a}{\sqrt{2}},w \right),
\end{equation}
for every $z,w\in \mathbb C$.
\end{corollary}
\begin{proof}
By the definition of the true polyanalytic Fock kernel (see~\eqref{tpoly}) with $ \alpha = 1 $, and by Lemma~\ref{L2}, we have:

\begin{align*}
\displaystyle I_{a,m}(z,w)&= m! 2^m e^{-\frac{a^2}{4}+\frac{za}{\sqrt{2}}}e^{(z-\frac{a}{\sqrt{2}})\overline{w}} L^0_m\left(\left|z-\frac{a}{\sqrt{2}}-w \right|^2\right)\\
&=m! 2^m e^{-\frac{a^2}{4}+\frac{za}{\sqrt{2}}}\mathsf{K}_{m+1,T}\left(z-\frac{a}{\sqrt{2}},w \right). 
\end{align*}
\end{proof}
As a consequence of the previous results, we can easily prove:
\begin{theorem}\label{TrueWeylR1}
Let  $a\in \mathbb R$, $\ell\in\mathbb N$, and  $f\in \mathcal{F}^{(\ell)}_T(\mathbb{C})$. Then

$$(W^a_\ell f)(z)=(\ell-1)!2^{\ell-1}e^{-\frac{a^2}{4}+\frac{za}{\sqrt{2}}} f(z-\frac{a}{\sqrt{2}}), \quad \text{ for all } z\in \mathbb C.$$

\end{theorem}
\begin{proof}
By applying Lemma \ref{L1} and  Corollary \ref{C1}, we get
\begin{align*}
\displaystyle (W^a_\ell f)(z)
&=\int_{\mathbb{C}}I_{a,\ell-1}(z,w) f(w) e^{-|w|^2}\frac{d\lambda(w)}{\pi} \\
&=(\ell-1)!2^{\ell-1}e^{-\frac{a^2}{4}+\frac{za}{\sqrt{2}}}\int_{\mathbb{C}} \mathsf{K}_{\ell-1,T}\left(z-\frac{a}{\sqrt{2}},w \right) f(w) e^{-|w|^2}\frac{d\lambda(w)}{\pi}\\
%&=m!2^me^{-\frac{a^2}{4}+\frac{za}{\sqrt{2}}}\int_{\mathbb{C}} \overline{\mathsf{K}_{m,T}\left(w,z-\frac{a}{\sqrt{2}} \right)} f(w) e^{-|w|^2}\frac{dA(w)}{\pi}\\
&=(\ell-1)!2^{\ell-1}e^{-\frac{a^2}{4}+\frac{za}{\sqrt{2}}} \int_\mathbb C f(w) \overline{\mathsf{K}_{\ell-1,T,z-\frac{a}{\sqrt{2}}} (w)} e^{-|w|^2}\frac{d\lambda(w)}{\pi}.
\end{align*}
Therefore, using the reproducing kernel property of the true-polyanalytic Fock kernel $\mathsf{K}_{T,\ell-1}$, we obtain
\begin{align*}
\displaystyle (W_\ell^af)(z)&= (\ell-1)!2^{\ell-1}e^{-\frac{a^2}{4}+\frac{za}{\sqrt{2}}} \langle f, \mathsf{K}_{\ell-1,T,z-\frac{a}{\sqrt{2}}} \rangle  \\
&=(\ell-1)!2^{\ell-1}e^{-\frac{a^2}{4}+\frac{za}{\sqrt{2}}} f(z-\frac{a}{\sqrt{2}})\\
&=(\ell-1)!2^{\ell-1} k_{\frac{a}{\sqrt{2}}}(z)f(z-\frac{a}{\sqrt{2}}).
\end{align*}
\end{proof}
\begin{corollary}
Let $\ell\in \mathbb N,$ and $a,b\in \mathbb R$. Then, the following semi-group property holds
$$W_\ell^aW_\ell^b=(\ell-1)!2^{\ell-1}W_{\ell}^{a+b}.$$
\end{corollary}
\begin{proof}
This property of the true polyanalytic Weyl operator $W_\ell^a$ is a direct application of the Definition \ref{pwI}.
\end{proof}

\begin{remark}
By taking $ f = H_{m,n}(z, \bar{z})$ in Theorem~\ref{TrueWeylR1}, we can express the translation of the complex Hermite polynomials in terms of an integral involving Laguerre polynomials:

 $$H_{m,n}\left(z-\frac{a}{\sqrt{2}},  \overline{z}-\frac{a}{\sqrt{2}}\right)=\int_\mathbb C e^{(z-\frac{a}{\sqrt{2}})\overline{w}} H_{m,n}(w,\overline{w}) L^0_m\left(\left|z-\frac{a}{\sqrt{2}}-w \right|^2\right) e^{-|w|^2}\frac{d\lambda(w)}{\pi}.$$
\end{remark}

\begin{remark}
Let $\ell \in \mathbb N$. Observe that, for the particular case $a=0$, we have 
$$I_{0,\ell}(z,w)= \int_\mathbb{R}H_{\ell-1}\left(\frac{z+\overline{z}}{\sqrt{2}}-x\right) H_{\ell-1}\left(\frac{w+\overline{w}}{\sqrt{2}}-x\right) A_z(x) \overline{A_w(x)} dx=(\ell-1)! 2^{\ell-1} \mathsf{K}_{\ell,T}\left(z,w \right).$$
Hence, by considering the feature map $\phi_\ell:z \longmapsto \phi_{\ell,z}\in L^2(\mathbb{R})$ defined by 
\begin{equation}
\phi_{\ell,z}(x)=\frac{1}{\sqrt{2^{\ell-1}(\ell-1)!}}H_{\ell-1}\left(\frac{z+\overline{z}}{\sqrt{2}}-x\right) A_z(x),\quad x\in \mathbb R, 
\end{equation}
we obtain the following factorization of the true polyanalytic Fock kernel: 
$$ \mathsf{K}_{\ell, T}(z,w)=\langle  \phi_{\ell,z},  \phi_{\ell,w} \rangle_{L^2(\mathbb R)}.$$
This allows us to conclude that the true polyanalytic Fock kernel is positive definite.
\end{remark}
\begin{remark}
The case $\ell=1$ corresponds to the well-known result for the classical Fock space presented in \cite{Z1}. 
\end{remark}

Now, in order to introduce a Weyl operator with a complex parameter, we need to add an additional step in the construction of the Weyl operator by composing the true Bargmann transform, its inverse, and the modulation operator.

\begin{definition}[Polyanalytic Weyl operator, Step II]
Let $b\in \mathbb{R}$ and $\ell\in\mathbb N$. 

The true polyanalytic Weyl operator of order $m$, denoted by $W^{ib}_{\ell}$, is the operator making the following diagram commute $$
W^{ib}_{\ell}: \xymatrix{
    \mathcal{F}^{(\ell)}_T(\mathbb{C}) \ar[r] \ar[d]_{B_{\ell}^{-1}} & \mathcal{F}^{(\ell)}_T(\mathbb{C})  \\ L^2(\mathbb{R}) \ar[r]_{M_b} & L^2(\mathbb{R}) \ar[u]_{B_\ell}
  }$$
so that 
$$W_{\ell}^{ib}:=B_\ell M_b B^{-1}_{\ell},$$
where $B_\ell$ is the true polyanalytic Bargmann transform of order $\ell$.
\end{definition}

Similarly to the case of translation, the action of the operator $ W_{\ell}^{ib} $ on a function in the true Fock space can be expressed as follows:

\begin{theorem}\label{Modulation}
Let $\ell \in \mathbb N$, $b\in \mathbb R$, and  $f\in \mathcal{F}^{(\ell)}_T(\mathbb{C})$. Then

$$(W^{ib}_\ell f)(z)=(\ell-1)!2^{\ell-1} e^{-\frac{b^2}{4}-\frac{ibz}{\sqrt{2}}} f\left(z+\frac{ib}{\sqrt{2}}\right), \quad \text{ for all } z\in \mathbb C.$$

\end{theorem}
\begin{remark}
We omit the proof of Theorem \ref{Modulation}, since it is based on arguments similar to those used in the proof of Theorem \ref{TrueWeylR1}.
\end{remark}

Finally, we have all the necessary tools to introduce the true polyanalytic Weyl operator for a generic complex parameter.

\begin{definition}[Polyanalytic Weyl operator, Step III]
	\label{threeW}
Let $a,b\in \mathbb{R}$ and $\ell \in\mathbb N$. 

The true polyanalytic Weyl operator of order $\ell$, denoted by $W^{a+ib}_{\ell}$, is the operator making the following diagram commute $$
W^{a+ib}_{\ell}: \xymatrix{
    \mathcal{F}^{(\ell)}_T(\mathbb{C}) \ar[r] \ar[d]_{B_{\ell}^{-1}} & \mathcal{F}^{(\ell)}_T(\mathbb{C})  \\ L^2(\mathbb{R}) \ar[r]_{M_b\tau_a} & L^2(\mathbb{R}) \ar[u]_{B_\ell}
  }$$
so that 
$$W_{\ell}^{a+ib}:=B_\ell M_b\tau_a B^{-1}_{\ell},$$
where $B_\ell$ is the true polyanalytic Bargmann transform of order $\ell$.
\end{definition}

More generally, the following theorem holds

\begin{theorem}\label{TrueWeylR2}
Let $\ell\in\mathbb N$, $a,b\in \mathbb R$, and $f\in \mathcal{F}^{(\ell)}_T(\mathbb{C})$. Then

$$(W^{a+ib}_\ell f)(z)=(\ell-1)!2^{\ell-1} k_{1,\frac{a+ib}{\sqrt{2}}}(z) f\left(z- \frac{a-ib}{\sqrt{2}}\right), \quad \text{ for all } z\in \mathbb{C},$$
where $k_{1,\frac{a+ib}{\sqrt{2}}}(z)$ is the normalized Fock kernel.
\end{theorem}

\begin{remark}
\label{vectorial}
It is possible to consider the full polyanalytic Bargmann transform, see~\eqref{fullpoly}, in order to define a polyanalytic Bargmann transform. In this case, however, one needs to work with the polyanalytic Fock space (see Definition~\ref{paraFock}) and use a vector-valued input $ \mathbf{f} = (f_1, \ldots, f_N) $ for the Weyl polyanalytic operator.
\end{remark}

\subsection{Applications: Christoffel-Darboux formula and Mehler's kernel}
In this part of the section we derive a first consequence of the previous section using the well-known Christoffel–Darboux formula; see~\cite{AS}. The building blocks for this construction are the integral functions $ I_{a,m}(z)$ (see~\eqref{binte}), which arise from the composition of the true Bargmann transform, its inverse, and the translation operator, see Lemma \ref{L1}.

\begin{definition}
Let $N\in \mathbb N_0$ and $a\in\mathbb R$. Define the partial sum
\begin{equation}
\label{SS}
S_{a,N}(z,w)= \sum_{m=0}^{N}\frac{1}{2^m m!} I_{a,m}(z,w),
\end{equation}
for every $z,w\in \mathbb C.$
\end{definition}

The above summation has the following closed expression.

\begin{proposition}\label{P1}
For every $z,w\in \mathbb C$ and $a \in \mathbb{R}$ we have

$$\displaystyle S_{a,N}(z,w)=e^{-\frac{a^2}{4}+\frac{za}{\sqrt{2}}}  K_{N+1}\left(z-\frac{a}{\sqrt{2}},w \right).$$

\end{proposition}
\begin{proof}
Using the identity for generalized Laguerre polynomials, (see \eqref{l1}), and
by Corollary \ref{C1} we have
\begin{align*}
\displaystyle \displaystyle S_{a,N}(z,w)&= \sum_{m=0}^{N}\frac{1}{2^m m!} I_{a,m}(z,w)\\
&=e^{-\frac{a^2}{4}+\frac{za}{\sqrt{2}}} \sum_{m=0}^{N} \mathsf{K}_{m+1,T}\left(z-\frac{a}{\sqrt{2}},w \right)\\
&=e^{-\frac{a^2}{4}+\frac{za}{\sqrt{2}}} \sum_{m=0}^{N}e^{(z-\frac{a}{\sqrt{2}})\overline{w}} L^0_m\left(\left|z-\frac{a}{\sqrt{2}}-w \right|^2\right)\\
&=e^{-\frac{a^2}{4}+\frac{za}{\sqrt{2}}} e^{(z-\frac{a}{\sqrt{2}})\overline{w}}    L^{1}_{N}\left(\left|z-\frac{a}{\sqrt{2}}-w \right|^2\right) \\
&=e^{-\frac{a^2}{4}+\frac{za}{\sqrt{2}}}  K_{N+1}\left(z-\frac{a}{\sqrt{2}},w \right).
\end{align*}

\end{proof}

\begin{remark}
Define 
\begin{equation}
\label{serie1}
M_N(u,v)=\sum_{m=0}^{N}\frac{1}{2^m m!} H_m(u)H_m(v), \quad u,v\in \mathbb R.
\end{equation}

Then the Christoffel-Darboux formula for Hermite polynomials, see \cite{AS}, gives
$$M_N(u,v)=\frac{1}{N!2^{N+1}} \frac{H_N(v)H_{N+1}(u)-H_N(u)H_{N+1}(v)}{u-v}.$$

\end{remark}

By using the series~\eqref{serie1}, we can provide an integral representation of the series~\eqref{SS}.

\begin{proposition}\label{P2}
Let $N\in \mathbb N_0$ and $a\in\mathbb R$. Then for all $z,w\in \mathbb C$, we have

$$ \displaystyle S_{a,N}(z,w)= \int_\mathbb{R} M_N\left(\frac{z+\overline{z}}{\sqrt{2}}-x, \frac{w+\overline{w}}{\sqrt{2}}-x+a\right)  A_z(x) \overline{A_w(x-a)} dx,$$
where $A_z(x)$ is the Segal-Bargmann kernel introduced in \eqref{SBkernel}.
\end{proposition}
\begin{proof}
By formula \eqref{binte} we have
\begin{align*}
 \displaystyle S_{a,N}(z,w)&= \sum_{m=0}^{N}\frac{1}{2^m m!} I_{a,m}(z,w)\\
&=  \int_\mathbb{R} \left( \sum_{m=0}^{N}\frac{1}{2^m m!} H_m\left(\frac{z+\overline{z}}{\sqrt{2}}-x\right) H_m\left(\frac{w+\overline{w}}{\sqrt{2}}-x+a\right) \right)A_z(x) \overline{A_w(x-a)} dx \\
&= \int_\mathbb{R} D_{N,a}(z,w,x) A_z(x) \overline{A_w(x-a)} dx,
\end{align*}
where we define
$$ D_{N,a}(z,w,x):=\sum_{m=0}^{N}\frac{1}{2^m m!} H_m\left(\frac{z+\overline{z}}{\sqrt{2}}-x\right) H_m\left(\frac{w+\overline{w}}{\sqrt{2}}-x+a\right), \quad z,w\in \mathbb C, \quad x,a\in \mathbb R.$$
Finally, observe that

\begin{equation}
D_{N,a}(z,w,x)=M_N\left(\frac{z+\overline{z}}{\sqrt{2}}-x, \frac{w+\overline{w}}{\sqrt{2}}-x+a\right)
\end{equation}
\end{proof}

\begin{corollary}
Let $N\in \mathbb N_0$ and $a\in\mathbb R$. Then the polyanalytic Fock kernel $K_{N+1}$ admits the following integral representation
$$K_{N+1}(z,w)=\frac{1}{\sqrt{\pi}} e^{z \overline{w}} \int_\mathbb{R}  M_N\left(\frac{z+\overline{z}}{\sqrt{2}}-x, \frac{w+\overline{w}}{\sqrt{2}}-x\right) e^{-\left(x-\frac{z+\overline{w}}{\sqrt{2}}\right)^2}dx,\quad z,w\in \mathbb C.$$
\end{corollary}
\begin{proof}
By setting $a=0$ in Proposition \ref{P2} we have
\begin{align*}
\displaystyle  S_{0,N}(z,w)&= \frac{1}{\sqrt{\pi}} e^{-\frac{1}{2}(z^2+\overline{w}^2)}  \int_\mathbb{R} M_N\left(\frac{z+\overline{z}}{\sqrt{2}}-x, \frac{w+\overline{w}}{\sqrt{2}}-x\right)e^{-x^2+\sqrt{2}x(z+\overline{w})}  dx\\
&= \frac{1}{\sqrt{\pi}} e^{-\frac{1}{2}(z^2+\overline{w}^2)+\frac{(z+\overline{w})^2}{2}}  \int_\mathbb{R} M_N\left(\frac{z+\overline{z}}{\sqrt{2}}-x, \frac{w+\overline{w}}{\sqrt{2}}-x\right) e^{-\left(x-\frac{z+\overline{w}}{\sqrt{2}}\right)^2}dx\\
&=\frac{1}{\sqrt{\pi}} e^{z\overline{w}} \int_\mathbb{R} M_N\left(\frac{z+\overline{z}}{\sqrt{2}}-x, \frac{w+\overline{w}}{\sqrt{2}}-x\right) e^{-\left(x-\frac{z+\overline{w}}{\sqrt{2}}\right)^2}dx.
\end{align*}
On the other hand, Proposition \ref{P1} gives

$$K_{N+1}(z,w)=S_{0,N}(z,w).$$
Combining both results, we get the integral representation of the kernel $K_{N+1}$ in terms of the function $M_N(u,v)$.
\end{proof}
\begin{remark}
Using the expression for the polyanalytic Fock kernel $K_{N+1}$ in terms of generalized Laguerre polynomials, we obtain the identity
$$\displaystyle L^1_N(|z-w|^2)=\frac{1}{\sqrt{\pi}} \int_\mathbb{R} M_N\left(\frac{z+\overline{z}}{\sqrt{2}}-x, \frac{w+\overline{w}}{\sqrt{2}}-x\right) e^{-\left(x-\frac{z+\overline{w}}{\sqrt{2}}\right)^2}dx,\quad z,w\in \mathbb C.$$
Note that when $z=w=0$, we get

$$\sqrt{\pi}(N+1)= \int_\mathbb{R} M_N\left(x,x\right) e^{-x^2}dx,$$
where $(h_k)_{k\in\mathbb N_0}$ denotes the sequence of Hermite functions.
Therefore, the function $$\phi_N(x):=\frac{M_N(x,x)e^{-x^2}}{\sqrt{\pi}(N+1)}$$ is a probability density function, since it satisfies 
$$\displaystyle \phi_N(x)\geq 0, \quad \text{ and }\quad \int_\mathbb R \phi_N(x)dx=1.$$
\end{remark}

Now we present a second application involving the Mehler kernel. First, we introduce the following:
\begin{definition}
Let $N\in\mathbb N$ and $|\rho|<1$. Define the sequence of function below
\begin{equation}
\label{ss1}
S_{\rho, a,N}(z,w)= \sum_{m=0}^{N}\frac{\rho^m}{2^m m!} I_{a,m}(z,w),\quad \text{ for all } z,w\in \mathbb C.
\end{equation}
\end{definition}

We want to compute the limit of the series \eqref{ss1}.
\begin{theorem}\label{AppMehler}
Let $z,w\in\mathbb C$ and $|\rho|<1$. Then 
\begin{equation}
\label{limit}
\displaystyle \lim_{N\to \infty} S_{\rho,a,N}(z,w):=\sum_{m=0}^{\infty}\frac{\rho^m}{2^m m!} I_{a,m}(z,w)=\frac{1}{1-\rho}e^{-\frac{a^2}{4}+\frac{za}{\sqrt{2}}} e^{(z-\frac{a}{\sqrt{2}})\overline{w}} e^{-\frac{\rho}{1-\rho}\left|z-\frac{a}{\sqrt{2}}-w \right|^2}.
\end{equation}
\end{theorem}
\begin{proof}
We begin with 
\begin{align*}
\displaystyle S_{\rho,a,N}(z,w)&= \sum_{m=0}^{N}\frac{\rho^m}{2^m m!} I_{a,m}(z,w)\\
&=e^{-\frac{a^2}{4}+\frac{za}{\sqrt{2}}} \sum_{m=0}^{N} \rho^m \mathsf{K}_{m,T}\left(z-\frac{a}{\sqrt{2}},w \right)\\
&=e^{-\frac{a^2}{4}+\frac{za}{\sqrt{2}}} e^{(z-\frac{a}{\sqrt{2}})\overline{w}} \sum_{m=0}^{N} \rho^m L^0_m\left(\left|z-\frac{a}{\sqrt{2}}-w \right|^2\right).
\end{align*}

We now recall the classical Laguerre expansion (see \cite[Example 2, pp 89]{Lebedev1972}): 

\begin{equation}\label{ax}
e^{-\gamma x}=(\gamma+1)^{-(\alpha+1)}\sum_{n=0}^{\infty}\left(\frac{\gamma}{\gamma+1}\right)^{n}L_n^\alpha(x),\quad x\geq 0, \quad \forall \gamma>0.
\end{equation}

In particular, for  $\displaystyle 0<\rho=\frac{\gamma}{\gamma+1}<1$ and $\alpha=0,$ we have

$$e^{-\frac{\rho}{1-\rho}x}=(1-\rho)\sum_{n=0}^{\infty}\rho^n L^0_n(x), \quad x\geq 0.$$

Taking the limit as $N$ tends to $\infty$, we obtain 

\begin{align*}
\displaystyle \lim_{N\to \infty} S_{\rho,a,N}(z,w)&=e^{-\frac{a^2}{4}+\frac{za}{\sqrt{2}}} e^{(z-\frac{a}{\sqrt{2}})\overline{w}} \sum_{m=0}^{\infty} \rho^m L^0_m\left(\left|z-\frac{a}{\sqrt{2}}-w \right|^2\right)\\
&=\frac{1}{1-\rho}e^{-\frac{a^2}{4}+\frac{za}{\sqrt{2}}} e^{(z-\frac{a}{\sqrt{2}})\overline{w}} e^{-\frac{\rho}{1-\rho}\left|z-\frac{a}{\sqrt{2}}-w \right|^2},
\end{align*}
for every $ z,w\in \mathbb{C}.$ 
\end{proof}

\begin{remark}
The above result also follows by using the integral expression of the building blocks $I_{a,m}$  and the Mehler kernel, which is defined by  $$E_{\rho}(u,v):=\lim_{N\to \infty}E_{\rho,N}(u,v)=\sum_{m=0}^{\infty}\frac{\rho^m}{2^m m!} H_m(u)H_m(v), \quad u,v\in \mathbb R,$$
where $E_{\rho,N}(u,v)$ denotes the truncated Mehler kernel  (i.e., the partial sum up to degree $N$). Precsiely, by \eqref{binte} we have

\begin{align*}
	\displaystyle \displaystyle S_{\rho,a,N}(z,w)&= \sum_{m=0}^{N}\frac{\rho^m}{2^m m!} I_{a,m}(z,w)\\
	&=  \int_\mathbb{R} \left( \sum_{m=0}^{N}\frac{\rho^m}{2^m m!} H_m\left(\frac{z+\overline{z}}{\sqrt{2}}-x\right) H_m\left(\frac{w+\overline{w}}{\sqrt{2}}-x+a\right) \right)A_z(x) \overline{A_w(x-a)} dx \\
	&= \int_\mathbb{R} E_{\rho, N}\left(\frac{z+\overline{z}}{\sqrt{2}}-x, \frac{w+\overline{w}}{\sqrt{2}}-x+a\right) A_z(x) \overline{A_w(x-a)} dx,
\end{align*}
As a consequence of Cramér's inequality for Hermite polynomials, see \cite{H},
p. 435, we can apply the dominated convergence theorem and get
\begin{align}
	\nonumber
 \lim_{N\to \infty} S_{\rho,a,N}(z,w)&=\int_\mathbb{R} E_{\rho}\left(\frac{z+\overline{z}}{\sqrt{2}}-x, \frac{w+\overline{w}}{\sqrt{2}}-x+a\right) A_z(x) \overline{A_w(x-a)} dx\\
	\label{newk}
	&=\frac{1}{1-\rho}e^{-\frac{a^2}{4}+\frac{za}{\sqrt{2}}} e^{(z-\frac{a}{\sqrt{2}})\overline{w}} e^{-\frac{\rho}{1-\rho}\left|z-\frac{a}{\sqrt{2}}-w \right|^2}.
\end{align}

\end{remark}
For $a = 0 $, formula~\eqref{limit} simplifies to the product of the classical Fock kernel and a notable exponential function.
\begin{corollary}
Let $z,w\in\mathbb C$ and $|\rho|<1$. 
\begin{equation}
\displaystyle \lim_{N\to \infty} S_{\rho,0,N}(z,w)=K(z,w)K_\rho(z,w)=:S_\rho(z,w),
\end{equation}
where $K(z,w)=e^{z\overline{w}}$ is the classical Fock kernel, and
\begin{equation}
\label{kknew}
K_\rho(z,w)=\frac{1}{1-\rho}e^{-\frac{\rho}{1-\rho}\left|z-w \right|^2}.
\end{equation}
\end{corollary}
\begin{proof}
This follows directly from Theorem \ref{AppMehler} by setting $a=0$ as
\begin{align*}
\displaystyle \lim_{N\to \infty} S_{\rho,0,N}(z,w)&:=\sum_{m=0}^{\infty}\frac{\rho^m}{2^m m!} I_{0,m}(z,w)\\
&=e^{z\overline{w}} \frac{1}{1-\rho} e^{-\frac{\rho}{1-\rho}\left|z-w \right|^2}\\ 
&=K(z,w)K_\rho(z,w).
\end{align*}
\end{proof}

\begin{remark}
The function $K_\rho(z,w)$ in \eqref{kknew} is a reproducing kernel of a specific Hilbert space, see \cite[Remark 5.2]{ACDSS1}.
\end{remark}

\begin{corollary}
Consider the family of reproducing kernels defined by $$G_{\rho,N}(z,w):=e^{-\frac{z^2}{2}}S_{\rho,0,N}(z,w)e^{-\frac{\overline{w}^2}{2}},\quad \text{ for all } 0<\rho<1, \quad z,w\in \mathbb C.$$ The kernels $(G_\rho(z,w))_{0<\rho<1}$ converge pointwise to the Gaussian RBF kernel $K_{RBF} (z,w) K_{\rho}(z,w)$ as $N$ goes to $\infty$. 
\end{corollary}

\begin{remark}
When $a = 0$ in~\eqref{newk}, we also obtain the identity

\begin{equation}
\displaystyle \frac{1}{\sqrt{\pi}} \int_\mathbb{R} E_{\rho}\left(\frac{z+\overline{z}}{\sqrt{2}}-x, \frac{w+\overline{w}}{\sqrt{2}}-x\right) e^{\sqrt{2}x(z+\overline{w})} e^{-x^2}dx=\frac{1}{1-\rho}e^{\frac{(z+\overline{w})^2}{2}+\frac{\rho}{\rho-1}\left|z-w \right|^2} 
\end{equation}

In particular, for $z=w=0$,  we get
$$\displaystyle \frac{(1-\rho)}{\sqrt{\pi}} \int_\mathbb{R} E_{\rho}(x,x) e^{-x^2}dx=1.$$
Therefore, the function $$\phi_\rho(x):=\frac{(1-\rho)E_\rho(x,x)e^{-x^2}}{\sqrt{\pi}}, $$ is a probablity density function, since it satisfies 
$$\displaystyle \phi_\rho(x)\geq 0, \quad \text{ and }\quad \int_\mathbb R \phi_\rho (x)dx=1.$$
\end{remark}
The examples above motivate the introduction of a new integral transform in two complex variables involving the Mehler kernel.
\begin{definition}
For any $\phi\in L^2(\mathbb R) $, define the integral transform 
$$ (S\phi)(z,w):= \frac{1}{\sqrt{\pi}} \int_\mathbb{R}  e^{-\frac{1}{2}(z^2+\overline{w}^2)+\sqrt{2}x(z+\overline{w})} E_{\rho}\left(\frac{z+\overline{z}}{\sqrt{2}}-x, \frac{w+\overline{w}}{\sqrt{2}}-x\right) \phi(x)dx, \quad z,w\in \mathbb C.
$$
\end{definition}
\begin{remark}
Observe that  $$ (S\phi)(0,0):= \frac{1}{\sqrt{\pi}} \int_\mathbb{R} E_{\rho}(x,x)\phi(x)dx.$$
It may be of interest to investigate the transforms $$ (S\phi)(z,\overline{z}), \quad (S\phi)(z,z), \quad (S\phi)(z,0)$$ in future research.
\end{remark}

\subsection{The Weyl operator on a polyanalytic Gaussian RBF space}

In this part of the section we are going to use all the tools developed in this paper to establish a Weyl operator on the true polyanalytic RBF space.

\begin{definition}
Let $\ell \in \mathbb{N}$ and $\beta \in \mathbb{C}$. The polyanalytic RBF-Weyl operator on the space $ \mathcal{H}_{T, \sqrt{2}}^{(\ell)}(\mathbb{C})$ is defined such that the following diagram is commutative

\begin{equation}
\label{diag}
\xymatrix{
	\mathcal{H}^{(\ell)}_{T, \sqrt{2}}(\mathbb{C})\ar[r]_{W_{RBF,\ell}^\beta} \ar[d]_{\mathbf{M}_2 } & \mathcal{H}^{(\ell)}_{T, \sqrt{2}}(\mathbb{C})  \\ \mathcal{F}_T^{(\ell)}(\mathbb{C}) \ar[r]_{W_\ell^\beta} & \mathcal{F}_T^{(\ell)}(\mathbb{C}) \ar[u]_{\mathbf{M}^{-1}_2}
}
\end{equation}
Therefore we define the polyanalytic RBF-Wyl operator as
$$ W_{RBF,\ell}^\beta=\mathbf{M}^{-1}_2 \circ W_\ell^\beta \circ \mathbf{M}_2,$$
where $W_\ell^\beta$ is the polyanalytic Weyl operator, see Definition \ref{threeW}, and $\mathbf{M}_2$ is the operator defined in \eqref{iso3}.
\end{definition}
\begin{remark}
In the above definition, we have considered the true polyanalytic RBF space with $\gamma = \sqrt{2}$ in order to establish a connection with the true polyanalytic Fock space with parameter $\alpha = 1 $, which is the setting where the polyanalytic Weyl operator is constructed.
\end{remark}
The polyanalytic RBF Weyl operator satisfies properties analogous to those established for the Weyl operator in the classical RBF setting; see~\cite[Theorem~5.3 and Proposition~5.4]{ACDSS}. For this reason, we omit the proof and simply state the  following result.

\begin{theorem}
Let $\ell \in \mathbb{N}$ and $\beta \in \mathbb{C}$. The polyanalytic RBF-Weyl operator is an isometric operator from the true polyanalytic Fock space $\mathcal{H}_{T, \sqrt{2}}^{(\ell)}(\mathbb{C})$ onto itself. The adjoint and inverse of the polyanalytic RBF-Weyl operator are given by
$$ (W_{RBF, \ell}^\beta)^*=(W_{RBF, \ell}^\beta)^{-1}=W_{RBF, \ell}^{-\beta}.$$
Moreover, if we consider $\beta_1 \in \mathbb{C}$ we have
\begin{equation}
\label{weyl}
W_{RBF, \ell}^{\beta}W_{RBF, \ell}^{\beta_1}= e^{-\frac{2i}{\gamma^2} \hbox{Im}(\beta\overline{\beta_1})} W_{RBF, \ell}^{\beta+\beta_1}.
\end{equation}
\end{theorem}

\begin{remark}
If we consider $\beta$, $\beta_1 \in \mathbb{R}$ we have that the polyanalytic RBF operator satisfy the semi-group property
$$ W_{RBF, \ell}^{\beta}W_{RBF, \ell}^{\beta_1}= W_{RBF, \ell}^{\beta+\beta_1}.$$
\end{remark}

In the following result we provide an explicit expression of the polyanalytic RBF-Weyl operator.

\begin{theorem}
Let $\ell \in \mathbb{N}$. Let $\beta \in \mathbb{C}$ and $f \in \mathcal{H}_{T, \sqrt{2}}^{(\ell)}(\mathbb{C})$. Then we have
\begin{equation}
	\label{RBFpoly}
W_{RBF, \ell}^{\beta}=(\ell-1)! 2^{\ell-1}  e^{\frac{z}{\sqrt{2}}(\beta-\bar{\beta})}e^{\frac{1}{4}(\bar{\beta}^2-|\beta|^2)}f\left(z-\frac{\bar{\beta}}{\sqrt{2}}\right), \quad z \in \mathbb{C}.
\end{equation}
\end{theorem}
\begin{proof}
By Theorem \ref{iso} and Theorem \ref{TrueWeylR2} we have
\begin{eqnarray*}
	W_{RBF, \ell}^{\beta}[f](z)&=&\left(\mathbf{M}^{-1}_2 \circ W_{\ell}^\beta\right) \left(e^{\frac{z^2}{2}} f\right)(z)\\
	&=& (\ell-1)! 2^{\ell-1} e^{- \frac{z^2}{2}} W_{\ell}^\beta \left(e^{\frac{z^2}{2}} f\right)(z)\\
	&=& (\ell-1)! 2^{\ell-1} e^{- \frac{z^2}{2}} e^{\frac{z \beta}{\sqrt{2}}} e^{- \frac{|\beta|^2}{4}} e^{\frac{1}{2} \left(z-\frac{\overline{\beta}}{\sqrt{2}}\right)^2}f\left(z-\frac{\bar{\beta}}{\sqrt{2}}\right)\\
	&=& (\ell-1)! 2^{\ell-1}  e^{\frac{z}{\sqrt{2}}(\beta-\bar{\beta})}e^{\frac{1}{4}(\bar{\beta}^2-|\beta|^2)}f\left(z-\frac{\bar{\beta}}{\sqrt{2}}\right).
\end{eqnarray*}
\end{proof}

\begin{remark}
If we consider $ \beta \in \mathbb{R}$ in~\eqref{RBFpoly}, we obtain that the polyanalytic RBF Weyl operator reduces to the standard translation operator. Namely, for $f \in \mathcal{H}_{T, \sqrt{2}}^{(\ell)}(\mathbb{C})$, we have

$$ W_{RBF, \ell}^{\beta}[f](z)=(\ell-1)! 2^{\ell-1}  f\left(z-\frac{\bar{\beta}}{\sqrt{2}}\right).$$
\end{remark}

\begin{remark}
If one considers the polyanalytic Gaussian RBF space (see Definition~\ref{polyRBF}) and the polyanalytic Fock space (see Definition~\ref{paraFock}) instead of the true polyanalytic RBF space and the true polyanalytic Fock space in \eqref{diag}, one obtains a vector-valued version of the polyanalytic RBF operator; see Remark~\ref{vectorial}.

\end{remark}

\section{Appendix}
\appendix
This appendix is divided into three parts.  In the first part, we provide a list of formulas and results of the paper in their parametric versions. In the second part, we present a proof of a generating formula for complex Hermite polynomials, based on reproducing kernel techniques. In the final part, we give an algebraic proof of the convolution formula for Hermite functions. \\ \\
The following table summarizes the main function spaces encountered in the paper: 

\begin{table}[h!]
	\caption{List of function spaces and their reproducing kernels.}
    \centering
    \begin{tabular}{|c|c|c|}
        \hline
        \textbf{Space} & \textbf{Notation} & \textbf{Reproducing Kernel} \\
        \hline
        Fock space (Definition~\ref{Fockdef}) 
        & $\mathcal{F}_\alpha(\mathbb{C})$ 
        & $K_\alpha(z,w)$ (Formula~\eqref{Kf}) \\
         \hline
        Fock space ($\alpha=1$)
        & $\mathcal{F}(\mathbb{C})$ 
        & $K(z,w)$ \\
        \hline
        Fock space ($\alpha=\frac{1}{2}$)
        & $\mathcal{F}_{1/2}(\mathbb{C})$ 
        & $K_{1/2}(z,w)$ \\
        \hline
        Gaussian RBF space (Definition~\ref{RBF}) 
        & $\mathcal{H}_\gamma^{RBF}(\mathbb{C})$  
        & $K^\gamma_{RBF}(z,w)$ (Formula~\eqref{rbfkernel}) \\
        \hline
        Gaussian RBF space ($\gamma=2$) 
        & $\mathcal{H}_2(\mathbb{C})$  
        & $K_{RBF}(z,w)$ \\
        \hline
        Polyanalytic Fock space (Definition~\ref{paraFock}) 
        & $\mathcal{F}_{\alpha,N}(\mathbb{C})$ 
        & $K_N^\alpha(z,w)$ (Formula~\eqref{Kn}) \\
         \hline
        Polyanalytic Fock space ($\alpha=1$) 
        & $\mathcal{F}_{N}(\mathbb{C})$ 
        & $K_N(z,w)$ \\
        \hline
        True polyanalytic Fock space (Definition~\ref{FT}) 
        & $\mathcal{F}_{T,\alpha}^{(\ell)}(\mathbb{C})$ 
        & $\mathsf{K}_{\ell,T}^\alpha(z,w)$ (Formula~\eqref{tpoly}) \\
         \hline
        True polyanalytic Fock space ($\alpha=1$) 
        & $\mathcal{F}_{T}^{(\ell)}(\mathbb{C})$ 
        & $\mathsf{K}_{\ell,T}(z,w)$ \\

        \hline
        True polyanalytic RBF space (Definition~\ref{polyRBF}) 
        & $\mathcal{H}_{T, \gamma}^{(\ell)}(\mathbb{C})$ 
        & $K_{T, \ell}(z,w)$ (Formula~\eqref{KtrueRBF}) \\
        
         \hline
        True polyanalytic RBF space ($\gamma=\sqrt{2}$) 
        & $ \mathcal{H}^{(\ell)}_{T, \sqrt{2}}(\mathbb{C})$ 
        & Does not appear\\

        \hline
        Polyanalytic Gaussian RBF space (Definition~\ref{trurbf}) 
        & $\mathcal{H}_{N, \gamma}(\mathbb{C})$ 
        & $K_{RBF,N}(z,w)$ (Formula~\eqref{KpolyRBF}) \\
        \hline
       
    \end{tabular}
\end{table}

\section{Results with a generic parameter}
Most of the results of the paper are presented for specific parameters: $\alpha = \frac{1}{2}$ (corresponding to the Fock space) and $\gamma=2$ for the RBF space in Sections 3 and 4, and $\alpha=1$ and $\gamma=\sqrt{2}$ in Section 6. In this part of the appendix, we present the main results for a general parameter $\gamma > 0$. This can be considered as a straightforward parametric extension that may be of independent interest. Indeed, we have already discussed a parametric version of the true polyanalytic RBF kernel and polyanalytic Gaussian RBF in Remark \ref{An}.

An orthonormal basis of the related true polyanalytic (parametric) RBF space $\mathcal{H}_{T, \gamma}^{(\ell)}(\mathbb{C})$ is given by

$$ e_{\ell,m}^\gamma(z, \bar{z}):= \frac{\gamma^{\ell+m}}{\sqrt{ (\ell-1)! m!} 2^{\frac{\ell+m}{2}}}H_{\ell-1,m}^{\frac{2}{\gamma^2}}(z, \bar{z}) e^{-\frac{z^2}{\gamma^2}},$$ 
where $H_{\ell-1,m}^{\frac{2}{\gamma^2}}$ are the complex Hermite, see \eqref{hermcom}.

The generalization of the operator $ \mathcal{A}_{z, \bar{z}}$ in terms of the parameter $\gamma>0$ is given by
$$
\mathcal{A}_{z, \bar{z}, \gamma}:=\left(- \frac{d}{dz}- \frac{2}{\gamma^2}M_z+\frac{2}{\gamma^2}M_{\bar{z}}\right).
$$

The application of $(\ell-1)$-th power of the operator $ \mathcal{A}_{z, \bar{z}, \gamma}$ to the basis of the RBF space $\mathcal{H}_{\gamma}^{RBF}(\mathbb{C})$ gives the following

$$
	\mathcal{A}_{z, \bar{z}, \gamma}^{\ell-1} \left( z^m e^{- \frac{z^2}{\gamma^2}}\right)= \frac{\gamma^{2m}}{2^m}H_{\ell-1,m}^{\frac{2}{\gamma^2}}(z, \bar{z}) e^{-\frac{z^2}{\gamma^2}}, \quad \ell \in \mathbb{N}
$$

The connection between the reproducing kernel of the RBF space and the true polyanalytic RBF space can be given through the  application of the operator $\mathcal{A}_{z, \bar{z}, \gamma}^{\ell-1}$. Thus we can write

	$$ \mathcal{A}_{z, \bar{z}, \gamma}^{\ell-1}  \overline{\mathcal{A}_{w, \bar{w}, \gamma}^{\ell-1}} \left( e^{- \frac{(z-\bar{w})^2}{\gamma^2}} \right)=\frac{2^{\ell-1}(\ell-1)!}{\gamma^{2(\ell-1)}}L_{\ell-1}^{0}\left(\frac{2|z-w|^2}{\gamma^2}\right)e^{- \frac{(z-\bar{w})^2}{\gamma^2}}, \qquad z,w \in \mathbb{C}.$$
The adjoint of the operator $\mathcal{A}_{z, \bar{z}, \gamma}$ is still $ \mathcal{A}_{z, \bar{z}, \gamma}= \frac{d}{d\bar{z}}$. So we can define the following operator

$$ \widetilde{\Box}_\gamma:=\mathcal{A}_{z, \bar{z}, \gamma} \mathcal{A}_{z, \bar{z}, \gamma}^{*}=- \frac{d^2}{dz d \bar{z}}-\frac{2z}{\gamma^2}\frac{d}{d \bar{z}}+\frac{2 \bar{z}}{\gamma^2} \frac{d}{d \bar{z}}, \quad  \, z \in \mathbb{C}.$$

One can prove that the functions
$$ \widetilde{H}_{\ell,m}^\gamma(z, \bar{z})= H_{\ell-1,m}^{\frac{2}{\gamma^2}}(z, \bar{z})e^{- \frac{z^2}{\gamma^2}},$$
are eigenfunctions of the the operator $\widetilde{\Box}_\gamma$, i.e.:
$$ \widetilde{\Box}_\gamma [\widetilde{H}_{\ell,m}^\gamma(z, \bar{z})]= \frac{2(\ell-1)}{\gamma^2} \widetilde{H}_{\ell,m}^{\frac{2}{\gamma^2}}(z, \bar{z}).$$	
The operator $\widetilde{\Box}_\gamma$ can be expressed in terms of the Schrödinger-type operator
$$
\mathcal{L}_\gamma := -\frac{1}{4} \Delta - \frac{2i y}{\gamma^2} \frac{\partial}{\partial x} + \frac{4y^2}{\gamma^4}, \qquad x, y \in \mathbb{R},
$$
via the RBF ground state transform (with parameter) given by
$$
G_\gamma \colon L^2(\mathbb{R}^2, dx\,dy) \to L^2(\mathbb{R}^2, g^2\,dx\,dy), \qquad f \mapsto \frac{f}{g},
$$
where the ground state function is given by
$$
g(y) := \left( \frac{2}{\pi \gamma^2} \right)^{1/2} e^{-\frac{2 y^2}{\gamma^2}}.
$$
Hence the relation between the operators is
$$
\widetilde{\Box}_\gamma = G_\gamma \circ \left( \mathcal{L}_\gamma - \frac{1}{\gamma^2} \right) \circ G_\gamma^{-1}.
$$
The above relation helps us to prove that the eigenvalues of the operator $ \mathcal{L}_\gamma$ are $\frac{2}{\gamma^2}\left(\ell-1\right)$, and corresponding eigenfunctions are $g(y) \widetilde{H}_{\ell,m}^\gamma(z,\bar{z})$. Moreover it is also possible to prove that the true polyanalytic RBF space $\mathcal{H}_{T, \gamma}^{(\ell)}(\mathbb{C})$ with a parameter $\gamma>0$ coincides with the RBF Landau level given by
 $$ \mathbf{H}_\gamma^n(\mathbb{C}):= \biggl \{f\in L^{2}\left(\mathbb{C}, e^{\frac{(z-\bar{z})^2}{\gamma^2}}\right), \qquad 	\widetilde{\Box}_\gamma(f(z))=\frac{2(\ell-1)}{\gamma^2}  f(z) \biggl\}.$$

\section{A generating formula for Itô-Hermite polynomials}\label{App2}
Although the following Lemma may be known, we provide a proof here, as we were not able to find an explicit proof in the literature.
\begin{lemma}\label{SumH}
	Let $m\in \mathbb{N}_0$ be fixed and $\alpha>0$. Then, for all $z,w\in \mathbb{C}$, we have 
	\begin{equation}
		\frac{1}{\alpha^m} \sum_{n=0}^{\infty} \frac{H_{n,m}^\alpha(z,\overline{z}) H_{m,n}^\alpha(w,\overline{w})}{m!n! \alpha^n}=e^{\alpha z\overline{w}}L_{m}^0\left(\alpha|z-w|^2\right).
	\end{equation}
	Here $L^{\beta}_{k}$ denotes the generalized Laguerre polynomials defined in \eqref{Laguerre}.
	
\end{lemma}
\begin{proof}
	Consider the kernel function
	\begin{equation}
		K_m(z,w)=K_m^w(z):=\frac{1}{\alpha^m} \sum_{n=0}^{\infty} \frac{H_{n,m}^\alpha(z,\overline{z}) H_{m,n}^\alpha(w,\overline{w})}{m!n! \alpha^n},
	\end{equation}
	for all $z,w\in \mathbb{C}$. It is known that $ \lbrace{ e_{m,n}(z,\overline{z})=H_{n,m}^\alpha (z,\overline{z})\mid n\in \mathbb{N}_0 \rbrace}$ forms an orthonormal basis of the true (parametric) polyanalytic Fock space, see \cite{Asampling}. We will show that $K_m$ is the reproducing kernel of $\mathcal{F}^{(m)}_{T,\alpha}(\mathbb{C})$. Let $f\in\mathcal{F}^{(m)}_{T,\alpha}(\mathbb{C})$ be given by $$\displaystyle f(w)=\sum_{k=0}^{\infty}\beta_kH_{k,m}^\alpha(w,\overline{w}),\quad w\in \mathbb{C}, \quad \{\beta_k\}_{k \in \mathbb{N}_0} \subseteq \mathbb{C}$$ 
	Then, we compute the inner product of $f$ and $K_m$. By the orthogonality of the complex Hermite polynomials, see \eqref{ort}, we have
	\begin{align*}
		\displaystyle \langle f, K_m^w \rangle_{\mathcal{F}^{(m)}_{T,\alpha}(\mathbb{C})}&= \frac{\alpha}{\pi}\int_{\mathbb{C}} \overline{K_m(z,w)} f(z) e^{-\alpha|z|^2}d \lambda(z)\\
		&=\sum_{k=0}^{\infty} \beta_k \frac{\alpha}{\pi} \int_{\mathbb{C}}\overline{K_m(z,w)} H_{k,m}^\alpha(z,\overline{z})  e^{-\alpha|z|^2}d \lambda(z)\\
		&=\sum_{k,n=0}^{\infty} \frac{ \beta_k}{m!n! \alpha^{m+n}} H_{n,m}^\alpha(w,\overline{w}) \frac{\alpha}{\pi}  \int_{\mathbb{C}}\overline{H_{n,m}^\alpha(z,\overline{z})} H_{k,m}^\alpha(z,\overline{z})  e^{-\alpha |z|^2}d \lambda(z)\\
		&=\sum_{n=0}^{\infty}\beta_nH_{n,m}^\alpha(w,\overline{w})\\
		&=f(w).
	\end{align*}
	Therefore, $K_m$ reproduces all elements of $\mathcal{F}^{(m)}_{T, \alpha}(\mathbb{C})$. Since it is known that $$e^{\alpha z\overline{w}}L_{m}^0\left(\alpha |z-w|^2\right)$$ is the reproducing kernel of $\mathcal{F}^{(m)}_{T,\alpha}(\mathbb{C})$ (see \cite{AF}), uniqueness of the reproducing kernel implies that  
	$$K_m(z,w)=e^{\alpha z\overline{w}}L_{m}^0\left(\alpha|z-w|^2\right),$$
	for all $z,w\in \mathbb{C}$. 
\end{proof}

\section{Convolution of Hermite functions: an algebraic proof}

In this part of the section, we provide an alternative proof of the convolution $h_k \ast h_\ell$ of two Hermite functions, where $k \leq \ell$ to fix the ideas. This is done in several steps, where the computation will be reduced to $h_0 \ast h_0$. We start by recalling an alternative definition of the Hermite functions of degree $k$:

$$
h_k(x) := (a^+)^k \, e^{-x^2 / 2}
$$
with
$$
a^+ := x - \frac{d}{dx}, \qquad a^- := x + \frac{d}{dx}.
$$
In the following result we will make use of the following definition of Fourier transform ($\mathcal{F}^+$) and its inverse ($\mathcal{F}^{-}$):
$$ \mathcal{F}^{\pm}(f)=\frac{1}{\sqrt{2 \pi}} \int_{\mathbb{R}} f(x)e^{\pm i \xi x}dx, \qquad f \in \mathcal{S}(\mathbb{R})$$
\begin{lemma}
Let $f$, $g \in \mathcal{S}(\mathbb{R})$ then we have
\begin{eqnarray}
	\label{one}
	x \left( f \ast g\right) &=& (x\,f) \ast g + f \ast (x\, g) \\
	\label{second}
	\frac{d}{dx} \left( f \ast g\right) &=& \left(\frac{d}{dx} \,f\right) \ast g = f \ast \left(\frac{d}{dx} \, g\right) .
\end{eqnarray}
Moreover we have
\begin{equation}\label{eq:actiona}
(a^+ f) \ast g = x \left( f \ast g\right) - f \ast \left(a^-\, g\right)
\end{equation}
\end{lemma}
\begin{proof}
Using formulas \eqref{one} and \eqref{second}, we find that
\begin{align*}
	\left[ \left( x - \frac{d}{dx}\right) f \right] \ast g &= \left( x \left( f \ast g\right) - f \ast (x\, g)\right)  - f \ast \left(\frac{d}{dx} \, g\right) \\
	&= x \left( f \ast g\right) - f \ast \left(\left(x + \frac{d}{dx}\right) g\right). 
\end{align*}
This proves formula \eqref{eq:actiona}.
\end{proof}
\begin{lemma}
Let $k \in \mathbb{N}$. For $f$, $g \in \mathcal{S}(\mathbb{R})$ we have
\begin{equation}
\label{kk}
\left( \left(a^+\right)^k f \right) \ast g = \sum_{j=0}^k (-1)^j  \binom{k}{j} \, x^{k-j} \left( f \ast ((a^-)^j \, g) \right),
\end{equation}
\end{lemma}
\begin{proof}
We prove the result by induction on $k$. For $k=1$, it follows from \eqref{eq:actiona} that

$$
\left( a^+\, f \right) \ast g =  \sum_{j=0}^1 (-1)^j  \binom{k}{j} \, x^{1-j} \left( f \ast ((a^-)^j \, g) \right)= x(f \ast g)-f \ast(a^{-}g).
$$
Thus the result follows by \eqref{eq:actiona}. For $k \geq 1$, by using the inductive hypothesis and the Stifel identity we have
 
\begin{align*}
	\left( \left(a^+\right)^{k+1} f \right) & \ast g = \left( \left(a^+\right)^{k} (a^+ \, f) \right) \ast g \\
	&= \sum_{j=0}^k (-1)^j \binom{k}{j} \, x^{k-j} \left( (a^+\,f) \ast ((a^-)^j \, g) \right) \\
	&\stackrel{(\ref{eq:actiona})}{=} \sum_{j=0}^k (-1)^j \binom{k}{j} \, x^{k-j} \left( x \left( f \ast  ((a^-)^j \, g)\right) - f \ast \left((a^-)^{j+1} \, g\right) \right) \\
	&= \sum_{j=0}^k(-1)^j \binom{k}{j} \, x^{k-j+1} \left( f \ast  ((a^-)^j \, g)\right) - \sum_{j=1}^{k+1} (-1)^j \binom{k}{j-1} \, x^{k-(j-1)} \left(f \ast \left((a^-)^j \, g\right) \right) \\
	&= x^{k+1}\, (f \ast g) + \sum_{j=1}^k (-1)^j \left( \binom{k}{j}+ \binom{k}{j-1}\right) x^{k-j+1} \left(f \ast \left((a^-)^j \, g\right) \right) - (-1)^k \ f \ast \left((a^-)^{k+1} \, g\right) \\
	&= \sum_{j=0}^{k+1} (-1)^j \binom{k+1}{j} x^{k+1-j} \left( f \ast ((a^-)^j \, g) \right).
\end{align*}
\end{proof}

\begin{lemma}
Let $k \in \mathbb{N}$. Then the action of the operator $a^{-}$ on the Hermite function leads to
\begin{equation}
\label{action}
a^- h_k=2k h_{k-1}.
\end{equation} 
\end{lemma}
\begin{proof}
To determine the action of $a^{-}$ on $h_k$, we first need to understand how $a^-$ and $a^+$ interact. We thus see that
\begin{align*}
	a^-\,h_k &= a^- \left(a^+\right)^k \, e^{-x^2 \slash 2} = \left(a^+ \, a^- + 2\right) \left(a^+\right)^{k-1} \, e^{-x^2 \slash 2} \\
	&= a^+ \left(2\,(k-1)\right) \left(a^+\right)^{k-2} \, e^{-x^2 \slash 2} + 2 \left(a^+\right)^{k-1} \, e^{-x^2 \slash 2} \\
	&= 2\,k \left(a^+\right)^{k-1} \, e^{-x^2 \slash 2} = 2kh_{k-1}.
\end{align*}
\end{proof}

Now, we have all the tools to show how to compute the convolution of two Hermite functions

\begin{theorem}
Let $k$, $\ell \in \mathbb{N}$ such that $k \leqslant \ell$, then we have
$$ (h_k \ast h_\ell)(x)=\sqrt{\pi} \, e^{-\frac{x^2}{4}}  \sum_{j=0}^k (-1)^j \, \binom{k}{j} \, 2^j \, \frac{\ell!}{(\ell-j)!}\, x^{k+\ell-2j}, \quad \text{for all } x\in\mathbb{R}. 
$$
\end{theorem}
\begin{proof}
By the definition of the Hermite functions, \eqref{kk} and \eqref{action} we have

\begin{align*}
	(h_k \ast h_\ell)(x) &= \left( a^+\right)^k (h_0 \ast h_\ell)(x) \\
	&= \sum_{j=0}^k (-1)^j \, \binom{k}{j} \, x^{k-j} \left( h_0 \ast \left((a^-)^j \, h_\ell\right) \right)(x)  \\
	&= \sum_{j=0}^k (-1)^j \, \binom{k}{j} \, (2\ell) \, (2\ell-2) \, \ldots (2\ell- 2j+2) \, x^{k-j} \left( h_0 \ast h_{\ell-j} \right)(x) \\
	&= \sum_{j=0}^k (-1)^j \, \binom{k}{j} \, 2^j \, \frac{\ell!}{(\ell-j)!}\, x^{k-j} \left( h_0 \ast h_{\ell-j} \right) (x)
\end{align*}

We now still determine  $h_0 \ast h_k$. By \eqref{eq:actiona} we have
\begin{align*}
	h_0 \ast h_k &= h_0 \ast \left(a^+\right)^k \,h_0 = h_0 \ast a^+ \left(a^+\right)^{k-1} \,h_0 \\ 
	&= x \left( h_0 \ast \left(a^+\right)^{k-1} \,h_0\right) - (a^-\,h_0) \ast \left(a^+\right)^{k-1} \,h_0.
\end{align*}
Since $
a^- \, h_0 = \left(x+ \frac{d}{dx}\right) e^{-x^2 \slash 2} = 0
$
we get
\begin{equation}
\label{final2}
	h_0 \ast h_k= x \left( h_0 \ast \left(a^+\right)^{k-1} \,h_0\right) = \ldots = x^k \left( h_0 \ast h_0\right).
\end{equation}
Thus we have
$$
	(h_k \ast h_\ell)(x) = \sum_{j=0}^k (-1)^j \, \binom{k}{j} \, 2^j \, \frac{\ell!}{(\ell-j)!}\, x^{k+\ell-2j}  \left( h_0 \ast h_0 \right)(x).
$$
The proof ends by observing that $$(h_0 \ast h_0)(x)= \sqrt{\pi} \, e^{-\frac{x^2}{4}}.$$

\end{proof}

\section*{Acknowledgments}
The authors are grateful to Professor Sundaram Thangavelu for pointing out the references on special Hermite functions \cite{Strichartz1989, Thangavelu1993}. The research of Kamal Diki is supported by the Research Foundation–Flanders (FWO) under grant number 1268123N.

\end{document}